\documentclass[10pt,twocolumn,journal]{IEEEtran}

\usepackage{float}
\usepackage{enumitem}
\usepackage{graphicx}
\usepackage{amssymb,amsfonts,amsthm,thmtools,thm-restate}
\usepackage{grffile}
\usepackage{bigints}
\usepackage{tikz}
\usetikzlibrary{shapes.geometric, arrows, positioning, calc}
\usepackage{mathdots}
\usepackage{cite}
\usepackage{caption,subcaption,url,hyperref,cleveref}
\usepackage{xcolor}

\usepackage[compact]{titlesec}

\newcommand{\ubar}[1]{\text{\b{$#1$}}}

\newcommand{\mcl}[1]{\mathcal{ #1}}
\newcommand{\mbf}[1]{\mathbf{ #1}}
\newcommand{\norm}[1]{\left\Vert #1\right\Vert}

\newcommand{\ip}[2]{\left\langle{#1},{#2}\right\rangle}

\newcommand{\bmat}[1]{\begin{bmatrix} #1\end{bmatrix}}

\newcommand{\mat}[1]{\begin{matrix}#1\end{matrix}}

\newcommand{\R}{\mathbb{R}}
\newcommand{\C}{\mathbb{C}}

\newcommand{\N}{\mathbb{N}}

\newcommand{\myint}{\int_{a}^{b}}
\newcommand{\myinta}[1]{\int_{a}^{#1}}
\newcommand{\myintb}[1]{\int_{#1}^{b}}

\newcommand\vlines{\fboxsep=-\fboxrule\!\!\!\fbox{\strut}\!\!\!}

\newtheorem{theorem}{Theorem}
\newtheorem{illus}{Illustration}
\newtheorem{definition}[theorem]{Definition}
\newtheorem{lemma}[theorem]{Lemma}

\newtheorem{corollary}[theorem]{Corollary}

\crefname{EqnBlock}{Block}{Blocks}

\newcommand{\pie}{\scalebox{0.9}[1.2]{$\mathit{\Pi}$}}
\newcommand{\picomp}[8]{\ensuremath{\mbf P_\times^4\left(\left[\footnotesize\begin{array}{c|c}
#1&#2\\\hline #3 & \{#4\}
\end{array}\right],\left[\footnotesize\begin{array}{c|c}
#5&#6\\\hline #7 & \{#8\}
\end{array}\right]\right)}}
\newcommand{\fourpi}[4]{\ensuremath{\pie\left[\footnotesize\begin{array}{c|c}
#1&#2\\\hline #3 & \{#4\}
\end{array}\right]}}
\newcommand{\threepi}[1]{\pie_{\{#1\}}}
\newcommand{\fourpiFull}[6]{\ensuremath{\pie\left[\footnotesize\begin{array}{c|c}
#1&#2\\\hline #3 & \left\lbrace#4,#5,#6\right\rbrace
\end{array}\right]}}

\title{Extension of the Partial Integral Equation Representation to GPDE Input-Output Systems\thanks{Submitted to IEEE TAC.}}

\author{Sachin~Shivakumar,~Amritam~Das,~Siep~Weiland,~and~Matthew~Peet,~\IEEEmembership{Member,~IEEE}
\thanks{S. Shivakumar and M. Peet are with School for Engineering of Matter, Transport and Energy, Arizona State University (email: sshivak8@asu.edu, mpeet@asu.edu).}
\thanks{A. Das and S. Weiland are with Department of Electrical Engineering, Eindhoven University of Technology (email: Am.Das@tue.nl, S.Weiland@tue.nl).}
\thanks{This work was supported by National Science Foundation under Grants No. CMMI-1935453 and CNS-1739990.}
}

\setlist{nosep}
\allowdisplaybreaks
\titlespacing{\section}{0pt}{2ex}{1ex}
\titlespacing{\subsection}{0pt}{1ex}{0ex}
\titlespacing{\subsubsection}{0pt}{0.5ex}{0ex}
\renewcommand\thesubsubsection{\thesubsection \arabic{subsubsection}}
\renewcommand\thesubsubsectiondis{\thesubsectiondis \arabic{subsubsection}}

\setlist[enumerate]{wide=0pt, widest=99,labelwidth=1pt,labelsep=2pt,leftmargin=4pt}
\setlist[itemize]{wide=0pt, widest=99,labelwidth=1pt,labelsep=2pt,leftmargin=4pt}
\begin{document}

\maketitle

	\begin{abstract}
	
It has been shown that the existence of a Partial Integral Equation (PIE) representation of a Partial Differential Equation (PDE) simplifies many numerical aspects of analysis, simulation, and optimal control. However, the PIE representation has not previously been extended to many of the complex, higher-order PDEs such as may be encountered in speculative or data-based models. In this paper, we propose PIE representations for a large class of such PDE models, including higher-order derivatives,  boundary-valued inputs, and coupling with Ordinary Differential Equations. The main technical contribution which enables this extension is a generalization of Cauchy's rule for repeated integration. The process of conversion of a complex PDE model to a PIE is simplified through a PDE modeling interface in the open-source software PIETOOLS. Several numerical tests and illustrations are used to demonstrate the results.

	\end{abstract}
	

	\section{Introduction}\label{sec:Introduction}	

	While Partial Differential Equations (PDEs) have historically been used to model spatially distributed physical phenomena dating back to Newton and Leibniz, the critical role of boundary conditions (BCs) in defining a PDE model was not formally recognized until Dirichlet's time; See~\cite{CHENG2005268} for a historical survey of PDEs and BCs. Even with BCs, a PDE model remains incomplete without `continuity' requirements on the solution. This means that spatial derivatives and boundary values must be well-defined. The formal mathematical framework for imposing continuity constraints was only established in the mid-20th century by Sergei Sobolev, introducing Sobolev spaces and allowing for the use of generalized functions or distributions to define weak solutions.
 
	The PDE, BCs, and continuity constraints create a `PDE model'-- a system relying on three types of constraints, none of which is individually sufficient. Yet, when combined, the 3 constraints establish a well-defined mapping from an initial state to a unique solution. In the late 20th century, this mapping and its continuity properties were formalized and generalized under the concept of a $C_0$-semigroup where these constraints are embedded in the `domain of the infinitesimal generator'~\cite{Engel2000,zwart-curtain}. 
    Thus, a well-posed PDE model requires three constraints: a) the differential equation or `PDE', constraining spatiotemporal evolution within the domain; b) the continuity condition, ensuring sufficient regularity of the solution; and c) the BCs, which may limit values or other properties of the solutions as allowed by the regularity.

	\subsection{Computational Challenges of Using a PDE Model}
	Although the 3-constraint PDE model is a natural representation due to its historical context and clear physical interpretation, using such a model becomes inconvenient when considering computational methods for analyzing, controlling, and simulating spatially distributed phenomena. The most significant inconveniences are as follows:
	\begin{enumerate}
		\item \textbf{Non-Algebraic Structure} All computations fundamentally involve algebraic operations addition and multiplication. PDE models, however, have operators such as spatial differentials and Dirac operators (evaluating limit points) -- neither of which can be incorporated into a *-algebra of bounded linear operators on a Hilbert space~\cite{segal_1947}.
		\item \textbf{No Universality} Computational methods are traditionally centered on the `PDE' part of the `PDE model', and are designed for a fixed set of BCs and continuity constraints. As a result, there are no generic/universal algorithms for the analysis, control, and simulation of PDEs.
	\end{enumerate}

	To illustrate, consider the problem of simulating a simple transport equation $u_t=u_{s}$ for a given initial condition. Using a finite-difference approximation of $u_{s}=\frac{u(s_{i+1})-u(s_i)}{\Delta_s}$ yields an ODE representation $\dot x(t)=\frac{1}{\Delta_s}Ax(t)$, where $x_i=u(s_i)$, $\Delta_s=s_{i+1}-s_i$, and $A$ is a bi-diagonal matrix of $\pm 1$ entries. An ideal ODE representation of the transport equation would have a matrix with infinitely large coefficients. Such problems associated with discretization are avoided by constructing an explicit basis for the domain of the infinitesimal generator and projecting our solution onto this basis -- e.g., as in Galerkin projection. However, these bases must satisfy the continuity constraints and BCs present in the 3-constraint model. Thus, any change in the set of BCs and continuity constraints necessitates a change in the bases -- an inconvenience and an obstacle to the design of general/universal simulation tools.

	To illustrate challenges in computational analysis of a PDE model, consider the stability analysis problem of a heat equation $u_t=u_{ss}$ with zero BCs, e.g. $u(t,0)=u_s(t,1)=0$. We can use a Lyapunov function $V(u)=\int_0^1 u(s)^2 ds$ that is uniformly decreasing with time to prove the stability. The challenge, however, is to use computation to prove this fact. By parameterizing positive operators using positive matrices, computation-based methods can recognize that $V(u)=\ip{u}{u}_{L_2}$ and hence a positive form~\cite{peet_SICON2008}. However, the method must also verify that $\dot V(u(t))\le 0$ for all $u(t)\in W_2$ satisfying the PDE model. Unfortunately, if we differentiate $V(u(t))$ in time along solutions of the PDE model we obtain $\dot V(u(t))=2\ip{u(t)}{\partial^2_s u(t)}=2\int_0^1 u(t,s)u_{ss}(t,s) ds$. Because differentiation is not embedded in a $^*$-algebra, we cannot simply parameterize a cone of positive quadratic forms involving differential operators, e.g., $\ip{\partial_s u}{\partial_s u}$. Moreover, since the derivative operator is unbounded, the functions $u$ and $u_{ss}$ are independent until the continuity constraints and BCs are enforced. However, accounting for the continuity and BCs is an ad-hoc process, using integration-by-parts or inequalities such as Wirtinger or Poincare. Such ad-hoc methods have been used to generate stability tests for specific classes of PDE models as in~\cite{Papa,datko1970extending,fridman2009exponential,valmorbida2016stability,AHMADI2016163,Gahlawat}, or for control design as in Backstepping methods~\cite{krstic1,aamo,susto2010control,Karafyllis2019,Privat2013,ZHANG}, late-lumping methods~\cite{lasiecka_triggiani_2000}, and port-Hamiltonian methods~\cite{villegas}). Despite these results, no universal approach to computational analysis and control of PDE models exists.


	To summarize, the presence of unbounded operators, continuity constraints, and BCs in a PDE model pose significant challenges to the development of a \textit{universal computational framework} for analysis, control, and simulation. 
 However, these limitations primarily appear in PDEs, and can be remedied by using an alternative modeling framework defined by Partial Integral Equations (PIEs).

	\subsection{The Partial Integral Equation (PIE) Framework} The PIE framework is a state-space approach to modeling spatially distributed systems. PIE models can be considered a generalization of the integro-differential systems that can model phenomena such as elasticity, mechanical fracture, etc.~\cite{appell2000partial, gil2015stability}. Unlike a PDE model, wherein the state is differentiated, consistent with continuity constraints, the state of a PIE model is the highest spatial derivative of the PDE model and this state is integrated in space in order to obtain the evolution equation. Consequently, a PIE model is defined by a single integro-differential equation, is parameterized by the $^*$-algebra of Partial Integral (PI) operators, and can be used to represent almost any well-posed PDE model.

	The simplest form of PIE, ignoring ODEs, inputs, and outputs, is defined by two PI operators, $\mcl T,\mcl A: L_2\rightarrow L_2$ as $\mcl{T}\dot{\mbf{v}}(t) = \mcl{A}\mbf{v}(t)$, where the state, $\mbf v(t) \in L_2$ admits no continuity constraints or BCs. An operator $\mcl P$ is said to be a 3-PI operator, denoted $\mcl P \in \Pi_3$ if there exist $R_0\in L_{\infty}$ and separable functions $R_1,R_2$ such that\vspace{-1mm}
	\[\left(\mcl P\mbf u\right)(s)\hspace{-.4mm} =\hspace{-.4mm} R_0(s)\mbf u(s) +\hspace{-.4mm} \int\limits_a^s \hspace{-1mm} R_1(s,\theta)\mbf u(\theta)\,d\theta +\hspace{-.4mm} \int\limits_s^b \hspace{-1mm} R_2(s,\theta)\mbf u(\theta)\, d\theta.
	\]

Since a PIE representation requires no auxiliary constraints on the state, we can develop universal algorithms for analysis, control, and simulation that apply to any well-posed PIE model. Furthermore, PIEs are parameterized by PI operators that inherit many properties of matrices and form a $^*$-algebra (See Appendix~\ref{app:algebra_pi}) affording us many conveniences in building computational tools using such algorithms.

To illustrate these advantages, consider the PIE model of the heat equation, $u_t=u_{ss}$ with BCs $u(t,0)=u_s(t,1)=0$, and continuity constraint $u \in W_2$. A PIE representation of this PDE model is given by
\begin{equation}
\mcl T v_t(t,s)=\int_{0}^s \theta\, v_t(t,\theta)\,d\theta+\int_s^1 s\, v_t(t,\theta)\,d\theta=-v(t,s).
\end{equation}
To prove the stability using the PIE model in the above equation, we can define a Lyapunov function $V=\ip{\mcl Tv}{\mcl T v}_{L_2}$ and differentiate in time to obtain
\begin{align*}
\dot V(v(t))&=\ip{\mcl T \dot v(t)}{\mcl T v(t)}_{L_2}+\ip{\mcl Tv(t)}{\mcl T \dot v(t)}_{L_2}\\
&=\ip{v(t)}{(\mcl T + \mcl T^*)v(t)}_{L_2}=\ip{v(t)}{\mcl Dv(t)}_{L_2}
\end{align*}
where $\mcl D \in \Pi_3^p$ is parameterized by $R_1(s,\theta)=-2\theta$,  $R_2(s,\theta)=-2s$ and $R_0=0$. We may now use convex optimization to find the PI operator $\mcl Q \in \Pi_3$ such that
$\mcl D  = -\mcl Q^*\mcl Q$. In this case $\mcl Q$ is parameterized by $R_1=\sqrt{2}$, $R_2=0$ and $R_0=0$. This proves that $\dot V(v)=\ip{v}{\mcl D v}=-\ip{Qv}{\mcl Q v}\le 0$.

\subsection{Contribution of this Paper}\label{subsec:intro-contribution}
Because the PIE representation is unified, any algorithm or method designed to analyze, control, or simulate PIE models can be applied to any well-posed system that admits such a representation. The impact of such algorithms and methods can be increased by expanding the class of PDE models for which an equivalent PIE model representation exists.
Unfortunately, however, the class of PDE models for which PDE-PIE conversion formulae exist is still rather limited. To illustrate, consider the following potential models which include integrals, higher order derivatives and disturbances:
\begin{enumerate}
	\item In \cite[Theorem 8]{day2013heat}, a diffusion equation model is given for the evolution of entropy in a 1D thermoelastic rod:
	\begin{flalign*}
	&\dot{\eta}(t,s) = \eta_{ss}(t,s),\qquad  \eta(t,\cdot)\in X&
	\end{flalign*}
 with the solution restricted to the state-space given by
 \[
 X =\left\lbrace\mbf x\in W_2[0,1]\;\mid\; \mbf x(0)=\mbf x(1)=-\int_0^1 \mbf x(s) ds\right\rbrace.
 \]
	\item If one includes boundary actuation ($u\in W_1[0,\infty)$) and multiple disturbances ($d_0, d_1\in L_2[0,\infty),d_2\in W_1[0,\infty)$) in the Timoshenko beam equation, then a possible model can be proposed as
	\begin{flalign*}
	&\ddot{x}(t) = - \mbf{x}_{sss}(t,0) + d_0(t),\quad \{x(t),\mbf x(t,\cdot)\}\in X_{(u(t),d_2(t))},\notag\\
	&\ddddot{\mbf{x}}(t,s) = -\ddot{\mbf{x}}(t,s)+ \ddot{\mbf{x}}_{s}(t,s) -\mbf{x}_{ssss}(t,s) + d_1(t,s),&
	\end{flalign*}
 with the solution restricted to the state-space given by
\[
 X_{(u,v)} =\left\lbrace \mat{x\in \R, \mbf x\in W_4[0,1]\;\mid\; \mbf x(0)=x, \partial_s \mbf x(0)=0,&\hspace{-2.5mm}\\ &\hspace{-6.5cm} \partial_s^2 \mbf x(1)=u, \partial_s^3 \mbf x(1)=v}\right\rbrace.
\]
\end{enumerate}
Such PDE models, with integral terms at the boundary or higher order spatial derivatives are rather simple extensions of standard PDEs. The task of simulation and control for such speculative extensions would require substantial effort because typical methods for such tasks require careful analysis. Assuming well-posedness, this effort could be substantially reduced if the PDE admits a PIE representation. Currently, however, the class of PDE models with known PIE conversion formulae does not include generators with spatial derivatives of order higher than 2, PDEs with inputs and outputs, PDEs coupled with ODEs, or BCs that combine boundary values with inputs and integrals of the state.

The goal of this paper, then, is to facilitate the rapid prototyping and analysis of new and speculative PDE models by extending the class of PDE models for which we have PDE to PIE conversion formulae to include the cases defined above. While we do not address the issue of well-posedness in the paper, this assumption allows one to apply the existing PIE framework for simulation, stability, gain analysis, control, and estimation to such models.

To achieve this goal, we: (a) define a parametric representation for the expanded class of PDE models; (b) define an appropriate state-space to be used in the corresponding PIE model; (c) find a unitary transformation from the PIE state-space to the state-space of the PDE model -- proving equivalence of solutions and equivalence of stability properties. In this context, the main contributions of the paper are:

\begin{enumerate}
	\item \emph{A unified class of PDE models}: We parameterize a class of linear PDE models, referred to as Generalized Partial Differential Equations (GPDEs), encompassing: ODEs coupled with PDEs, $N$th-order spatial derivatives, integrals of the state, control inputs and disturbances, and sensed and regulated outputs.	
	
	\item \emph{Formulae to convert GPDE models to PIEs}:  Given a sufficiently well-posed GPDE model, we generalize Fundamental Theorem of Calculus (FTC) and give formulae (with a Matlab modeling interface) for conversion to a PIE -- \Cref{sec:pde2pie,sec:gpde2pie}.

	\item \emph{Equivalence of GPDEs and PIEs}:  We show the solution of a GPDE uniquely maps to a PIE solution via a unitary map. Additionally, we prove input-output and internal stability of the GPDE model are equivalent to that of the associated PIE -- \Cref{sec:equivalence}.
\end{enumerate}
Before looking at the class of models considered in this paper, we will look at the notation used throughtout the paper and then discuss the main results. For brevity, we only provide proofs for the main results and an outline for other Corollaries/Lemmas. Refer to the Appendices included in the supplementary documents or the full version~\cite{shivakumar_2022GPDE} for extended proofs and non-essential definitions.



	\section{Notation, PI operators and PIEs}\label{sec:notation}
	The empty set, $\emptyset$, is occasionally used to denote a matrix or matrix-valued function with either zero row or column dimension and whose non-zero dimension can be inferred from context. $0_{m,n} \in \R^{m\times n}$ is the matrix of all zeros, $0_n:=0_{n,n}$, and $I_n \in \R^{n\times n}$ is the identity matrix. We use $0$ and $I$ for these matrices when dimensions are clear from context. $\R_+$ is the set of non-negative real numbers.
	The set of $k$-times continuously differentiable n-dimensional vector-valued functions on the interval $[a,b]$ is denoted by $C_k^{n}[a,b]$. $L_2^n[a,b]$ is the Hilbert space of $n$-dimensional vector-valued Lebesgue square-integrable functions on the interval $[a,b]$ equipped with the standard inner product. $L_{\infty}^{m,n}[a,b]$ is the Banach space of $m\times n$-dimensional essentially bounded measurable matrix-valued functions on $[a,b]$ equipped with the essential supremum singular value norm.

	We utilize different fonts to denote different mathematical objects: $u$ or $u(t)$ typically denotes a scalar or finite-dimensional vector; and, $\mbf x$ or $\mbf x(t)$ typically denotes a scalar or vector-valued function. For a function, $\mbf x$, of spatial variable $s$, we use $\partial_s^j \mbf x$ to denote the $j$th order partial derivative $\frac{\partial^j\mbf x}{\partial s^j}$. For a function of time and possibly space, we denote  $\dot{\mbf x}(t)=\frac{\partial }{\partial t}\mbf x(t)$. We use $W_{k}^n$ to denote the Sobolev spaces $W_{k}^n[a,b] := \{\mbf u\in L_2^n[a,b]\mid \partial^i_s \mbf u\in L_2^n[a,b] ~\forall ~i\le k\}$ with inner product $\ip{\mbf u}{\mbf v}_{W_{k}^n}=\sum\nolimits_{i=0}^k \ip{\partial_s^i\mbf u}{\partial_s^i \mbf v}_{L_2^n}$.
	Clearly, $W_{0}^n[a,b] = L_2^n[a,b]$. For given $n =\{n_0, \cdots, n_N\} \in \N^{N+1}$, we define the Cartesian product space $W^{n}:=\prod_{i=0}^{N}W_i^{n_i}$ and for $\mbf u=\{\mbf u_0,\cdots,\mbf u_N\}\in W^{n}$ and $\mbf v=\{\mbf v_0,\cdots,\mbf v_N\}\in W^{n}$ we define the associated inner product as $\ip{\mbf u}{\mbf v}_{W^{n}}=\sum\nolimits_{i=0}^N \ip{\mbf u_i}{\mbf v_i}_{W_i^{n_i}}$.
	We use $\R L_2^{m,n}[a,b]$ to denote the space $\R^{m}\times L_2^{n}[a,b]$ and for $\mbf x=\bmat{x_1\\\mbf x_2}\in \R L_2^{m,n}$ and $ \mbf y=\bmat{y_1\\\mbf y_2}\in \R L_2^{m,n}$, we define the associated inner product as
\[\ip{\bmat{x_1\\\mbf x_2}}{\bmat{y_1\\\mbf y_2}}_{\R L_2^{m,n}} = x_1^Ty_1 + \ip{\mbf x_2}{\mbf y_2}_{L_2^n}.
\]
Frequently, we omit the domain $[a,b]$ and simply write $L_2^n$, $W_{k}^n$, $W^n$, or $\R L_2^{m,n}$.	For functions of time only ($L_2[\R_+]$ and $W_k[\R_+]$), 
to denote the extended subspaces of such functions by $L_{2e}[\R_+]$ and $W_{ke}[\R_+]$ respectively as
\[L_{2e}[\R_+] := \left\lbrace x\;\mid \; x \in L_2[0,T] \; \forall\; T\ge 0 \right\rbrace,\]
\[W_{ke}[\R_+] :=\left\lbrace x\;\mid\; x\in W_{k}[0,T] \; \forall \; T\ge0\right\rbrace.\vspace{-2mm}\]
Finally, for normed spaces $A,B$, $\mcl L(A,B)$ denotes the space of bounded linear operators from $A$ to $B$ equipped with the induced operator norm. $\mcl L(A) := \mcl L(A,A)$.
\subsection{PI Operators: A $^*$-algebra of bounded linear operators}\label{sec:PI}
PI algebras are parameterized classes of operators on $\R L_2^{m,n}$. The first of these is the algebra of 3-PI operators that map $L_2^n \rightarrow L_2^n$ with separable functions in $L_{\infty}$ as parameters.

\begin{definition}[Separable Function]\label{def:separable}
We say $R:[a,b]^2 \rightarrow \R^{p\times q}$ is separable if there exist $r \in \N$, $F\in L_{\infty}^{r\times p}[a,b]$  and $G\in L_{\infty}^{r\times q}[a,b]$  such that $R(s,\theta)=F(s)^T G(\theta)$.
\end{definition}

\begin{definition}[3-PI operators, $\Pi_3$]\label{def:3PI}
Given $R_0\in L_{\infty}^{p \times q}[a,b]$ and separable functions $R_1, R_2: [a,b]^{2} \rightarrow \R^{p \times q}$, we define the operator $\mcl P=\threepi{R_i}$ for $\mbf v \in L_2[a,b]$ as
\begin{align}
&\left(\threepi{R_i}\mbf v\right)(s):=\label{eq:3pi}\\
& R_0(s) \mbf v(s) +\int_{a}^s R_1(s,\theta)\mbf v(\theta)d \theta+\int_s^bR_2(s,\theta)\mbf v(\theta)d \theta.\notag
\end{align}
Furthermore, we say an operator, $\mcl P$, is 3-PI of dimension $p \times q$, denoted $\mcl P \in [\Pi_3]_{p,q}\subset \mcl L(L_2^q,L_2^p)$, if there exist functions $R_0$ and separable functions $R_1,R_2$ such that $\mcl P=\threepi{R_i}$.
\end{definition}

For any $p\in \N$, $[\Pi_3]_{p,p}$ is a $^*$-algebra, being closed under addition, composition, scalar multiplication, and adjoint (See \Cref{app:algebra_pi}~\cite{shivakumar_2022GPDE}). 
%
%
Likewise, the algebra of 3-PI operators can be extended to $\mcl L(\R L_2)$ as follows.

\begin{definition}[4-PI operators]\label{def:4PI} Given $P\in\R^{m\times n}$, $Q_1 \in L_{\infty}^{m\times q}$, $Q_2\in L_{\infty}^{p\times n}$, and $R_0,R_1,R_2$ with $\threepi{R_i}\in [\Pi_3]_{p,q}$, we say $\mcl P=\fourpi{P}{Q_1}{Q_2}{R_i}\in \mcl L(\R L_2^{m,p},\R L_2^{n,q})$ if 
\begin{align}\label{eq:4pi}
\left(\mcl P\bmat{u\\\mbf{v}}\right)(s)
& := \bmat{Pu + \int_{a}^{b}Q_1(\theta)\mbf{v}(\theta)d\theta\\Q_2(s)u+\left(\threepi{R_i}\mbf{v}\right) (s)}.
\end{align}
Furthermore, we say $\mcl P$, is 4-PI, denoted $\mcl P \in [\Pi_4]^{m,n}_{p,q}$, if there exist $P,Q_1,Q_1,R_0,R_1,R_2$  such that $\mcl P=\fourpi{P}{Q_1}{Q_2}{R_i}$.
\end{definition}

\begin{definition}[*-subalgebras of $\Pi_i$ with polynomial parameters]\label{def:4PI_poly} We say $\mcl P \in [\Pi_3^p]_{p,q}$ if there exist polynomials $R_i$ of appropriate dimension such that $\mcl P=\threepi{R_i}$. We say $\mcl P \in [\Pi_4^p]_{p,q}^{m,n}$ if there exist matrix $P$ and polynomials $Q_i,R_i$ of appropriate dimension such that $\mcl P=\fourpi{P}{Q_1}{Q_2}{R_i}$.
\end{definition}

Algebraic operations on $\Pi_i$ are defined by algebraic operations on the parameters that represent these operators. To illustrate, consider how composition of two operators in $\Pi_4$ defines a map on the parameters which define those operators (this map will be used later in the article). To define this map, we first note that any $\Pi_4$ operator has an associated set of matrix and polynomial parameters which lie in the space
\begin{align*}
[\Gamma_4]^{m,p}_{n,q}&:=\left\lbrace \mat{{\small\left[\begin{array}{c|c}P & Q_1 \\\hline  Q_2 & \{R_0,R_{1a}R_{1b},R_{2a}R_{2b}\}\end{array}\right]} \;: \\\hspace{5mm} P\in \R^{m\times n},\,Q_1 \in L_{\infty}^{m\times q},\,Q_2\in L_{\infty}^{p\times n},\\\hspace{5mm}R_0\in L_{\infty}^{q\times n},\,R_{ia}\in L_{\infty}^{q\times n_b},R_{ib}\in L_{\infty}^{n_b\times n}}\right\rbrace.
\end{align*}
Then, given any $\left[\footnotesize\begin{array}{c|c}A&B_1\\\hline B_2&\{C_i\}\end{array}\right]\in [\Gamma_4]^{m,p}_{n,q}$ and $\left[\footnotesize\begin{array}{c|c}P&Q_1\\\hline Q_2&\{R_i\}\end{array}\right]\in [\Gamma_4]^{p,l}_{q,k}$, if the parametric map $\mbf P^4_\times:[\Gamma_4]^{m,p}_{n,q} \times [\Gamma_4]^{p,l}_{q,k}\rightarrow [\Gamma_4]^{m,l}_{n,k}$, is defined as
\begin{align}\label{eq:pi_comp_map}
\left[\footnotesize\begin{array}{c|c}\hat{P}&\hat{Q}_1\\\hline \hat{Q}_2&\{\hat{R}_i\}\end{array}\right]=\picomp{A}{B_1}{B_2}{C_i}{P}{Q_1}{Q_2}{R_i}
\end{align}
where
\begin{align*}
&\hat{P} = AP + \int_a^b B_1(s)Q_2(s)ds,\quad \hat{R}_0(s) = C_0(s)R_0(s),\\
&\hat{Q}_1(s) = AQ_1(s) + B_1(s)R_0(s)+\int_{s}^b B_1(\eta)R_1(\eta,s)d\eta\\
&\qquad\qquad+\myinta{s}B_1(\eta)R_2(\eta,s)d\eta,\\
&\hat{Q}_2(s) = B_2(s)P + C_0(s)Q_2(s) + \int_a^s C_1(s,\eta)Q_2(\eta)d\eta\\
&\qquad\qquad+\int_s^b C_2(s,\eta)Q_2(\eta)d\eta,\\
&\hat{R}_1(s,\eta) =B_2(s)Q_1(\eta)+C_0(s)R_1(s,\eta)+C_1(s,\eta)R_0(\eta)\\
&+\int_a^{\eta} C_1(s,\theta)R_2(\theta,\eta)d\theta+\int_{\eta}^{s}C_1(s,\theta)R_1(\theta,\eta)d\theta\\
&\qquad\qquad+\int_{s}^bC_2(s,\theta)R_1(\theta,\eta)d\theta,\\
&\hat{R}_2(s,\eta) =B_2(s)Q_1(\eta)+C_0(s)R_2(s,\eta)+C_2(s,\eta)R_0(\eta)\\
&+\int_a^{s} C_1(s,\theta)R_2(\theta,\eta)d\theta+\int_{s}^{\eta}C_2(s,\theta)R_2(\theta,\eta)d\theta\\
&\qquad\qquad+\int_{\eta}^bC_2(s,\theta)R_1(\theta,\eta)d\theta,
\end{align*} we have 
 \begin{align*}
&\pie\left[\footnotesize\begin{array}{c|c}\hat{P}&\hat{Q}_1\\\hline \hat{Q}_2&\{\hat{R}_i\}\end{array}\right]=\fourpi{A}{B_1}{B_2}{C_i}\left(\fourpi{P}{Q_1}{Q_2}{R_i}\right).
\end{align*}

\subsection{Partial Integral Equations}\label{sec:PIE}
A Partial Integral Equation (PIE) is an extension of the state-space representation of ODEs (vector-valued first-order differential equations on $\R^n$) to spatially distributed states on the product space $\R L_2$. Mirroring the 9-matrix optimal control framework developed for state-space systems, a PIE system is parameterized by twelve 4-PI operators as\vspace{-2mm}

{\small
\begin{align}\label{eq:PIE_full}		
\bmat{\mcl{T}\dot{\ubar{\mbf{x}}}(t)\\z(t)\\y(t)}&=\bmat{\mcl{A}&\mcl{B}_1&\mcl{B}_2\\\mcl{C}_1&\mcl{D}_{11}&\mcl{D}_{12}\\\mcl{C}_2&\mcl{D}_{21}&\mcl{D}_{22}}\bmat{\ubar{\mbf x}(t)\\w(t)\\u(t)}-\bmat{\mcl{T}_{w}\dot{w}(t)+\mcl{T}_{u}\dot{u}(t)\\0\\0},\notag\\
&\qquad \ubar{\mbf{x}}(0) = \ubar{\mbf{x}}^0\in \R L_2^{m,n}[a,b],
\end{align}}
where $z(t)\in \R^{n_z}$ is the regulated output, $y(t)\in \R^{n_y}$ is the sensed output, $w(t)\in \R^{n_w}$ is the disturbance, $u(t)\in \R^{n_u}$ is the control input, and $\ubar{\mbf x}(t)\in \R L_{2}^{n_x, n_{\hat{\mbf x}}}$ is the internal state.


Let us note two significant features of this model. First, through some slight abuse of notation, in this paper, we will use expressions such as $\mcl T \dot{\ubar{\mbf x}}$ to represent $\partial_t (\mcl T \ubar{\mbf x})$. Second, we observe that PIEs allow for the dynamics to depend on the time-derivative of the input signals: $\partial_t (\mcl T_w w)$ and $\partial_t (\mcl T_u u)$. These terms are included to allow for the PIEs to represent certain classes of PDEs wherein signals enter through the BCs.


\noindent\textbf{Notation:} Finally, we collect the 12 PI parameters which define a PIE system in Eq.~\eqref{eq:PIE_full} and introduce the shorthand notation $\mbf G_{\mathrm{PIE}}$ which represents the labeled tuple of such system parameters as
\begin{align*}
	&\mbf G_{\mathrm{PIE}} = \left\{\mcl T, \mcl T_w, \mcl T_u, \mcl A, \mcl B_1, \mcl B_2, \mcl C_1, \mcl C_2, \mcl D_{11}, \mcl D_{12}, \mcl D_{21}, \mcl D_{22}\right\}.
	\end{align*}
	When this shorthand notation is used, it is presumed that all parameters have appropriate dimensions.

	We now define minimal requirements that a solution of a PIE system should satisfy.

	\begin{definition}[Admissible solution for a PIE system]\label{def:piesolution}
	For given inputs $u\in L_{2e}^{n_u}[\R_+]$, $w\in L_{2e}^{n_w}[\R_+]$ with $(\mcl T_u u)(\cdot,s)\in W_{1e}^{n_x+n_{\hat{\mbf x}}}[\R_+]$ and $(\mcl{T}_{w} w)(\cdot,s)\in W_{1e}^{n_x+n_{\hat{\mbf x}}}[\R_+]$ for all $s \in [a,b]$ and $\ubar{\mbf x}^0(t)\in \R L_2^{n_x, n_{\hat{\mbf x}}}$, we say that $\{\ubar{\mbf x}, z, y\}$ satisfies the PIE defined by $\mbf G_{\mathrm{PIE}}$ with initial condition $\ubar{\mbf x}^0$ and input $\{w, u\}$ if $z\in L_{2e}^{n_z}[\R_+]$, $y\in L_{2e}^{n_y}[\R_+]$, $\ubar{\mbf x}(t)\in \R L_2^{n_x, n_{\hat{\mbf x}}}[a,b]$ for all $t\ge 0$, $\mcl T\ubar{\mbf x}$ is Frech\'et differentiable on $\R_+$, $\mcl T\ubar{\mbf x}(0) = \mcl T\ubar{\mbf x}^0$, and Eq.~\eqref{eq:PIE_full} is satisfied for all $t\in\R_+$.
	\end{definition}

 Note that the Frech\'et (and not Gateaux) derivative is used in Defn.~\ref{def:piesolution} because we restrict our analysis to linear systems.

	\section{GPDEs: A Generalized Class of Linear Models}\label{sec:ODEPDE}
	Now, we parameterize the class of PDE models for which we may define associated PIE systems. To simplify the notation and analysis, we will represent these models as the interconnection of ODE and PDE subsystems -- See \Cref{fig:block_diagram}. This class of PDE models will be referred to as Generalized Partial Differential Equations (GPDEs). 

	\subsection{ODE Subsystem}\label{subsec:ODE}
\begin{figure}[t]
	\centering
	\includegraphics[width=0.75\columnwidth]{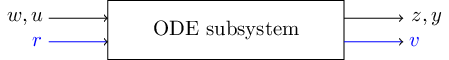}
	\caption{Depiction of the ODE subsystem for use in defining a GPDE. All external input signals in the GPDE model pass through the ODE subsystem and are labeled as $u$ and $w$, corresponding to control input and disturbance/forcing input. All external outputs pass through the ODE subsystem and are labeled $y$ and $z$, corresponding to measured output and regulated output. All interaction with the PDE subsystem is routed through two vector-valued signals: $r$ the sole output of the PDE subsystem and $v$ the sole input to the PDE subsystem.}
	\label{fig:ode_subsystem}
	\end{figure}

	The ODE subsystem of a GPDE model, illustrated in \Cref{fig:ode_subsystem}, is a typical state-space representation with real-valued finite-dimensional inputs and outputs. These inputs and outputs include both the interconnection with the PDE subsystem and the inputs and outputs of the GPDE model as a whole. We partition both the input and output signals into 3 components, differentiating these channels by function. The input channels are: the control input to the GPDE ($u(t)\in \R^{n_u}$); the exogenous disturbance ($w(t)\in \R^{n_w}$); and the internal feedback input ($r(t)\in \R^{n_r}$) which is the output of the PDE subsystem. The output channels of the ODE subsystem are: the regulated output of the GPDE ($z(t)\in \R^{n_z}$); the sensed outputs of the GPDE ($y(t)\in \R^{n_y}$); and the output from the ODE subsystem which becomes the input to the PDE subsystem ($v(t)\in \R^{n_v}$).

	\begin{definition}[Admissible solution for an ODE subsystem]
	Given matrices $A$, $B_{xw}$, $B_{xu}$, $B_{xr}$, $C_z$, $D_{zw}$, $D_{zu}$, $D_{zr}$, $C_y$, $D_{yw}$, $D_{yu}$, $D_{yr}$, $C_v$, $D_{vw}$, $D_{vu}$ of appropriate dimension, we say $\{x,z,y,v\}$ with $\{x(t),z(t),y(t),v(t)\}\in \R^{n_x}\times \R^{n_z}\times \R^{n_y}\times \R^{n_v}$ satisfies the ODE with initial condition $x^0 \in \R^{n_x}$ and input $\{w,u,r\}$ if $x$ is differentiable, $x(0)=x^0$  and for $t\ge 0$
	\begin{flalign}\label{eq:odepde_general}
	\bmat{\dot{x}(t)\\\hline z(t)\\y(t)\\v(t)} &= \bmat{A&\vlines&B_{xw}&B_{xu} & B_{xr}\\\hline C_z&\vlines&D_{zw}&D_{zu}&D_{zr}\\C_y&\vlines&D_{yw}&D_{yu}&D_{yr}\\C_v&\vlines&D_{vw}&D_{vu}&0}\bmat{x(t)\\\hline w(t)\\u(t)\\r(t) }.
	\end{flalign}
	\end{definition}

	\noindent\textbf{Notation:} We collect all matrix parameters from the ODE subsystem in~\eqref{eq:odepde_general} and introduce the shorthand notation $\mbf G_{\mathrm{o}}$ which represents the labelled tuple of such parameters as
	\begin{align}
	\mbf G_{\mathrm{o}} &= \left\{A, B_{xw}, B_{xu}, B_{xr}, C_z, D_{zw}, D_{zu}, D_{zr}, C_y, D_{yw}, D_{yu},\right. \notag\\
	&\qquad\left. D_{yr}, C_v, D_{vw}, D_{vu}\right\}.\label{eq:ode-params}
	\end{align}
	When this shorthand notation is used, it is presumed that all parameters have appropriate dimensions. 

	\subsection{PDE Subsystem}\label{subsec:PDE} 
\begin{figure}[t]
	\centering
	\includegraphics[width=0.75\columnwidth]{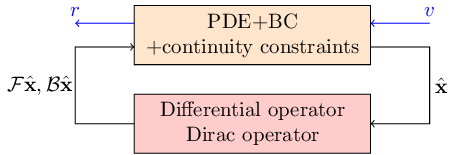}
	\caption{Depiction of the PDE subsystem for use in defining a GPDE. All interaction of the PDE subsystem with the ODE subsystem is routed through the two vector-valued signals: $r(t)$ an output of the PDE subsystem (and input to the ODE subsystem) and $v(t)$ is input to the PDE subsystem (and output from the ODE subsystem).} 
	\label{fig:pde_subsystem}
	\end{figure}

Our parameterization of the PDE subsystem, illustrated in Fig.~\ref{fig:pde_subsystem}, is divided into three parts: the continuity constraints, the in-domain dynamics, and the BCs. The continuity constraints specify the existence of partial derivatives and boundary values for each state as required by the in-domain dynamics and BCs. The in-domain dynamics (or generating equation) specify the time derivative of the state at every point in the interior of the domain, and are expressed using integral, dirac, and $N^{\text{th}}$-order spatial derivative operators. The BCs are represented as a real-valued algebraic constraint subsystem that maps the distributed state and inputs to a vector of boundary values. External inputs or outputs are not defined for the PDE subsystem and instead are routed through the ODE subsystem using the interconnection signals.

	\subsubsection{The continuity constraint} The `continuity constraint' partitions the state vector of the PDE subsystem, $\hat{\mbf x}(t,\cdot)$, and specifies the differentiability properties of each partition as required for existence of the required partial derivatives and limit values. This partition is defined by the parameter $n\in \N^{N+1}=\{n_0,\cdots n_N\}$, wherein $n_i$ specifies the dimension of the $i$th partition vector so that $\hat{\mbf x}_i(t,s)\in \R^{n_i}$. The partitions are ordered by increasing differentiability so that
	\[
	\hat{\mbf x}(t,\cdot)=\bmat{\hat{\mbf x}_0(t,\cdot)\\ \vdots\\\hat{\mbf x}_N(t,\cdot)} \in W^{n}:=\bmat{W_0^{n_0}\\ \vdots \\ W_N^{n_N}}.
	\]
	Given the partition defined by $n\in \N^{N+1}$, and given $\hat{\mbf x}\in W^n$, we would like to list all well-defined partial derivatives of $\hat{\mbf x}$. \textbf{Note} that the Sobolev spaces $W_N$ have a continuous embedding into the space of continuous functions, $W_N[a,b]\subset C_{N-1}[a,b]$ -- implying existence of boundary values $\partial_s^i \mbf x(a)$ and classical partial derivatives $\partial_s^i \mbf x$ for $i< N$ above and also in Lem.~\ref{lem:cauchy}.
 
    To do this, we first define $n_{\hat{\mbf x}}:=|n|_1=\sum_{i=0}^N n_i$ to be the number of states in $\hat{\mbf x}$, $n_{S_i}:=\sum_{j=i}^N n_{j}\le n_{\hat{\mbf x}}$ to be the total number of $i$-times differentiable states, and $n_{S}=\sum_{i=1}^N n_{S_i}$ to be the total number of possible partial derivatives of $\hat{\mbf x}$ as permitted by the continuity constraint.

	\textbf{Notation:} For indexed vectors (such as $n$ or $\hat{\mbf x}$) we use $\hat{\mbf x}_{i:j}$ to denote the components $i$ to $j$. Specifically, $\hat{\mbf x}_{i:j} = \text{col}(\hat{\mbf x}_i,\cdots, \hat{\mbf x}_j)$,  $n_{i:j}: = \sum_{k=i}^j n_k$ and $n_{S_{i:j}} = \sum_{k=i}^j n_{S_k}$.

	Next, we define the selection operator  $S^i: \R^{n_{\hat{\mbf x}}} \rightarrow \R^{n_{Si}}$ which is used to select only those states in $\hat{\mbf x}$ which are at least $i$-times differentiable. Specifically, for $\hat{\mbf x} \in W^n$, we have
	$S^i = \bmat{0_{{n_{S_i}} \times n_{\hat{\mbf x}}-n_{S_i}}& I_{n_{S_i}}}$, so that	$(S^i\hat{\mbf x})(s)=\bmat{\hat{\mbf x}_i(s)\\ \vdots\\\hat{\mbf x}_N(s)}.$
	We may now conveniently represent all well-defined $i$th-order partial derivatives of $\hat{\mbf x}$ as $\partial_s^i S^i\hat{\mbf x}$ so that\\
	$(\partial_s^i S^i\hat{\mbf x})(s)=\bmat{\partial_s^i \hat{\mbf x}_i(s)\\ \vdots\\ \partial_s^i \hat{\mbf x}_N(s)}$ and $(\mcl F\hat{\mbf x})(s):=\bmat{\hat{\mbf x}(s)\\ (\partial_s S\hat{\mbf x})(s) \\ \vdots \\ (\partial_s^{N}S^{N}\hat{\mbf x})(s)}$
	where $\mcl F$ concatenates all the $\partial_s^i S^i\hat{\mbf x}$ for $i=0,\cdots,N$ --- creating an ordered list including both the PDE state, $\hat{\mbf x}$, as well as all $n_S$ possible partial derivatives of $\hat{\mbf x}$ as permitted by the continuity constraint and the vector $(\mcl F\hat{\mbf x})(s)\in \R^{n_S+n_{\mbf x}}$.

	Lastly, we may construct $(\mcl C\hat{\mbf x})(s) \in \R^{n_S}$, the vector of all absolutely continuous functions generated by $\hat{\mbf x}$ and its partial derivatives.
	Using $\mcl C \hat{\mbf x}$, we may then construct $\mcl B\hat{\mbf x} \in \R^{2n_S}$, the list all possible boundary values of $\hat{\mbf x} \in W^n$. Specifically,  $\mcl C \hat{\mbf x}$ and $\mcl B\hat{\mbf x}$ are defined as
	\begin{flalign}\label{eq:odepde_bc}
		\mcl C\hat{\mbf x}(s)&=\bmat{(S\hat{\mbf x})(s)\\ (\partial_s S^2\hat{\mbf x})(s) \\ \vdots \\ (\partial_s^{N-1}S^{N}\hat{\mbf x})(s)} ~ \text{and} ~
		\mcl B\hat{\mbf x}=\bmat{(\mcl C\hat{\mbf x})(a)\\(\mcl C\hat{\mbf x})(b)}.
		\end{flalign}
  This notation also allows us to specify all well-defined boundary values of $\hat{\mbf x}\in W^n$ and of its partial derivatives. 

		\subsubsection{Boundary Conditions (BCs)}
		Given the notational framework afforded by the continuity condition, and equipped with our list of well-defined terms ($\mcl F\hat{\mbf x}$ and $\mcl B \hat{\mbf x}$), we may now parameterize a generalized class of BCs consisting of a combination of boundary values, integrals of the PDE state, and the effect of the input signal, $v$. Specifically, the BCs are parameterized by the square-integrable function $B_{I}:[a,b]\to \R^{n_{BC} \times (n_{S}+n_{\hat{\mbf x}})}$ and matrices $B_v\in \R^{n_{BC}\times n_v}$ and $B\in \R^{n_{BC}\times 2n_S}$ as
		\begin{equation}\label{eq:odepde_b}
		\int_{a}^{b} B_{I}(s)(\mcl F\hat{\mbf x}(t))(s)ds +\bmat{B_v & -B}\bmat{v(t)\\ \mcl B \hat{\mbf x}(t)}=0
		\end{equation}
		where $n_{BC}$ is the number of specified BCs. For reasons of well-posedness, as discussed in \Cref{sec:pde2pie}, we typically require $n_{BC}=n_{S}$. If fewer BCs are available, it is likely that the continuity constraint is too strong.

		Now that we have parameterized a general set of BCs, we embed these BCs in what is typically referred to as the domain of the infinitesimal generator -- which combines the BCs and continuity constraints into a set of acceptable states.
		\begin{equation}\label{eq:odepde_general_domain}
		X_{v} := \left\lbrace \mat{\hat{\mbf x}\in W^{n}[a,b]: &\\
		&\hspace{-25mm} \int_{a}^{b} B_{I}(s)(\mcl F\hat{\mbf x})(s)ds +\bmat{B_v & -B}\bmat{v\\ \mcl B \hat{\mbf x}}=0}\right\rbrace.
		\end{equation}
		The set $X_{v}$ is used to restrict the state and initial conditions as $\hat{\mbf x}(t)\in X_{v(t)}$ and $\hat{\mbf x}(0) = \hat{\mbf x}^0 \in X_{v(0)}$.

		\noindent\textbf{Notation:} We collect all the parameters which define the constraint in Eq.~\eqref{eq:odepde_b} in $\mbf G_{\mathrm{b}}$, as
		\begin{equation}
		\mbf G_{\mathrm{b}} = \left\{B,~B_{I},~B_v\right\}.\label{eq:BC-parms}
		\end{equation}
		When this shorthand notation is used, it is presumed that all parameters have appropriate dimensions.

		\subsubsection{In-Domain Dynamics of the PDE Subsystem}
		Having specified the continuity constraint and BCs using $\{n,\mbf G_{\mathrm{b}}\}$, we once again use our list of well-defined terms ($\mcl F\hat{\mbf x}$ and $\mcl B \hat{\mbf x}$) to define the in-domain dynamics of the PDE subsystem and the output to the ODE subsystem. These dynamics are parameterized by the functions $A_{0}(s),$ $A_{1}(s,\theta),$ $A_{2}(s,\theta)$ $\in \R^{n_{\hat{\mbf x}} \times (n_{S}+n_{\hat{\mbf x}})}$, $C_{r}(s)$ $\in\R^{n_r\times (n_{S}+n_{\hat{\mbf x}})}$, $B_{xv}(s)$ $\in\R^{n_{\hat{\mbf x}}\times n_v}$, $B_{xb}(s)$ $\in\R^{n_{\hat{\mbf x}}\times 2n_S}$, and matrices $D_{rv}$ $\in\R^{n_r \times n_v}$ and $D_{rb}(s)$ $\in\R^{n_r \times 2 n_S}$ as
		\begin{flalign}
		&\bmat{\dot{\hat{\mbf x}}(t,s)\\ r(t)}\hspace*{-1pt} =\hspace*{-1pt}  \bmat{A_{0}(s)(\mcl F\hat{\mbf x}(t))(s) \\0}+\bmat{B_{xv}(s) & \hspace*{-5pt}B_{xb}(s)\\ 0&\hspace*{-5pt} D_{rb} }\bmat{v(t)\\ \mcl B \hat{\mbf x}(t)}\notag\\
		&+\hspace{-.5mm}\bmat{\int\limits_{a}^{s}\hspace{-.5mm}A_{1}(s,\theta)(\mcl F\hat{\mbf x}(t))(\theta) d\theta+\hspace{-.75mm}\int\limits_{s}^{b} \hspace{-.5mm}A_{2}(s,\theta)(\mcl F\hat{\mbf x}(t))(\theta) d\theta\hspace{-.25mm}\\	
		\int_{a}^{b} C_{r}(\theta)(\mcl F\hat{\mbf x}(t))(\theta) d\theta}.\hspace{-1.5mm}\label{eq:general_pde_subsystem}
		\end{flalign}

		Many commonly used PDE models are defined solely by $A_{0}$. For example, if we consider $u_t=\lambda u+u_{ss}$, then $A_0=\bmat{\lambda&0&1}$ and all other parameters are zero. While the parameters other than $A_0$, appear infrequently in application, they can be used to model non-local effects of distributed state, forcing function, boundary values, etc.


		\noindent\textbf{Notation:} We collect all parameters from the in-domain dynamics of the PDE subsystem (Eq.~\eqref{eq:general_pde_subsystem}) in $\mbf G_{\mathrm{p}}$ as
		\begin{align}
	&\mbf G_{\mathrm{p}} = \left\{A_{0},~A_{1},~A_{2},~B_{xv},~B_{xb},~C_{r},~D_{rb}\right\}.\label{eq:pde-params}
	\end{align}
	When this shorthand notation is used, it is presumed that all parameters have appropriate dimensions. 

	\begin{definition}[Admissible solution for a PDE subsystem]\label{defn:PDE}
	For given $\hat{\mbf x}^0 \in X_{v(0)}$ and $v\in L_{2e}^{n_v}[\R_+]$ with $B_v v\in W_{1e}^{2n_S}[\R_+]$, we say that $\{\hat{\mbf x}, r\}$ satisfies the PDE subsystem defined by $n\in \N^{N+1}$ and $\{\mbf G_{\mathrm b}, \mbf G_{\mathrm{p}}\}$ (defined in Eqs.~\eqref{eq:BC-parms} and \eqref{eq:pde-params}) with initial condition $\hat{\mbf x}^0$ and input $v$ if $r\in L_{2e}^{n_r}[\R_+]$, $\hat{\mbf x}(t) \in X_{v(t)}$ for all $t\ge 0$, $\hat{\mbf x}$ is Frech\'et differentiable with respect to the $L_2$-norm on $\R_+$, $\hat{\mbf x}(0) = \hat{\mbf x}^0$, and Eq.~\eqref{eq:general_pde_subsystem} is satisfied for all $t\ge 0$.
	\end{definition}

Note that the Defn.~\ref{defn:PDE} defines a minimal criterion for the admissibility of solutions of the PDE to guarantee the existence of an associated PIE solution. While Defn.~\ref{defn:PDE} is not intended as a formal statement of the Cauchy problem, the required properties of the solution are similar to the regularity conditions used in classical definitions of solution such as in, e.g., see \cite[Theorem 3.3.3]{zwart-curtain}, \cite[Proposition 10.1.8]{marius}, and \cite{cheng_wellposedness}. 

	\subsection{GPDE: Interconnection of ODE and PDE Subsystems}\label{subsec:interconnection}
\begin{figure}[t]
	\centering
	\includegraphics[width=0.75\columnwidth]{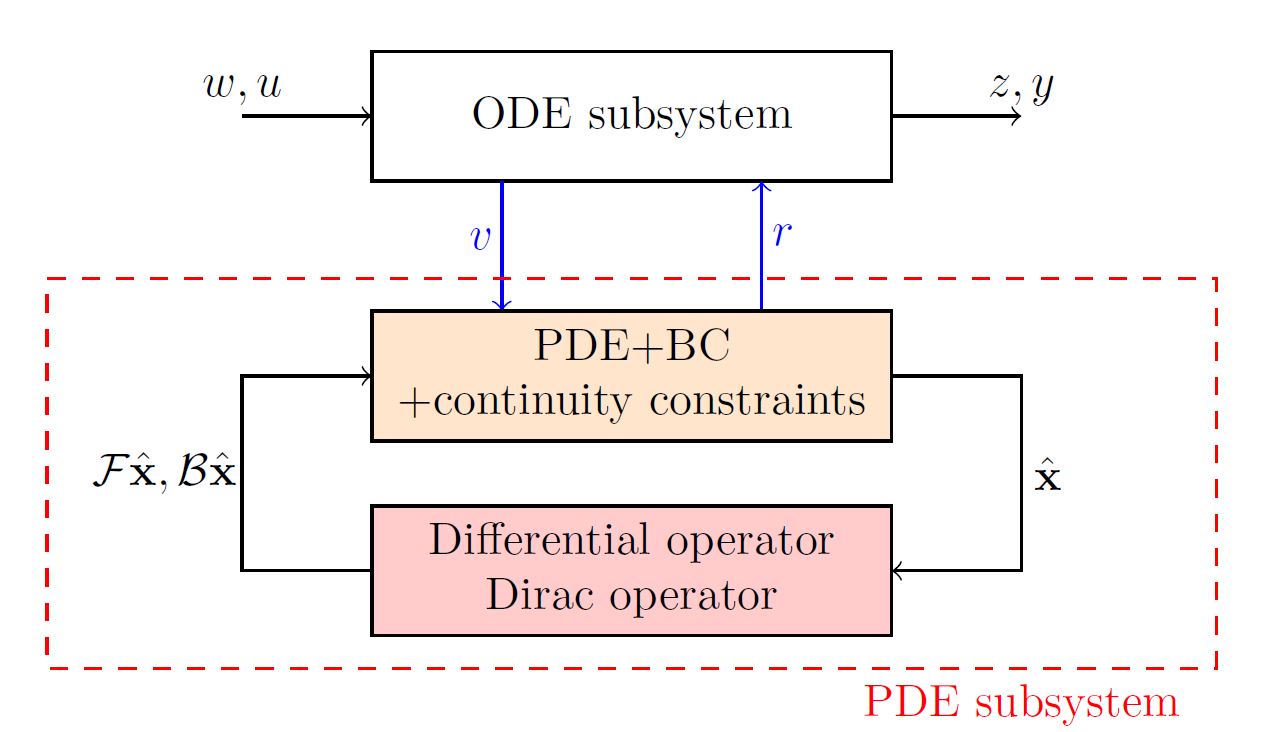}
	\caption{A GPDE is the interconnection of an ODE subsystem (an ODE with finite-dimensional inputs $w,u,v$ and outputs $z,y,r$) with a PDE subsystem ($N$th-order PDEs and BCs with finite-dimensional input $r$ and output $v$). The BCs and internal dynamics of the PDE subsystem are specified in terms of all well-defined spatially distributed terms as encoded in $\mcl F \hat{\mbf x}(t)$ and all well-defined limit values as encoded in $\mcl B \hat{\mbf x}(t)$.}
	\label{fig:block_diagram}
	\end{figure}

	A GPDE model, illustrated in \Cref{fig:block_diagram}, is the mutual interconnection of the ODE and PDE subsystems through the interconnection signals $(r, v)$ and is collectively defined by Eqs.~\eqref{eq:odepde_general}-\eqref{eq:general_pde_subsystem}.
	Given suitable inputs $w$,$u$, for a GPDE model, parameterized by $\{n,\mbf G_{\mathrm{o}}, \mbf G_{\mathrm b},\mbf G_{\mathrm{p}}\}$, we define the continuity constraint and time-varying BCs by $\{x(t),\hat{\mbf x}(t)\} \in \mcl X_{w(t),u(t)}$ where
	\begin{equation}\label{eq:GPDE_domain}
	\mcl X_{w,u} \hspace{-0.5mm}:= \hspace{-0.5mm}\left\lbrace \bmat{x\\\hat{\mbf x}}\hspace{-0.5mm}\in \hspace{-0.5mm}\R^{n_x}\hspace{-0.5mm}\times\hspace{-0.5mm} X_{v} \mid v = C_v x+ D_{vw} w + D_{vu} u\right\rbrace.
	\end{equation}
 
	\begin{definition}[Admissible solution for a GPDE model]\label{defn:GPDE}
	For given $\{x^0,~\hat{\mbf x}^0\}\in \mcl X_{w(0),u(0)}$ and $w\in L_{2e}^{n_w}[\R_+]$, $u\in L_{2e}^{n_u}[\R_+]$ with $B_vD_{vw} w\in W_{1e}^{2n_S}[\R_+]$ and $B_vD_{vu}u\in W_{1e}^{2n_S}[\R_+]$, we say that $\{x,$ $\hat{\mbf x},$ $z,$ $y,$ $v,$ $r\}$ satisfies the GPDE defined by $\{n,$ $\mbf G_{\mathrm{o}},$ $\mbf G_{\mathrm b},$ $\mbf G_{\mathrm{p}}\}$ (defined in Eqs.~\ref{eq:ode-params},~\ref{eq:BC-parms}, and~\ref{eq:pde-params}) with initial condition $\{x^0,~\hat{\mbf x}^0\}$ and input $\{w,~u\}$ if $z\in L_{2e}^{n_z}[\R_+]$, $y\in L_{2e}^{n_y}[\R_+]$, $v\in L_{2e}^{n_v}[\R_+]$, $r\in L_{2e}^{n_r}[\R_+]$, $\{x(t),\hat{\mbf x}(t)\} \in \mcl X_{w(t),u(t)}$ for all $t\ge 0$, $x$ is differentiable  on $\R_+$, $\hat{\mbf x}$ is Frech\'et differentiable with respect to the $L_2$-norm on $\R_+$, $x(0)=x^0$, $\hat{\mbf x}(0) = \hat{\mbf x}^0$, and Eqs.~\eqref{eq:odepde_general}-\eqref{eq:general_pde_subsystem} are satisfied for all $t\ge 0$.
	\end{definition}

	\subsection{Illustrative Example of the GPDE Representation}
	In this subsection, we illustrate the process of identifying the GPDE parameters of a given system. 
 Firstly, we divide the system into ODE and PDE subsystems and focus on identifying the continuity constraint for the PDE subsystem -- always the least restrictive constraint necessary for existence of the partial derivatives and boundary values. We then proceed to identify the remaining parameters.

	\textbf{Illustration 1} (Damped Wave equation with delay and motor dynamics)\label{ill:Datko}
    Let us consider a wave equation $\ddot{\eta}(t,s) = \partial_s^2 \eta (t,s)$, defined on the interval $s \in [0,1]$, to which we apply a control using a DC motor and where the output from the DC motor experiences a distributed delay and some disturbance modeled as $\eta_s(t,1)=w(t)+\int_{-\tau}^0\mu (s/\tau)T(t+s)ds$ where $T(t)$ is the output of the DC motor and $\mu(s)$ is a given multiplier. The delay is represented using a transport equation with distributed state $p(t,s)$ on the interval $[-1,0]$ so that
	\[
	\dot p(t,s)=\frac{1}{\tau}p_s(t,s), p(t,0)=T(t),\eta(t,1)= \int\limits_{-1}^0\mu(s)p(t,s)ds.
	\]
	The DC motor dynamics relate the voltage input, $u(t)$ to the torque $T(t)$ through the current, $i(t)$ as
	\[
	\partial_t i(t)=\frac{-R}{L}i(t)+u(t) \qquad T(t)=K_t i(t).
	\]
	The sensed output is the feedback signal $\eta_t(1,t)$ and the regulated output is a combination of the integral of the displacement and controller effort so that
	\[
	z(t)=\bmat{\int_{0}^1 \eta(t,s)ds\\u(t)},\qquad y(t)=\eta_t(1,t).
	\]
	Firstly, we introduce the change of variables $\zeta_1 = \eta$, $\zeta_2 = \dot{\eta}$, $\zeta_3(t,s) =p(t,s-1)$. A complete list of equations is now $\dot i(t)=\frac{-R}{L}i(t)+u(t)$ and
	\begin{align*}
	&\dot{\zeta}_1(t,s) = \zeta_2(t,s),\; \dot{\zeta}_2(t,s) = \partial_s^2 \zeta_1 (t,s),\notag\\
	&\dot \zeta_3(t,s)=\frac{1}{\tau}\partial_s \zeta_3(t,s),\;\zeta_1(t,0)=0,\;\zeta_3(t,1)=K_t i(t), \\[-1mm]
	&\partial_s \zeta_1(t,1)=w(t)+ \int_0^{1}\mu(s-1)\zeta_3(t,s)ds,\notag\\
	&z(t)= \bmat{\int_0^1 \zeta_1(t,s) ds\\u(t)}, \; y(t) = \zeta_2(t,1),\; s\in[0,1], \quad t\ge 0.
	\end{align*}
	\noindent\textbf{ODE Subsystem:} Since $i(t)$ is the only finite dimensional state we set $x(t) = i(t)$ to get $\dot x(t)=\frac{-R}{L}x(t)+u(t)$. The ODE subsystem influences the PDE subsystem via signals $w(t)$ and $T(t)$. The effect of the PDE subsystem on the regulated and observed outputs ($z$ and $y$, respectively) is routed through $r(t)$. The outputs, $z,y$ and internal signals, $v,r$, are now defined as
	\begin{align*}
	&v(t)=\bmat{T(t)\\w(t)} = \bmat{K_t\\0} x(t) +\bmat{0\\1}w(t), \\[-1mm]
	&r(t)=\bmat{\int_0^1 \zeta_1(t,s) ds\\ \zeta_2(t,1)},\quad \bmat{z(t)\\y(t)}=\bmat{0\\u(t)\\0}+ \bmat{1&0\\0&0\\0&1}r(t).
\end{align*}
Expressing these equations in the form of Eq.~\eqref{eq:odepde_general}, we find
$\mbf G_{\mathrm o}$ has the following nonzero parameters: $A = \frac{-R}{L}$, $B_{xu} = 1$, $D_{yr} = \bmat{0&\hspace{-1mm}1}$,
\begin{align*}
&\,D_{zu} = \bmat{0\\1}\hspace{-.5mm},\,  \, C_v = \bmat{K_t\\0}\hspace{-.5mm}, \, D_{vw} = \bmat{0\\1}\hspace{-.5mm},\, D_{zr} = \bmat{1&\hspace{-1mm}0\\0&\hspace{-1mm}0}.
\end{align*}
\noindent\textbf{PDE subsystem:} Next, we need to define: $n$, $\mbf G_{\mathrm b}$, and $\mbf G_{\mathrm p}$.

\noindent\textbf{Continuity Constraint:} To identify, $n$, we consider the required partial derivatives and limit values for the three distributed states: $\zeta_1$, $\zeta_2$ and $\zeta_3$. Since $\partial_s^2\zeta_{1}$ appears in the in-domain dynamics and the BCs involve $\zeta_1(t,0)$ and $\partial_s\zeta_1(t,1)$, the least restrictive constraint which guarantees existence of all three terms is $\zeta_1 \in \hat{\mbf x}_2$. Next, no partial derivatives of $\zeta_2$ are needed, but the limit value $\zeta(t,1)$ appears in the BCs -- so we restrict $\zeta_2 \in \hat{\mbf x}_1$. Finally, $\partial_s\zeta_3$ appears in the in-domain dynamics and $\zeta_3(t,1)$ appears in the BCs -- implying $\zeta_3 \in \hat{\mbf x}_1$. We conclude that \textbf{$n=\{n_0,n_1,n_2\}=\{0,2,1\}$} and the GPDE state is $\hat{\mbf x}=col(\hat{\mbf x}_1,\hat{\mbf x}_2)$ where $\hat{\mbf x}_1 = col(\zeta_2,\zeta_3)$ and $\hat{\mbf x}_2 = \zeta_1$.

\noindent\textbf{Boundary Conditions:} For this, $n$, we have $n_{\hat{\mbf x}}=3, n_{S_0}=3, n_{S_1}=3,n_{S_2}=1$ and $n_{S}=4$. Recalling the definitions of $\mcl F, \mcl C,$ and $\mcl B$, we have
\begin{align*}
&\mcl F \hat{\mbf x}=\text{col}(\zeta_2,\zeta_3,\zeta_1,\partial_s\zeta_2,\partial_s \zeta_3,\partial_s \zeta_1,\partial_s^2\zeta_{1}),\\
&\mcl C \hat{\mbf x}=\text{col}(\zeta_2,\zeta_3,\zeta_1,\partial_s\zeta_1),\quad \mcl B \hat{\mbf x}=\text{col}(\mcl C \hat{\mbf x}(0), \mcl C \hat{\mbf x}(1)).
\end{align*}
Checking our BCs, we note that $\zeta_1(t,0)=0$ implies $\zeta_2(t,0)=0$. Placing these BCs in the required form, we have
\[
\int_{0}^1\bmat{
0\\
0\\
0\\
\mu(s-1)\zeta_3(s)
}ds
=\bmat{
\zeta_1(0)\\
\zeta_2(0)\\
\zeta_3(1)\\
\partial_s\zeta_1(1)
}+\bmat{
0\\
0\\
-v_1\\
-v_2}.
\]
Recalling the expansions of $\mcl F \hat{\mbf x}$ and $\mcl B \hat{\mbf x}$ and Eq.~\eqref{eq:odepde_b}, we may identify the parameters in $\mbf G_{\mathrm b}$ as
\begin{align}
&B = \bmat{0&0&1&0_{1,2}&0&0&0\\1&0&0&0_{1,2}&0&0&0\\0&0&0&0_{1,2}&1&0&0\\0&0&0&0_{1,2}&0&0&1}, \notag\\
&B_v =\bmat{0_2\\I_2},\;B_{I}(s) = \bmat{0_{3,1}&\hspace{-2mm}0_{3,1}&\hspace{-2mm}0_{3,5}\\0&\hspace{-2mm}\mu(s-1)&\hspace{-2mm}0_{1,5}}.\label{eq:illus_wave_G_b}
\end{align}
\textbf{In-Domain Dynamics:} Recalling the form of the PDE dynamics from Eq.~\eqref{eq:general_pde_subsystem},
we can represent the dynamics as
\begin{align*}
\dot{\hat{\mbf x}}(t,s)
\hspace{-.5mm}=\hspace{-.5mm}\bmat{\partial_s^2 \zeta_1(t,s)\\1/\tau \partial_s \zeta_3(t,s)\\\zeta_2(t,s)}\hspace{-.5mm}= \hspace{-.75mm} \underbrace{\bmat{0&0_{1,3}&0&0&1\\0&0_{1,3}&\frac{1}{\tau}&0&0\\1&0_{1,3}&0&0&0}}_{A_0}(\mcl F\hat{\mbf x}(t))(s).
\end{align*}
Likewise, from the definition of $r(t)$, we have
\begin{align*}
r(t)=\bmat{\int_0^1 \zeta_1(t,s) ds\\ \zeta_2(t,1)}&= \int\limits_{0}^1\overbrace{\bmat{0_{1,2}&1&0_{1,4}\\0_{1,2}&0&0_{1,4}}}^{C_r}(\mcl F\hat{\mbf x}(t))(\theta)d\theta\\
& \quad +\underbrace{\bmat{0_{1,4}&0&0_{1,3}\\0_{1,4}&1&0_{1,3}}}_{D_{rb}}\mcl B\hat{\mbf x}(t).
\end{align*}
Thus we have $A_0,C_r,D_{rb}$ -- the only nonzero terms in $\mbf G_{\mathrm p}$.
	\section{Representing a PDE Subsystem as a PIE}\label{sec:pde2pie}
	In \Cref{sec:ODEPDE}, we proposed a GPDE representation for a broad class of coupled ODE-PDEs Systems. Now, we beign the process of finding an alternative representation of such a GPDE model as a PIE by focusing on conversion of the PDE subsystem to a restricted class of PIE subsystem of the form
	\begin{equation}\label{eq:PIE_subsystem}			
	\bmat{\hat{\mcl{T}}\dot{\hat{\ubar{\mbf x}}}(t)\\r(t)}=\bmat{\hat{\mcl{A}}&\mcl{B}_v\\\mcl{C}_r&\mcl{D}_{rv}}\bmat{\ubar{\hat{\mbf x}}(t)\\v(t)}-\bmat{\mcl{T}_{v}\dot{v}(t)\\0},
	\end{equation}
	with initial condition $\ubar{\hat{\mbf{x}}}(0) = \ubar{\hat{\mbf{x}}}^0\in L_2^{m}$. Such PIE subsystems are a special case of \Cref{def:piesolution} with parameter set given by
	\[
	\mbf G_{\mathrm{PIE}_s} := \left\lbrace\hat{\mcl T}, \mcl T_v, \emptyset, \hat{\mcl A}, \mcl B_v, \emptyset, \mcl C_{rv}, \emptyset, \mcl D_{rv}, \emptyset,\emptyset,\emptyset\right\rbrace.
	\]
	In this section, we will show that for any PIE-compatible PDE subsystem, there exists a corresponding PIE subsystem such that the existence of a solution for one ensures the existence of a solution for the other.

	\subsection{PIE-compatibility}\label{subsec:assumptions}
	Before we map the PDE subsystem to an associated PIE subsystem, we first define a notion of `PIE-compatibility'. This definition imposes a notion of well-posedness on $X_v$ in the sense that it guarantees the existence of a unitary map to the solution of a PIE. However, it does not guarantee the existence of a solution itself. This condition ensures, e.g., that there are a sufficient number of independent BCs to establish a mapping between the distributed state and its partial derivatives. Without such a mapping, the solution to the PDE may not exist (too many BCs) or may not be unique (too few BCs).
	\begin{definition}[PIE-compatible] \label{as:invertible}
	Given a GPDE parameter set $\{n\in \N^{N+1}$, $\mbf G_{\mathrm b}:=\{B,\, B_{I},\, B_v\}$, $\mbf G_{\mathrm o}$, $\mbf G_{\mathrm p}\}$, we say the \textbf{GPDE is PIE-compatible}, or alternatively \textbf{$\{n,\mbf G_{\mathrm b}\}$ is PIE-compatible}, if $B_T$ is invertible where
	\[
	B_T := B\bmat{T(0)\\T(b-a)}-\int_{a}^{b}B_{I}(s)U_2 T(s-a)ds\in \R^{n_{BC}\times n_S},
	\]	
	and where $T$ and $U_2$ are defined (See also~\Cref{fig:Gb_definitions}) as
	\begin{align}& U_{2i} = \bmat{0_{n_i\times n_{i+1:N}}\\ I_{n_{i+1:N}}} \in \R^{n_{S_{i}} \times n_{S_{i+1}}},\\
	&U_2 = \bmat{\text{diag}(U_{20},\cdots,U_{2(N-1)})\\0_{n_N\times n_S}}\in \R^{ \left(n_{\hat{\mbf x}}+n_S\right)\times n_S},\label{eq:U1}\\
	&T_{i,j}(s) = \frac{s^{(j-i)}}{(j-i)!} \bmat{0_{n_{Si}-n_{Sj}\times n_{Sj}}\\I_{n_{Sj}}}\in \R^{n_{Si} \times n_{Sj}},\label{eq:Tij}\\
	& T(s) = \bmat{T_{1,1}(s)&T_{1,2}(s)&\cdots&T_{1,N}(s)\\0&T_{2,2}(s)&\cdots&T_{2,N}(s)\\\vdots&\vdots&\ddots&\vdots\\0&0&\cdots&T_{N,N}(s)} \in  \R^{n_S \times n_S}.\label{eq:Ts}
	\end{align}
	\end{definition}
	Invertibility of $B_T$ implies the matrix is square requiring $n_{BC}=n_S$. This requirement can be explained by relating to the loss of boundary value information during the differentiation of a function. If a function is differentiated $n_S$-times, we need $n_S$ BCs to relate all the partial derivatives to the original function. Furthermore, $n_{BC}=n_S$ is a necessary but not sufficient requirement -- the BCs must be both independent and provide enough information to allow us to find the original function. See Subsection 3.2.2 in~\cite{peet_2020aut} for an enumeration of pathological cases, including periodic BCs. 


\subsection{A map between PIE and PDE states}
Given a PIE-compatible GPDE, we may construct a PIE subsystem, which we will associate with the PDE subsystem defined by those parameters. The first step is to map $\hat{\mbf x}(t) \in X_v$, the state of the PDE subsystem, to $\hat{\ubar{\mbf x}}(t)\in L_2$, the state of the PIE subsystem using
\begin{align*}	
&\hat{\ubar{\mbf x}} =
\mcl{D}\hat{\mbf{x}}=\bmat{\hat{\mbf x}_0\\\partial_s\hat{\mbf x}_1\\\vdots\\\partial_s^N\hat{\mbf x}_N} \in L_2^{n_{\hat{\mbf x}}},
\end{align*}
where $\mcl D := \text{diag}(\partial_s^0 I_{n_0},\,\cdots,\, \partial_s^N I_{n_N})$. The following theorem shows that this mapping is invertible and the inverse is defined by PI operators.

\begin{restatable}{theorem}{Tmap}\label{thm:T_map}
Given $\{n\in \N^{N+1},$ $\mbf G_{\mathrm b}\}$  PIE-compatible, let $\{\hat{\mcl{T}},$ $\mcl{T}_{v}\}$ be as defined in \cref{fig:Gb_definitions}, $X_v$ as defined in Eq.~\eqref{eq:odepde_general_domain} and $\mcl D$ $:=$diag$(\partial_s^0 I_{n_0},$ $\cdots,$ $\partial_s^N I_{n_N})$. Then we have the following: (a) For any $v\in \R^{n_v}$, if $\hat{\mbf x} \in X_v$, then  $\mcl D \hat{\mbf x} \in L_2^{n_{\hat{\mbf x}}}$ and $\hat{\mbf{x}} = \hat{\mcl{T}}\mcl D \hat{\mbf x}+\mcl{T}_{v}v$; and (b) For any $v\in \R^{n_v}$ and $\ubar{\hat{\mbf x}} \in L_2^{n_{\hat{\mbf x}}}$, $\hat{\mcl{T}}  \ubar{\hat{\mbf{x}}}+\mcl{T}_{v}v \in X_v$ and $\ubar{\hat{\mbf{x}}} = \mcl D (\hat{\mcl{T}}  \ubar{\hat{\mbf{x}}}+\mcl{T}_{v}v)$.
\end{restatable}

First, we generalize Cauchy's classical rule (can also be viewed as a form of Taylor's expansion) for repeated integration to include dependence on boundary values.
\begin{restatable}{lemma}{CauchyLem}\label{lem:cauchy}
Suppose $\mbf x\in W_N^n[a,b]$ for any $N \in \N$. Then
\begin{equation*}
\mbf x(s) \hspace{-.5mm}= \mbf x(a) + \sum\limits_{j=1}^{N-1}\frac{(s\hspace{-.5mm}-\hspace{-.5mm}a)^j}{j!}\partial_s^j\mbf x(a)+\int\limits_{a}^{s}\hspace{-.5mm}\frac{(s\hspace{-.5mm}-\hspace{-.5mm}\theta)^{N-1}}{(N-1)!}\partial_s^N\mbf x(\theta)d\theta.\vspace{-2mm}
\end{equation*}
\end{restatable}
\begin{proof}
This result is proved by induction. Assume the identity holds is true for some $k<N-1$. We can substitute $\partial_s^k \mbf x(s)$ with $\partial_s^k\mbf x(a) +\int\limits_{a}^{s}\hspace{-.5mm}\frac{(s\hspace{-.5mm}-\hspace{-.5mm}\theta)^{k+1}}{(k+1)!}\partial_s^{k+1}\mbf x(\theta)d\theta$ because of the Fundamental Theorem of Calculus. Then, by changing the order of integration, we obtain the identity presented in the Lemma statement. A detailed proof can be found in Appendix~\ref{app:equivalence}.
\end{proof}
The following proof of Thm.~\ref{thm:T_map} applies Lemma~\ref{lem:cauchy} to the PDE model and uses the admissibility criterion to eliminate boundary conditions.
\begin{proof}


\noindent\textbf{Proof of Thm.~\ref{thm:T_map}: Statement a)} Let $\hat{\mbf x}\in X_v$ for some $v\in \R^{q}$. Clearly, by definition of $X_v$, $\partial_s^i\hat{\mbf x}_i\in L_2^{n_i}$. Therefore, $\mcl D\hat{\mbf x}\in L_2^{n_{\hat{\mbf x}}}$. Next we need to express $\hat{\mbf x}$ in terms of $\ubar{\hat{\mbf x}}:=\mcl D\hat{\mbf x}$ and $v$. For that, we will first express $(\mcl C\hat{\mbf x})(a)$ solely in terms of $\ubar{\hat{\mbf x}}$ and $v$. Using Lem.~\ref{lem:cauchy}, we can show that if $\{n,\mbf G_{\mathrm b}\}$ is PIE compatible, then
\[
(\mcl C\hat{\mbf x})(a)= \int_{a}^{b}B_Q(\theta)\ubar{\hat{\mbf x}}(\theta)d\theta+B_T^{-1}B_v v.
\]
This is obtained by applying the $\mcl C$ operator on a $\hat{\mbf x}\in X_v$, followed by boundary conditions in $X_v$ to eliminate all evaluations at $s=b$. For details, see Cor.~\ref{cor:x_D} in ~\Cref{app:T_map} of the supplementary material or ~\cite{shivakumar_2022GPDE}.

Now that we have an expression for $(\mcl C\hat{\mbf x})(a)$, we simply apply Lem.~\ref{lem:cauchy} to $\mbf x\in X_v$, and substitute $(\mcl C\hat{\mbf x})(a)$ into the expression to obtain
\begin{align}
\hat{\mbf x}(s)=\bmat{\hat{\mbf x}_0(s)\\\hat{\mbf x}_{1:N}(s)}&=G_0\ubar{\hat{\mbf x}}(s)
+\int_{s}^{b}G_2(s,\theta)\ubar{\hat{\mbf x}}(\theta)d\theta\notag\\
&\quad + \int_{a}^{s}G_1(s,\theta)\ubar{\hat{\mbf x}}(\theta)d\theta+G_v(s) v\notag\\
&= (\hat{\mcl T}\ubar{\hat{\mbf x}})(s) + (\mcl T_v v) (s).&\notag
\end{align}
\textbf{Proof of Thm.~\ref{thm:T_map}: Statement b)} First, we reverse the formula in Lemma~\ref{lem:cauchy} by differentiating to show that $\mcl D(\hat{\mcl{T}}  \ubar{\hat{\mbf{x}}}+\mcl{T}_{v}v)=\ubar{\hat{\mbf{x}}}$. Next, we can use a similar process to find the boundary values of $\mbf y:=\hat{\mcl{T}}  \ubar{\hat{\mbf{x}}}+\mcl{T}_{v}v$ and show that $\mbf y \in X_v$. 
For details, see~\Cref{app:T_map} in the supplementary material, or~\cite{shivakumar_2022GPDE}.
\end{proof}

For any given $v \in \R^{n_v}$, Theorem~\ref{thm:T_map} provides an invertible map between the state of the PIE subsystem, $\ubar{\hat{\mbf x}}(t)\in L_2^{n_{\hat{\mbf x}}}$ and the state of the PDE subsystem, $\hat{\mbf x}(t)\in X_v$. In the following subsection, we apply this mapping to the internal dynamics of the PDE subsystem in order to obtain an equivalent PIE representation of this subsystem.

\begin{figure}[ht]\vspace{-3mm}
\renewcommand{\figurename}{Block}%
\hspace*{-3mm}
\scalebox{0.9}{
\framebox{\begin{minipage}{0.5\textwidth}
\begin{align*}
& n_{\hat{\mbf x}}=\sum_{i=0}^N n_i,\; n_{S_i}=\sum_{j=i}^N n_{j},\; n_{S}=\sum_{i=1}^N n_{S_i} \; n_{i:j}= \sum_{k=i}^j n_k, \\
&\tau_i(s) = \frac{s^i}{i!},\quad T_{ij} = \tau_{(j-i)}\bmat{0_{(n_{Si}-n_{Sj}), n_{Sj}}\\I_{n_{Sj}}},\; U_{1i} = \bmat{I_{n_i} \\ 0_{n_{i+1:N}, n_i}},\\
&T(s)=\bmat{T_1(s)\\\vdots\\T_N(s)} = \bmat{T_{11}(s)&\cdots&T_{1N}(s)\\\vdots&\ddots&\vdots\\0&\cdots&T_{NN}(s)},\; U_{2i} = \bmat{0_{n_i, n_{i+1:N}}\\ I_{n_{i+1:N}}},\\
&U_1 = \text{diag}(U_{10},\cdots,U_{1N}),\;U_2 = \bmat{\text{diag}(U_{20},\cdots,U_{2(N-1)})\\0_{n_N, n_S}},\\
&B_T = B\bmat{T(0)\\T(b-a)}-\int_{a}^{b}B_{I}(s)U_2 T(s-a)ds,\\
& Q_i(s) =\bmat{0&\hspace{-2mm}\tau_0(s)I_{n_i}&&&\\0&&\hspace{-4mm}\tau_1(s)I_{n_{i+1}}&&\\
&&&\hspace{-4mm}\ddots&\\0&&&&\hspace{-3mm}\tau_{N-i}(s)I_{n_N}}\hspace{-1mm},Q(s)= \bmat{Q_{1}(s)\\\vdots\\Q_N(s)}\\
&B_Q(s)\hspace{-.5mm} = \hspace{-.5mm}B_T^{-1}\left(\hspace{-.5mm}B_I(s)U_1\hspace{-.5mm}+\hspace{-.5mm}\hspace{-.5mm}\int\limits_s^b \hspace{-.5mm}B_I(\theta)U_2Q(\theta-s)d\theta-\hspace{-.5mm}B\bmat{0\\Q(b-s)}\hspace{-.5mm}\right)\hspace{-.5mm},\\
&G_0 = \bmat{I_{n_0}&\\&0_{(n_{\hat{\mbf x}}-n_0)}},\; G_2(s,\theta) = \bmat{0\\T_1(s-a)B_Q(\theta)},\\
&G_1(s,\theta) = \bmat{0\\Q_1(s-\theta)}+G_2(s,\theta),\; G_v(s) = \bmat{0\\T_1(s-a)B_T^{-1}B_v},\\	
&\hat{\mcl{T}} = \fourpi{\emptyset}{\emptyset}{\emptyset}{G_i}, \qquad \mcl{T}_{v}=\fourpi{\emptyset}{\emptyset}{G_v}{\emptyset}.
\end{align*}
\end{minipage}}}
\caption{Definitions based on $n \in \N^{N+1}$ and the parameters of $\mbf G_{\mathrm b}:=\{B$, $B_{I}$, $B_v\}$ used in~\Cref{thm:T_map}.}\label[EqnBlock]{fig:Gb_definitions}
\end{figure}

\subsection{PIE representation of a PDE Subsystem}\label{subsec:pde_equivalence}

For finite-dimensional state-space systems, similarity transforms are used to construct equivalent representations of the input-output map. Specifically, for any invertible $T$,  the system $G:=\{A,B,C,D\}$ with internal state $x$ may be equivalently represented as $G:=\{T^{-1}AT,T^{-1}B,CT,D\}$ with internal state $\hat x=T^{-1}x$. In this subsection, we apply this approach to PDE subsystems. Specifically, now that we have obtained an invertible transformation from $L_2^{n_{\hat{\mbf x}}}$ to $X_v$, we apply the logic of the similarity transform to the internal dynamics of the PDE subsystem to obtain an equivalent PIE subsystem representation. Specifically, in Theorem~\ref{thm:equivalence}, we substitute $\hat{\mbf x} =\hat{\mcl T} \ubar{\hat{\mbf x}}+\mcl T_v v$ in the internal dynamics of the PDE subsystem. The result is a set of equations parameterized entirely using PI operators. These PI operators, as defined in \cref{fig:PIE_subsystem_equation}, specify a PIE subsystem whose input-output behavior mirrors that of the PDE subsystem and whose solution can be constructed using the solution of the PDE subsystem. Conversely, any solution of the associated PIE subsystem can be used to construct a solution for the PDE subsystem.

\begin{restatable}{theorem}{thmEquivalence}\label{thm:equivalence}
Given a PIE-compatible GPDE with parameter set $\{n\in \N^{N+1},$ $\mbf G_{\mathrm b},$ $\mbf G_{\mathrm o},$ $\mbf G_{\mathrm p}\}$ as defined in \cref{eq:BC-parms,eq:ode-params,eq:pde-params}, suppose $v\in L_{2e}^{n_v}[\R_+]$ with $B_v v\in W_{1e}^{2n_S}[\R_+]$,  $\{\hat{\mcl{T}},~\mcl{T}_{v}\}$ are as defined in \cref{fig:Gb_definitions}, $\{\hat{\mcl A},\,\mcl{B}_{v},\,\mcl{C}_r,\,\mcl{D}_{rv}\}$ are as defined in \cref{fig:PIE_subsystem_equation} and\vspace{-.5mm}
\[
\mbf G_{\mathrm{PIE}} = \left\{\hat{\mcl T},\mcl T_v, \emptyset, \hat{\mcl A},\mcl B_v,\emptyset,\mcl C_r,\emptyset,\mcl D_{rv},\emptyset,\emptyset,\emptyset\right\}.\vspace{-.5mm}
\]
Then we have the following.\vspace{-1mm}
\begin{enumerate}
	\item For any $\hat{\mbf{x}}^0$ $\in$ $X_{v(0)}$ ($X_v$ is as defined in \Cref{eq:odepde_general_domain}), if $\{\hat{\mbf{x}},r\}$ satisfies the PDE defined by $\{n, \mbf G_{\mathrm b}, \mbf G_{\mathrm{p}}\}$ with initial condition $\hat{\mbf{x}}^0$ and input $v$, then $\{\mcl D\hat{\mbf{x}}, r, \emptyset\}$ satisfies the PIE defined by $\mbf G_{\mathrm{PIE}}$ with initial condition $\mcl D\hat{\mbf{x}}^0 \in L_2^{n_{\hat{\mbf x}}}$ and input $\{v,\emptyset\}$ where $\mcl D\hat{\mbf x} = \text{col}(\partial_s^0\hat{\mbf x}_0,\cdots,\partial_s^N\hat{\mbf x}_N)$.
	\item For any $\ubar{\hat{\mbf{x}}}^0 \in L_2^{n_{\hat{\mbf x}}}$, if $\{\ubar{\hat{\mbf{x}}},r,\emptyset\}$ satisfies the PIE defined by
	$\mbf G_{\mathrm{PIE}}$ for initial condition $\ubar{\hat{\mbf{x}}}^0$ and input $\{v,\emptyset\}$, then $\{\hat{\mcl T}\ubar{\hat{\mbf{x}}}+\mcl T_v v,~r\}$ satisfies the PDE defined by $\{n, \mbf G_{\mathrm b}, \mbf G_{\mathrm{p}}\}$ with initial condition $\hat{\mbf x}^0 = \hat{\mcl T}\ubar{\hat{\mbf x}}^0 + \mcl T_v v(0)$ and input $v$.\vspace{-1mm}
\end{enumerate}
\end{restatable}

\begin{proof}
Here we will prove 1) implies 2). Let $\{\hat{\mbf x},~r\}$ be as stated in 1). Then by Definition~\ref{defn:PDE}: a)	$r\in L_{2e}^{n_r}[\R_+]$; b) $\hat{\mbf x}(t) \in X_{v(t)}$ for all $t\ge 0$; c) $\hat{\mbf x}$ is Frech\'et differentiable with respect to the $L_2$-norm on $\R_+$; d) \Cref{eq:general_pde_subsystem} is satisfied for all $t\ge 0$; and e) $\hat{\mbf x}(0) = \hat{\mbf x}^0$.

Let $\ubar{\hat{\mbf x}} = \mcl{D}\hat{\mbf x}$, $\ubar{\hat{\mbf x}}^0 = \mcl{D}\hat{\mbf x}^0$, $n=n_{\hat{\mbf x}}$ and $m=0$. Our goal is to show that, $\{\ubar{\hat{\mbf x}},r,\emptyset\}$ satisfies the PIE defined by $\mbf G_{PIE}$ for initial condition $\ubar{\hat{\mbf x}}^0$ and input $\{v,\emptyset\}$. Since $v\in L_{2e}^{n_v}[\R_+]$ from the theorem statement and by the definition of $\mcl T_w$, $B_v v\in W_{1e}^{2n_S}[\R_+]$ we have
\[
(\mcl{T}_{w} w)(s) = (\mcl T_v v)(s) = \bmat{0\\T_1(s-a)}B_T^{-1}B_v v\in W_{1e}^{n_{\hat{\mbf x}}}[\R_+].
\]
From \Cref{thm:T_map}a we have that $\hat{\mbf x}(t) \in X_{v(t)}$ implies $\ubar{\hat{\mbf x}}(t)=\mcl{D}\hat{\mbf x}(t) \in RL_2^{0,n_{\hat{\mbf x}}}=L_2^{n_{\hat{\mbf x}}}$ for all $t \ge 0$. Furthermore, from the \Cref{defn:PDE}, $r\in L_{2e}^{n_r}[\R_+]$. From \Cref{thm:T_map}a we have that $\hat{\mbf x}^0 \in X_{v(0)}$ implies $\ubar{\hat{\mbf x}}^0=\mcl{D}\hat{\mbf x}^0 \in RL_2^{0,n_{\hat{\mbf x}}}=L_2^{n_{\hat{\mbf x}}}$. Furthermore, since $\ubar{\hat{\mbf x}}(t) = \mcl{D}\hat{\mbf x}(t)$ for all $t\ge 0$, we have $\ubar{\hat{\mbf x}}(0) = \mcl{D}\hat{\mbf x}(0)=\mcl D \hat{\mbf x}^0 = \ubar{\mbf x}^0$. Since $\hat{\mbf x}$ is Frech\'et differentiable, the limit of $\frac{\hat{\mbf x}(t+h)-\hat{\mbf x}(t)}{h}$ as $h\to0^+$ exists for any $t\ge 0$ when convergence is defined with respect to the $L_2$ norm. This, and the fact that $\mcl T_v v \in W_{1e}^{n_v}$ implies that
$\lim_{h\rightarrow 0^+}\frac{\hat{\mcl T}\ubar{\hat{\mbf x}}(t+h)-\hat{\mcl T}\ubar{\hat{\mbf x}}(t) }{h}$ exists for all $t\ge 0$. Thus, $\hat{\mcl T} \ubar{\hat{\mbf x}}$ is Frech\'et differentiable with respect to $L_2$-norm. Lastly, since $\hat{\mbf x}(t)$ satisfies \eqref{eq:odepde_b}-\eqref{eq:general_pde_subsystem} for all $t\ge 0$, and since $\hat{\mbf x}(t)\in X_{v(t)}$ and $\ubar{\hat{\mbf x}}(t)=\mcl D \hat{\mbf x}(t)$ for all $t\ge 0$, from \Cref{thm:T_map}, we have that
$\hat{\mbf x}(t) = \hat{\mcl{T}}\ubar{\hat{\mbf x}}(t)+\mcl{T}_{v}v(t)$ implies $\dot{\hat{\mbf x}}(t) = \hat{\mcl{T}}\ubar{\dot{\hat{\mbf{x}}}}(t)+\mcl{T}_{v}\dot{v}(t)$.
We can substitute this into Eq.~\ref{eq:general_pde_subsystem2} and re-write Eq.~\ref{eq:general_pde_subsystem2} using the PI operator notation to get the compact relation
\begin{align*}
\bmat{r(t)\\\hat{\mcl{T}}\ubar{\dot{\hat{\mbf{x}}}}(t)+\mcl{T}_{v}\dot{v}(t)} &= \fourpi{\bmat{0&D_{rb}}}{C_r}{\bmat{B_{xv}&B_{xb}}}{A_i}\bmat{\bmat{v(t)\\(\mcl B \hat{\mbf x})(t)}\\ (\mcl F\hat{\mbf x})(t)}.
\end{align*}
Using \Cref{thm:T_map} and above relations, we can show that $\{\ubar{\hat{\mbf x}},r,\emptyset\}$ satisfies the PIE defined by $\mbf G_{\mathrm{PIE}}$ for initial condition $\ubar{\hat{\mbf x}}^0$ and input $\{v,\emptyset\}$.

Likewise, one can prove 2) implies 1) using the Def.~\ref{def:piesolution}, Def.~\ref{defn:PDE}, and Thm.~\ref{thm:T_map} -- details may be found in~\Cref{app:equivalence} of the supplementary material or~\cite{shivakumar_2022GPDE}. 
\end{proof}
The two parts of Theorem~\ref{thm:equivalence} show that the PDE subsystem's well-posedness guarantees the PIE subsystem's well-posedness and that the PIE subsystem's input-output behavior mirrors that of the PDE subsystem and vice versa. 
Because PIEs are potentially easier to numerically analyze, control, and simulate, this result suggests that the tasks of analysis, control, and simulation of a PDE subsystem may be more readily accomplished by performing the desired task on the PIE subsystem and then applying the result to the original PDE subsystem.


\begin{figure}[ht]
\renewcommand{\figurename}{Block}%
\hspace*{-2.5mm}
\scalebox{0.9}{
\framebox{\begin{minipage}{0.5\textwidth}
\begin{align*}
&R_{D,2}(s,\theta) = U_2T(s-a)B_Q(\theta), \\
&R_{D,1}(s,\theta) = R_{D,2}(s,\theta)+U_2Q(s-\theta),\\
&\Upsilon = \left[\footnotesize\begin{array}{c|c}\bmat{I_{n_v}\\B_T^{-1}B_v\\T(b-a)B_T^{-1}B_v}& \bmat{0_{n_r \times n_x}\\B_Q(s)\\T(b-a)B_Q(s)+Q(b-s)}\\\hline U_2T(s-a)B_T^{-1}B_v&\{U_1,R_{D,1},R_{D,2}\}\end{array}\right],\\
&\Xi = \left[\footnotesize\begin{array}{c|c}
\bmat{0&D_{rb}}& C_r\\\hline \bmat{B_{xv}&B_{xb}}&\{A_i\}
\end{array}\right],\quad \left[\begin{array}{c|c} D_{rv}& C_{rx}\\\hline  B_{v}& \{\hat{A}_i\}\end{array}\right]= \mbf P^4_\times\left(\Xi, \Upsilon \right),\\
&\hat{\mcl A} = \fourpi{\emptyset}{\emptyset}{\emptyset}{\hat A_i}, \; \mcl B_v=\fourpi{\emptyset}{\emptyset}{ B_{v}}{\emptyset},\; \mcl C_{r} =\fourpi{\emptyset}{ C_{rx}}{\emptyset}{\emptyset},\\
&\mcl D_{rv} =\fourpi{ D_{rv}}{\emptyset}{\emptyset}{\emptyset},\;\mcl T =\bmat{I_{n_x}&0\\G_v C_v&\hat{\mcl T}}, \; \mcl T_w =\bmat{0&0\\G_v D_{vw}&0}, \\
&\mcl T_u =\bmat{0&0\\G_v D_{vu}&0}, ,\;\mcl A =\bmat{A+B_{xr}\mcl D_{rv}C_v&B_{xr}\mcl C_{r}\\\mcl B_vC_v&\hat{\mcl A}}, \\
&\mcl B_1= \bmat{B_{xw}+B_{xr}\mcl D_{rv}D_{vw}\\\mcl B_vD_{vw}},\;\mcl B_2 = \bmat{B_{xu}+B_{xr}\mcl D_{rv}D_{vu}\\\mcl B_vD_{vu}},\\
&\mcl C_1 =\bmat{C_z+D_{zr}\mcl D_{rv}C_v&\hspace{-2mm}D_{zr}\mcl C_r},\; \mcl C_2 = \bmat{C_y+D_{yr}\mcl D_{rv}C_v&\hspace{-2mm}D_{yr}\mcl C_r}, \\
&\mcl D_{11} = D_{zw}+D_{zr}\mcl D_{rv}D_{vw}, \quad \mcl D_{12} = D_{zu}+D_{zr}\mcl D_{rv}D_{vu},\\
&\mcl D_{21} =  D_{yw}+D_{yr}\mcl D_{rv}D_{vw}, \quad\mcl D_{22} = D_{yu}+D_{yr}\mcl D_{rv}D_{vu}.
\end{align*}
\end{minipage}}}
\caption{Definitions based on the PDE and GPDE parameters in $\mbf G_{\mathrm{p}}$$=$ $\{A_{0},$ $A_{1},$ $A_{2},$ $B_{xv},$ $B_{xb},$ $C_{r},$ $D_{rb}\}$ and $\mbf G_{\mathrm{o}}$ $=$ $\{A,$ $B_{xw},$ $B_{xu},$ $B_{xr},$ $C_z,$ $D_{zw},$ $D_{zu},$ $D_{zr},$ $C_y,$ $D_{yw},$ $D_{yu},$ $D_{yr},$ $C_v,$ $D_{vw},$ $D_{vu}\}$, the Definitions from $\mbf G_{\mathrm b}$ as listed in \cref{fig:Gb_definitions} and the map $\mbf P_\times^4$ as defined in \Cref{eq:pi_comp_map}.}
\label[EqnBlock]{fig:PIE_subsystem_equation}
\end{figure}

To summarize, finding the PIE representation of a GPDE system involves $4$ major steps. 
To illustrate these steps we go back to the GPDE model of the entropy evolution introduced in Subsec.~\ref{subsec:intro-contribution}. 

\textbf{Illustration 3} \label{ill:entropy-PIE}
Referring back to the PDE model of entropy change from Subsec.~\ref{subsec:intro-contribution}, we have 
\begin{align*}
&\dot{\eta}(t,s) = \partial_s^2 \eta (t,s),\;\eta(t,0) =\eta(t,1)=-\int_{0}^1 \eta(t,s) ds.
\end{align*}
The GPDE representation of this model is defined by $n=\{0,0,2\}$, $\mbf G_{\mathrm p}=\{A_{0} = \bmat{0&0&1}\}$, and
\begin{align*}
&\mbf G_{\mathrm b}=\left\{B = \bmat{1&0&0&0\\0&0&1&0}, \; B_{I} = -\bmat{1&0&0\\1&0&0}\right\}.
\end{align*}
Using the formulae in~\Cref{fig:Gb_definitions,fig:PIE_subsystem_equation}, we find the PIE subsystem as follows (we neglect interconnection to the ODE subsystem as there are no ODEs, inputs, or outputs). 
\begin{align*}
&U_2 = \bmat{1&0\\0&1\\0&0}, \,  U_1 = \bmat{0\\0\\1},\, T(s) = \bmat{1&s\\0&1}, \, Q(s) = \bmat{s\\1}, \\
&B_T = \bmat{2&1/2\\2&3/2},\quad B_Q(s) = (1-s)\bmat{\frac{s}{4}\\-1}, \quad G_0(s) = 0, \\[-1mm]
&G_1(s,\theta) = G_2(s,\theta)+(s-\theta), \quad G_2(s,\theta) = 3s\frac{(s-1)}{4}.
\end{align*}
The PIE form ($\ubar \eta = \partial_2^2 \eta$) of the entropy PDE is then given by
\[
\int\limits_0^s \left(s^2+\frac{s}{4}-\theta\right) \dot{\ubar \eta}(t,\theta)d\theta
+\int\limits_s^1 \frac{3}{4}(s^2-s)  \dot{\ubar \eta}(t,\theta)d\theta =  \ubar \eta(t,s).\vspace{-2mm}
\]


\section{PIE Representation of a GPDE}\label{sec:gpde2pie}
Having converted the PDE subsystem to a PIE, integration of the ODE dynamics is a simple matter of augmenting the PIE subsystem (\Cref{eq:PIE_subsystem}) with the differential equations that define the ODE (\Cref{eq:odepde_general}), followed by elimination of the interconnection signals $v$ and $r$. The result is an augmented PIE system, as defined in~\Cref{eq:PIE_full} whose parameters are 4-PI operators, as defined in \Cref{fig:Gb_definitions,fig:PIE_subsystem_equation}.
%

\subsection{A map between PIE and GPDE states}
Our first step in constructing the augmented PIE system that will be associated with a given GPDE model is to construct the augmented map from the GPDE state (defined on $\mcl X_{w,u}$) to the associated PIE state (defined on $\R L_{2}^{n_x, n_{\hat{\mbf x}}}$). Specifically, given a GPDE model $\{n,\mbf G_{\mathrm b},\mbf G_{\mathrm o}, \mbf G_{\mathrm p}\}$ with $\{n,\mbf G_{\mathrm b}\}$ PIE compatible and state $\mbf x=\bmat{x\\\hat{\mbf x}}\in \mcl X_{w,u}$, the associated PIE system state is $\ubar{\mbf x} = \bmat{x\\\mcl{D}\hat{\mbf{x}}}\in \R L_2^{n_x,n_{\hat{\mbf x}}}$ where $\mcl D:=\text{diag}(\partial_s^0 I_{n_0},$ $\cdots,$ $\partial_s^N I_{n_N})$. Using this definition,~\Cref{cor:T_map_GPDE} shows that if $\{\mcl{T},$ $\mcl{T}_{w},$ $\mcl T_u\}$ are as defined in \cref{fig:PIE_subsystem_equation}, then the map $\mbf x\rightarrow \ubar{\mbf x} $ can be inverted as $\mbf x=\mcl T{\ubar{\mbf x}} + \mcl T_w w+\mcl T_u u$.
\begin{restatable}[Corollary of \Cref{thm:T_map}]{corollary}{TmapGPDE}\label{cor:T_map_GPDE}
Given $\{n\in \N^{N+1},$ $\mbf G_{\mathrm b}\}$ PIE-compatible,  let $\{\mcl{T},$ $\mcl{T}_{w},$ $\mcl T_u\}$ be as defined in \cref{fig:PIE_subsystem_equation}, $\mcl X_{w,u}$ as defined in Eq.~\eqref{eq:GPDE_domain} and $\mcl D$ $:=$diag$(\partial_s^0 I_{n_0},$ $\cdots,$ $\partial_s^N I_{n_N})$. Then for any $w\in \R^{n_w}$ and $u\in \R^{n_u}$ we have:
\begin{enumerate}[label=(\alph*)]
\item If $\mbf x:=\{x,\hat{\mbf{x}}\} \in \mcl X_{w,u}$, then  $\ubar{\mbf x}:=\{x,\mcl D \hat{\mbf x}\} \in \R L_2^{n_x,n_{\hat{\mbf x}}}$ and $\mbf x= \mcl T\ubar{\mbf x}+\mcl{T}_{w}w + \mcl T_u u$.
\item If $\ubar{\mbf x} \in \R L_2^{n_x, n_{\hat{\mbf x}}}$, then $\mbf x:=\mcl{T}  \ubar{\mbf{x}}+\mcl{T}_{w}w+\mcl T_u u \in \mcl X_{w,u}$ and $\ubar{\mbf{x}} = \bmat{I_{n_x}&0\\0&\mcl D} \mbf x$.
\end{enumerate}
\end{restatable}

\begin{proof}
This follows from the application of Thm.~\ref{thm:T_map} along with the definitions of $\mbf x$, $\ubar{\mbf x}$, and $v$. One can augment the ODE state $x$ to the PDE subsystem state $\hat{\mbf x}$ as well as the PIE subsystem state $\hat{\ubar{\mbf x}}$. Since the interconnection signals and definition of $\mcl X_{w,u}$ are consistent with the assumptions of Thm.~\ref{thm:T_map}, $\hat{\mbf x}$ and $\hat{\ubar{\mbf x}}$ are as stated in Thm.~\ref{thm:T_map}. Using this, we can simply eliminate the interconnection signals $v$ to get direct maps between $\mbf x$, $\ubar{\mbf x}$, $w$, and $u$--- See~\Cref{app:T_map_GPDE} of the supplementary material or~\cite{shivakumar_2022GPDE} for more details.
\end{proof}
Thus, for any given $w, u$, we have an invertible transformation from $\R L_2^{n_x,n_{\hat{\mbf x}}}$ to $\mcl X_{w,u}$.

\subsection{Representation of a GPDE model as a PIE system} In this subsection, we define the PIE system associated with a given GPDE model. This associated PIE system is defined by 4-PI parameters as defined in~\Cref{fig:Gb_definitions,fig:PIE_subsystem_equation}. For convenience, we use $\mbf M: \{n, \mbf G_{\mathrm b}, \mbf G_{\mathrm o}, \mbf G_{\mathrm p}\} \mapsto \{\mcl T, \mcl T_w,\mcl T_u,\mcl A, \mcl B_1,\mcl B_2, \mcl C_1, \mcl C_2,\mcl D_{11},\mcl D_{12},\mcl D_{21},\mcl D_{22}\}$ to represent the several formulae used to map GPDE parameters to PIE parameters.
\begin{definition}\label{defn:param-map}
Given $\{n, \mbf G_{\mathrm b}, \mbf G_{\mathrm o}, \mbf G_{\mathrm p}\}$ where
\begin{align*}
\mbf G_{\mathrm{b}} &= \left\{B,B_{I},B_v\right\},~~ \mbf G_{\mathrm{p}} = \left\{A_{0},A_{1},A_{2},B_{xv},B_{xb},C_{r},D_{rb}\right\}\\
\mbf G_{\mathrm{o}} &= \left\{A, B_{xw}, B_{xu}, B_{xr}, C_z, D_{zw}, D_{zu}, D_{zr}, C_y, D_{yw}, D_{yu},\right. \\
&\quad\left. D_{yr}, C_v, D_{vw}, D_{vu}\right\}
\end{align*}
we say that $\mbf G_{\mathrm{PIE}}=\mbf M (\{n, \mbf G_{\mathrm b}, \mbf G_{\mathrm o}, \mbf G_{\mathrm p}\})$ if $\mbf G_{\mathrm{PIE}}=\{\mcl T, \mcl T_w,\mcl T_u,\mcl A, \mcl B_1,\mcl B_2, \mcl C_1, \mcl C_2,\mcl D_{11},\mcl D_{12},\mcl D_{21},\mcl D_{22}\}$ where $\{\mcl T,$ $\mcl T_w,$ $\mcl T_u,$ $\mcl A,$ $\mcl B_1,$ $\mcl B_2,$ $\mcl C_1,$ $\mcl C_2,$ $\mcl D_{11},$ $\mcl D_{12},$ $\mcl D_{21},$ $\mcl D_{22}\}$ are as defined in \Cref{fig:Gb_definitions,fig:PIE_subsystem_equation}.
\end{definition}

Having specified the PIE system associated with a given GPDE model, we now extend the results of~\Cref{thm:equivalence} to show that the map $\mbf{x} \mapsto \bmat{I&0\\0&\mcl{D}}\mbf{x}$ proposed in~\Cref{cor:T_map_GPDE} maps a solution of a given GPDE model to a solution of the associated PIE system and that the inverse map $\ubar{\mbf x}\mapsto \mcl T\ubar{\mbf x}+\mcl{T}_{w}w + \mcl T_u u$ maps a solution of the associated PIE to a solution of the given GPDE model.

%

\begin{restatable}[Corollary of \Cref{thm:equivalence}]{corollary}{corEquivalenceReverse}\label{thm:equivalence_2_reverse}
Given a PIE-compatible GPDE with parameter set $\{n\in \N^{N+1},$ $\mbf G_{\mathrm{o}},$ $\mbf G_{\mathrm b},$ $\mbf G_{\mathrm{p}}\}$ as defined in \Cref{eq:ode-params,,eq:BC-parms,eq:pde-params}, let $w\in L_{2e}^{n_w}[\R_+]$ with $B_vD_{vw} w\in W_{1e}^{2n_S}[\R_+]$, $u \in L_{2e}^{n_u}[\R_+]$ with $B_vD_{vu}u\in W_{1e}^{2n_S}[\R_+]$. Define
\begin{align*}
\mbf G_{\mathrm{PIE}}&=\{\mcl T, \mcl T_w,\mcl T_u,\mcl A, \mcl B_1,\mcl B_2, \mcl C_1, \mcl C_2,\mcl D_{11},\mcl D_{12},\mcl D_{21},\mcl D_{22}\}\\&=  \mbf M(\{n,\mbf G_{\mathrm b}, \mbf G_{\mathrm o}, \mbf G_{\mathrm p}\}.
\end{align*}
Then we have the following:
\begin{enumerate}
	\item For any $\{x^0,\hat{\mbf{x}}^0\}\in \mcl X_{w(0),u(0)}$ (where $\mcl X_{w,u}$ is as defined in \Cref{eq:GPDE_domain}), if $\{x,$ $\hat{\mbf{x}},$ $z,$ $y,$ $v,$ $r\}$ satisfies the GPDE defined by $\{n,$ $\mbf G_{\mathrm{o}},$ $\mbf G_{\mathrm b},$ $\mbf G_{\mathrm{p}}\}$ with initial condition $\{x^0,\hat{\mbf{x}}^0\}$ and input $\{w,u\}$, then $\left\{\bmat{x\\\mcl D\hat{\mbf x}},\, z,\,y\right\}$ satisfies the PIE defined by $\mbf G_{\mathrm{PIE}}$
	with initial condition $\bmat{x^0\\\mcl D\hat{\mbf x}^0}$ and input $\{w,u\}$ where $\mcl D\hat{\mbf x}$ $=$ $\text{col}$$($$\partial_s^0\hat{\mbf x}_0,$$\cdots,$$\partial_s^N\hat{\mbf x}_N$$)$.
	\item For any $\ubar{\mbf x}^0 \in \R L_2^{n_x, n_{\hat{\mbf x}}}$, if $\{\ubar{\mbf x},z,y\}$ satisfies the PIE defined by $\mbf G_{\mathrm{PIE}}$ with initial condition $\ubar{\mbf x}^0$ and input $\{w,u\}$, then $\{x,\hat{\mbf x},z,y,v,r\}$ satisfies the GPDE defined by $\{n, \mbf G_{\mathrm{o}}, \mbf G_{\mathrm b}, \mbf G_{\mathrm{p}}\}$ with initial condition $\bmat{x^0\\\hat{\mbf x}^0} = \mcl T\ubar{\mbf x}^0+\mcl T_w w(0)+\mcl T_uu(0)$ and input $\{w,u\}$ where \begin{align*}
	\bmat{x(t)\\\hat{\mbf{x}}(t)} &:= \mcl{T}\ubar{\mbf x}(t)+\mcl{T}_{w}w(t)+\mcl T_uu(t),\\
	v(t)&:= C_vx(t)+D_{vw}w(t)+D_{vu}u(t),\\
	r(t) &:= \bmat{0_{n_{\hat{\mbf x}}\times n_x}&\mcl C_{r}}\ubar{\mbf x}(t) + \mcl D_{rv} v(t),
	\end{align*}
	and where $\mcl C_{r}$ and $\mcl D_{rv}$ are as defined in \Cref{fig:PIE_subsystem_equation}.
\end{enumerate}

\end{restatable}
To prove Cor.~\ref{thm:equivalence_2_reverse} one simply uses \Cref{thm:equivalence} to construct a PIE representation of the PDE subsystem and then augments this PIE using the dynamics of the ODE subsystem. For details, we may refer to~\Cref{app:equivalence_2} of the supplementary material or extended Arxiv version of this paper~\cite{shivakumar_2022GPDE}. Several examples of the conversion of GPDE models to PIE systems can be found in~\Cref{sec:numerical}.


	\section{Equivalence of Properties of GPDE and PIE}\label{sec:equivalence}
	We have motivated the construction of PIE representations of GPDE models by stating that many analysis, control, and simulation tasks may be more readily accomplished in the PIE framework. However, this motivation is predicated on the assumption that the results of analysis, control and simulation of a PIE system somehow translate to analysis, control and simulation of the original GPDE model. For simulation, the conversion of a numerical solution of a PIE system to the numerical solution of the GPDE is trivial, as per~\Cref{thm:equivalence_2_reverse} through the mapping $\ubar{\mbf{x}}(t) \mapsto \mcl{T}\ubar{\mbf{x}}(t)+\mcl{T}_{w}w(t)+\mcl T_uu(t)$. In this section, we show that analysis and control of the PIE system may also be translated to the GPDE model. For input-output properties, this translation is trivial. For internal stability and control, additional mathematical structure is required.

%
\subsection{Equivalence of Input-Output Properties}
Because the translation of PIE solution to GPDE solution is limited to the internal state of the PIE (inputs and outputs are unchanged),~\Cref{thm:equivalence_2_reverse} implies that all input-output (I/O) properties of the GPDE model are inherited by the PIE system and vice versa. As a result, we have the following.

\begin{corollary}[Input-Output Properties]\label{thm:input_output}
Given a PIE compatible GPDE with parameter set $\{n\in \N^{N+1},$ $\mbf G_{\mathrm{o}},$ $\mbf G_{\mathrm b},$ $\mbf G_{\mathrm{p}}\}$ as defined in \Cref{eq:ode-params,,eq:BC-parms,eq:pde-params}, let $w\in L_{2e}^{n_w}[\R_+]$ with $B_vD_{vw} w\in W_{1e}^{2n_S}[\R_+]$. Let $\mbf G_{\mathrm{PIE}}=  \mbf M(\{n,\mbf G_{\mathrm b}, \mbf G_{\mathrm o}, \mbf G_{\mathrm p}\}$. Suppose $\{x^0,\hat{\mbf{x}}^0\}=\{0,0\}$. Then the following are equivalent.
\begin{enumerate}
	\item If $\{x,$ $\hat{\mbf{x}},$ $z,$ $y,$ $v,$ $r\}$ satisfies the GPDE defined by $\{n,$ $\mbf G_{\mathrm{o}},$ $\mbf G_{\mathrm b},$ $\mbf G_{\mathrm{p}}\}$ with initial condition $\{0,0\}$ and input $\{w,0\}$, then $\norm{z}_{L_2}\le \gamma \norm{w}_{L_2}$.
	\item If $\{\ubar{\mbf x},~z,~y\}$ satisfies the PIE defined by $\mbf G_{\mathrm{PIE}}$ with initial condition $0$ and input $\{w,0\}$, then $\norm{z}_{L_2}\le \gamma \norm{w}_{L_2}$.
\end{enumerate}
Suppose $\mcl K:\in L_{2e}^{n_y}\rightarrow L_{2e}^{n_u}$. Then the following are equivalent.
\begin{enumerate}[label=\alph*)]
	\item If $\{x,$ $\hat{\mbf{x}},$ $z,$ $y,$ $v,$ $r\}$ satisfies the GPDE defined by $\{n,$ $\mbf G_{\mathrm{o}},$ $\mbf G_{\mathrm b},$ $\mbf G_{\mathrm{p}}\}$ with initial condition $\{0,0\}$ and input $\{w,\mcl Ky\}$, then $\norm{z}_{L_2}\le \gamma \norm{w}_{L_2}$.
	\item If $\{\ubar{\mbf x},~z,~y\}$ satisfies the PIE defined by $\mbf G_{\mathrm{PIE}}$ with initial condition $0$ and input $\{w,\mcl Ky\}$, then $\norm{z}_{L_2}\le \gamma \norm{w}_{L_2}$.
\end{enumerate}
\end{corollary}
\begin{proof}
Clearly, the change in representation only changes the internal state of the solutions of a GPDE and its PIE without affecting the inputs or outputs. Thus, \Cref{thm:input_output} follows directly from~\Cref{thm:equivalence_2_reverse}.
\end{proof}

	\subsection{Equivalence of Internal Stability}\label{sec:stab_pde}
	Unlike I/O properties, the question of internal stability of a GPDE model is complicated by the fact that there is no universally accepted definition of stability for such models. Specifically, if the state-space of a GPDE model is defined to be $\mcl X_{u,w}$ (a subspace of the Sobolev space $W^n$), then the obvious norm is the Sobolev norm -- implying that exponential stability requires exponential decay with respect to the Sobolev norm. However, many results on stability of PDE models use the $L_2$ norm as a simpler notion of size of the state.

	In this section, we show that while both notions of stability are reasonable, the use of the Sobolev norm and associated inner product confers significant advantages in terms of mathematical structure on the GPDE model and offers a clear equivalence between internal stability of the GPDE model and associated PIE system. In particular, we first show that $\mcl X_{0,0}$ is a Hilbert space when equipped with the Sobolev inner product and furthermore, exponential stability of the GPDE model with respect to the Sobolev norm is equivalent to exponential stability of the PIE system with respect to the $L_2$ norm.

	\subsubsection{Topology of $\mcl X_{0,0}$ }
	Before we begin, for $n \in \N^N$, let us recall the standard inner product on $\R^{n_x}\times W^n$
	\begin{align*}
	&\ip{\bmat{u\\\mbf u}}{\bmat{v\\\mbf v}}_{\R^{n_x}\times W^n}=u^T v + \sum\nolimits_{i=0}^N \ip{\mbf u_i}{\mbf v_i}_{W^{n_i}_i},\\
	&\ip{\mbf u_i}{\mbf v_i}_{W^{n_i}_i}: = \sum\nolimits_{j=0}^{i}\ip{\partial_s^j \mbf u_i}{\partial_s^j \mbf u_i}_{L_2}
	\end{align*}
	with associated norms $\norm{\mbf u_i}_{W^{n_i}_i}: = \sum\nolimits_{j=0}^{i}\norm{\partial_s^j \mbf x_i}_{L_2^{n_i}}$ and
	\begin{align*}
	&\norm{\bmat{u\\\mbf u}}_{\R^{n_x}\times W^n} = \norm{u}_{\R}+\sum\nolimits_{i=0}^N \norm{\mbf u_i}_{W_i^{n_i}}.
	\end{align*}
	As we will see, however, the standard inner product $\R^{n_x}\times W^n$ is not quite the right inner product for $\mcl X_{0,0}$. For this reason, we propose a slightly modified inner product which we will denote $\ip{\cdot}{\cdot}_{X^n}$, and show that this new inner product is equivalent to the standard inner product on $W^n$. Specifically, we have
	\begin{align}\label{eq:ip_xv}
	\ip{\mbf{u}}{\mbf{v}}_{X^n}:= \sum\nolimits_{i=0}^N \ip{\partial_s^i\mbf{u}_i}{\partial_s^i\mbf{v}_i}_{L_2^{n_i}}= \ip{\mcl D\mbf u}{\mcl D \mbf v}_{L_2^{n_{\mathrm x}}}
	\end{align}	
	and define the obvious extension
	\[
	\ip{\bmat{u\\\mbf u}}{\bmat{v\\\mbf v}}_{\R^{n_x}\times X^n}:=u^T v +\ip{\mbf{u}}{\mbf{v}}_{X^n}.
	\]

	We now show that the norms $\norm{\cdot}_{\R^{n_x} \times W^n}$ and $\norm{\cdot}_{\R^{n_x} \times X^n}$ are equivalent on the subspace $\mcl X_{0,0}$.

	\begin{restatable}{lemma}{lemNormEquivalence}\label{lem:norm_equivalence}
	Suppose $\{n,\mbf G_{\mathrm b}\}$ is PIE-compatible. Then $\norm{\mbf{u}}_{\R^{n_x} \times X^n}\le \norm{\mbf{u}}_{\R^{n_x}\times W^n}$ and there exists $c_0>0$ such that for any $\mbf{u}\in \mcl X_{0,0}$, we have $\norm{\mbf{u}}_{\R^{n_x}\times W^n}\le c_0\norm{\mbf{u}}_{\R^{n_x} \times X^n}$.
	\end{restatable}
	\begin{proof}
	Clearly, the $\R^{n_x}\times W^n$-norm upper bounds $\R^{n_x} \times X^n$-norm by definition of these norms. Next, the map $\ubar{\mbf x} \rightarrow \mbf x$ is a PI operator and, hence, bounded, which allows a bound on all terms in the Sobolev norm. See~\Cref{app:norm_equivalence} of the supplementary material or \cite{shivakumar_2022GPDE} for a complete proof.
	\end{proof}
	Trivially, using $n_x=0$, this result also extends to equivalence of $\norm{\cdot}_{W^n}$ and $\norm{\cdot}_{X^n}$ on $X_0$.

	Next, we will show that $\hat{\mcl T}$ and $\mcl T$ are isometric when $X_0$ and $\mcl X_{0,0}$ are endowed with the inner products $\ip{\cdot}{\cdot}_{\R^{n_x}\times W^n}$ and $\ip{\cdot}{\cdot}_{\R^{n_x} \times X^n}$, respectively. This implies that these spaces are complete with respect to both $\norm{\cdot}_{\R^{n_x} \times X^n}$ ($\norm{\cdot}_{X^n}$) and $\norm{\cdot}_{\R^{n_x} \times W^n}$ ($\norm{\cdot}_{W^n}$).

	\subsubsection{$\mcl{X}_{0,0}$ is Hilbert and $\mcl T$ is unitary}\label{subsec:unitary}
First, recall $X_0$ and $\mcl X_{0,0}$ are defined by $\{n,\mbf G_{\mathrm b},\mbf G_{\mathrm o},\mbf G_{\mathrm p}\}$ as
\begin{align*}
X_{0} &:= \left\lbrace \mat{\hat{\mbf{x}}\in W^{n}[a,b]: B \mcl B \hat{\mbf x} = \int_{a}^{b} B_{I}(s)(\mcl F\hat{\mbf x})(s)ds}\right\rbrace,\\
\mcl X_{0,0}&:=\left\{\bmat{x \\ \hat{\mbf x}} \in \R \times X_{v}\; : \; v=C_v x\right\}.
\end{align*}
The sets $X_{0}$ and $\mcl X_{0,0}$ are the subspaces of valid PDE subsystem and GPDE model states when $v=0$ and when $u=0,w=0$, respectively. Previously, in \Cref{thm:T_map}, we have shown that $\hat{\mcl T}$ is a bijective map. In~\Cref{thm:unitary_T} we extend this result to show that $\hat{\mcl T}:L_2^{n_{\hat{\mbf x}}} \rightarrow X^n$ and ${\mcl T}:\R L_2^{n_x, n_{\hat{\mbf x}}} \rightarrow \R^{n_x} \times X^n$ are unitary in that the respective inner products are preserved under these transformations.

\begin{restatable}{theorem}{thmUnitary}\label{thm:unitary_T}
Suppose $\{n,\mbf G_{\mathrm b}\}$ is PIE-compatible, $\{\hat{\mcl T},\mcl T_v\}$ are as defined in \Cref{fig:Gb_definitions}, and $\{\mcl T$, $\mcl T_w$, $\mcl T_u \}$ are as defined in \Cref{fig:PIE_subsystem_equation} for some matrices $C_v$, $D_{vw}$ and $D_{vu}$. If $\ip{\cdot}{\cdot}_{X^n}$ is as defined in \Cref{eq:ip_xv}, then we have the following:
\begin{enumerate}[label=\alph*)]
\item for any $v_1, v_2 \in \R^{n_v}$ and $\ubar{\hat{\mbf{x}}},~ \ubar{\hat{\mbf{y}}}\in L_2^{n_{\hat{\mbf x}}}$
\begin{align}
&\ip{\left(\hat{\mcl{T}}\ubar{\hat{\mbf{x}}}+\mcl{T}_v v_1 \right)}{\left(\hat{\mcl{T}}\ubar{\hat{\mbf{y}}}+\mcl{T}_v v_2\right)}_{X^n} = \ip{\ubar{\hat{\mbf{x}}}}{\ubar{\hat{\mbf{y}}}}_{L_2^{n_{\hat{\mbf x}}}}.
\end{align}
\item for any $w_1, w_2 \in \R^{n_w}$, $u_1, u_2 \in \R^{n_u}$, $\ubar{\mbf x}, \ubar{\mbf y}\in \R L_2^{n_x, n_{\hat{\mbf x}}}$,
\begin{align}
&\ip{\left(\mcl{T}\ubar{\mbf x}+\mcl T_w w_1+\mcl T_u u_1\right)}{\left(\mcl{T}\ubar{\mbf y}+\mcl T_w w_2+\mcl T_u u_2\right)}_{\R^{n_x} \times X^n} \notag\\
&\quad= \ip{\ubar{\mbf x}}{\ubar{\mbf y}}_{\R L_2^{n_x, n_{\hat{\mbf x}}}}.
\end{align}
\end{enumerate}
\end{restatable}
\begin{proof}
Let $\ubar{\hat{\mbf{x}}}, \ubar{\hat{\mbf y}}\in L_2^{n_{\hat{\mbf x}}}$  and $v_1, v_2\in \R^{n_v}$. Then, from \Cref{thm:T_map}, we have
\begin{align*}
\hat{\mcl{T}}\ubar{\hat{\mbf x}} +\mcl{T}_{v}v_1 \in X_{v_1},\quad \hat{\mcl{T}}\ubar{\hat{\mbf y}} +\mcl{T}_{v}v_2 \in X_{v_2}.
\end{align*}
Therefore, by \Cref{eq:ip_xv} and the result in \Cref{thm:T_map}b,
\begin{align*}
&\ip{\left(\hat{\mcl{T}}\ubar{\hat{\mbf x}}+\mcl{T}_v v_1 \right)}{\left(\hat{\mcl{T}}\ubar{\hat{\mbf y}}+\mcl{T}_v v_2\right)}_{X^n} \\
&= \ip{\mcl D\left(\hat{\mcl{T}}\ubar{\hat{\mbf x}}+\mcl{T}_v v_1 \right)}{\mcl D\left(\hat{\mcl{T}}\ubar{\hat{\mbf y}}+\mcl{T}_v v_2 \right)}_{L_2^{n_{\hat{\mbf x}}}} = \ip{\ubar{\hat{\mbf x}}}{\ubar{\hat{\mbf y}}}_{L_2^{n_{\hat{\mbf x}}}}.
\end{align*}

For b), let $\ubar{\mbf x}, \ubar{\mbf y} \in \R L_2^{n_x,n_{\hat{\mbf x}}}$ and $w_1, w_2\in \R^{n_w}$, $u_1, u_2\in \R^{n_u}$. Then, from \Cref{cor:T_map_GPDE}, we have $\mcl T\ubar{\mbf x} +\mcl{T}_{w}w_1+\mcl T_u u_1 \in \mcl X_{w_1,u_1}$ and $\mcl T \ubar{\mbf y} +\mcl{T}_{w}w_2 +\mcl T_u u_2\in \mcl X_{w_2, u_2}$.
Since $\R^{n_x}\times X^n$ inner product is just sum of $\R$ and $X^n$ inner products, using definitions of $\mcl T$, $\mcl T_w$, and $\mcl T_u$ and the result in \Cref{cor:T_map_GPDE}b, we have
\begin{align*}
&\ip{\left(\mcl T\ubar{\mbf x} +\mcl{T}_{w}w_1+\mcl T_u u_1\right)}{\left(\mcl T\ubar{\mbf y} +\mcl{T}_{w}w_2+\mcl T_u u_2\right)}_{\R^{n_x} \times X^n} \\
&=\ip{\ubar{\mbf x}}{\ubar{\mbf y}}_{\R L_2^{n_x, n_{\hat{\mbf x}}}}.
\end{align*}
\end{proof}

\begin{restatable}{corollary}{corUnitary}\label{cor:unitary_T}
Suppose $\{n,\mbf G_{\mathrm b}\}$ is PIE-compatible, $\hat{\mcl T}$ is as defined in \Cref{fig:Gb_definitions}, $\mcl T$ is as defined in \Cref{fig:PIE_subsystem_equation}, $X_v$ is as defined in Eq.~\eqref{eq:odepde_general_domain} and, for any matrices $C_v$, $D_{vw}$ and $D_{vu}$, $\mcl X_{w,u}$ is as defined in Eq.~\eqref{eq:GPDE_domain}. Then $X_{0}$ is complete with respect to $\norm{\cdot}_{X^n}$ and $\mcl X_{0,0}$ is complete with respect to $\norm{\cdot}_{\R^{n_x}\times X^n}$. Furthermore, $\hat{\mcl T}:L_2^{n_{\hat{\mbf x}}}\to X_0$ and $\mcl T: \R L_2^{n_x,n_{\hat{\mbf x}}} \to \mcl X_{0,0}$ are unitary (isometric surjective mappings between Hilbert spaces).
\end{restatable}
\begin{proof}
From~\Cref{thm:T_map} and~\Cref{cor:T_map_GPDE}, we have that $\mcl T$ is a bijective mapping from $\R L_2^{n_x, n_{\hat{\mbf x}}}$ to $\mcl X_{0,0}$.  From~\Cref{thm:unitary_T}, we have that $\mcl T$ is isometric with respect to the $\R^{n_x}\times X^n$ inner product. Since $\R L_2^{n_x, n_{\hat{\mbf x}}}$ is complete, we conclude that $\mcl X_{0,0}$ is complete with respect to the $\R^{n_x}\times X^n$ norm. Completeness of $X_0$ follows trivially from the special case $n_x=0$.
%
\end{proof}
As a direct consequence of~\Cref{cor:unitary_T} and~\Cref{lem:norm_equivalence}, $X_0$ and $\mcl X_{0,0}$ are also complete with respect to $\norm{\cdot}_{W^n}$ and $\norm{\cdot}_{\R^{n_x}\times W^n}$, respectively.

\subsubsection{Equivalence of Internal Stability Properties}
Although the natural definition of exponential stability of a GPDE model is with respect to the $\R^{n_x}\times X^n$ norm, since the $\R^{n_x}\times X^n$ norm is equivalent to the $\R^{n_x}\times W^n$ norm, we can also claim stability with respect to the $\R^{n_x}\times W^n$ norm.
\begin{definition}[Exponential Stability of a GPDE model]
We say a GPDE model defined by $\{n, \mbf G_{\mathrm o}, \mbf G_{\mathrm b}, \mbf G_{\mathrm p}\}$ is exponentially stable if there exist constants $M$, $\alpha>0$ such that for any $\{x^0,\hat{\mbf{x}}^0\}\in \mcl X_{0,0}$, if $\{x$, $\hat{\mbf{x}}, z, y, v, r\}$ satisfies the GPDE defined by $\{n, \mbf G_{\mathrm o}, \mbf G_{\mathrm b}, \mbf G_{\mathrm p}\}$ with initial condition $\{x^0,\hat{\mbf x}^0\}$ and input $\{0,0\}$, then
\begin{align*}
\norm{\bmat{x(t)\\\hat{\mbf{x}}(t)}}_{\R^{n_x}\times W^n}\le M\norm{\bmat{x^0\\\hat{\mbf{x}}^0}}_{\R^{n_x}\times W^n}e^{-\alpha t} \quad \text{for all} ~t\ge 0.
\end{align*}
\end{definition}

Clearly, internal stability of a PIE system is with respect to the $\R L_2$ norm.

\begin{definition}[Exponential Stability of a PIE system]
We say a PIE defined by $\mbf G_{\mathrm{PIE}}$ is exponentially stable if there exist $M$, $\alpha>0$ such that for any $\ubar{\mbf{x}}^0\in \R L_2^{n_x, n_{\hat{\mbf x}}}$, if $\{\ubar{\mbf{x}}, z, y\}$ satisfies the PIE defined by $\mbf G_{\mathrm{PIE}}$ with initial condition $\ubar{\mbf x}^0$ and input $\{0,0\}$, then $\norm{\ubar{\mbf{x}}(t)}_{\R L_2}\le M\norm{\ubar{\mbf x}^0}_{\R L_2}e^{-\alpha t}$ for all $t\ge 0$.
\end{definition}

Exponential stability of a GPDE model is equivalent to exponential stability of the associated PIE system.

\begin{restatable}{theorem}{thmStability}\label{thm:stability_equivalence}
Given $\{n,\mbf G_{\mathrm{o}}, \mbf G_{\mathrm b}, \mbf G_{\mathrm{p}}\}$  PIE-compatible, the GPDE model defined by $\{n,\mbf G_{\mathrm{o}}, \mbf G_{\mathrm b}, \mbf G_{\mathrm{p}}\}$ is exponentially stable if and only if the PIE defined by $\mbf G_{\mathrm{PIE}}:=\mbf M(\{n,\mbf G_{\mathrm b},\mbf G_{\mathrm o},\mbf G_{\mathrm p}\})$ is exponentially stable.
\end{restatable}
\begin{proof}
Suppose GPDE defined by $\{n, \mbf G_{\mathrm{o}}, \mbf G_{\mathrm b}, \mbf G_{\mathrm{p}}\}$ is exponentially stable. Then, there exist constants $M$, $\alpha>0$ such that for any $\{x^0,\hat{\mbf{x}}^0\}\in \mcl X_{0,0}$, if $\{x$, $\hat{\mbf{x}},z,y,v,r\}$ satisfies the GPDE defined $\{n, \mbf G_{\mathrm o}, \mbf G_{\mathrm b}, \mbf G_{\mathrm p}\}$ with initial condition $\{x^0,\hat{\mbf x}^0\}$ and input $\{0,0\}$, we have
\begin{align*}
\norm{\bmat{x(t)\\\hat{\mbf x}(t)}}_{\R^{n_x}\times W^n}\le M\norm{\bmat{x^0\\\hat{\mbf{x}}^0}}_{\R^{n_x}\times W^n}e^{-\alpha t} \quad \text{for all} ~t\ge 0.
\end{align*}

For any $\ubar{\mbf x}^0\in \R L_2^{n_x,n_{\hat{\mbf x}}}$, let $\{\ubar{\mbf x},z,y\}$ satisfy the PIE defined by $\mbf G_{\mathrm{PIE}}$ with initial condition $\ubar{\mbf x}^0 \in \R L_2^{n_x, n_{\hat{\mbf x}}}$ and input $\{0,0\}$. Then, from \Cref{thm:equivalence_2_reverse}, $\{x,\hat{\mbf x},z,y,v,r\}$ satisfies the GPDE defined by $\{n, \mbf G_{\mathrm{o}}, \mbf G_{\mathrm b}, \mbf G_{\mathrm{p}}\}$ with initial condition $\bmat{x^0\\\hat{\mbf x}^0} := \mcl T\ubar{\mbf x}^0 \in \mcl X_{0,0}$ and input $\{0,0\}$ for some $v$ and $r$ where $\bmat{x(t)\\\hat{\mbf x}(t)} := \mcl T \ubar{\mbf x}(t)$.
Then, by the exponential stability of the GPDE, we have
\[
\norm{\bmat{x(t)\\\hat{\mbf x}(t)}}_{\R^{n_x}\times W^n}\le M\norm{\bmat{x^0\\\hat{\mbf x}^0}}_{\R^{n_x}\times W^n}e^{-\alpha t} \quad \text{for all} ~t\ge 0.
\]

By \cref{thm:unitary_T,lem:norm_equivalence}, for any $\mbf x\in \R L_2$ we have $\norm{\mbf x}_{\R L_2} = \norm{\mcl T\mbf x}_{\R^{n_x}\times X^n}$. Thus, we have the following:
\begin{align*}
\norm{\ubar{\mbf x}(t)}_{\R L_2} &= \norm{\mcl T\ubar{\mbf x}(t)}_{\R^{n_x}\times X^n} \\
&\le \norm{\bmat{x(t)\\\hat{\mbf x}(t)}}_{\R^{n_x}\times W^n} \le M\norm{\bmat{x^0\\\hat{\mbf x}^0}}_{\R^{n_x}\times W^n} e^{-\alpha t} \\
&\hspace{-5mm}\le c_0M \norm{\bmat{x^0\\\hat{\mbf x}^0}}_{\R^{n_x}\times X^n} e^{-\alpha t} = c_0 M \norm{\ubar{\mbf x}^0}_{\R L_2}e^{-\alpha t}.
\end{align*}
Therefore, the PIE defined by $\mbf G_{\mathrm{PIE}}$ is exponentially stable.

For the reverse implication, we start by assuming the exponential stability of the PIE. Then, we take an arbitrary solution of the GPDE, and find the associated solution of the PIE using \Cref{thm:equivalence_2_reverse,thm:T_map}. Finally, by applying, \cref{lem:norm_equivalence,thm:unitary_T}, we get the exponentially stability of the GPDE.

\end{proof}		

The results of~\Cref{thm:stability_equivalence} also imply that Lyapunov and asymptotic stability of the GPDE model in the $\R^{n_x}\times W^n$ norm are equivalent to Lyapunov and asymptotic stability of the associated PIE system in the $\R L_2$ norm --- See \Cref{app:stability_equivalence} for details. 


	\subsection{Convex Conditions for Internal Stability of a GPDE model}\label{sec:appl}

Next, we look at a test for stability of the PDE in $\R L_2$ norm (as opposed to the $\R \times W^n$ norm) using the PIE framework. While we do not establish stability of the PIE system itself. The stability test is defined in terms of the PIE system representation but is not actually a test for the stability of the PIE system. This stability test is defined in terms of the existence of positive semidefinite PI operators subject to affine inequality constraints. Such forms of convex optimization are labeled Linear PI Inequalities (LPIs), and LMI-based methods for the feasibility of LPIs have been discussed in, e.g.~\cite{peet_2020aut}. 

\begin{theorem}\label{thm:Lyapunov}
Given $\{n,\mbf G_{\mathrm{o}}, \mbf G_{\mathrm b}, \mbf G_{\mathrm{p}}\}$  PIE-compatible, let $\mbf G_{\mathrm{PIE}}:=\mbf M(\{n,\mbf G_{\mathrm b},\mbf G_{\mathrm o},\mbf G_{\mathrm p}\})$. Suppose there exist $\epsilon,\delta>0$,  and $\mcl P \in [\Pi_4]_{n_x,n_x}^{n_{\hat{\mbf x}},n_{\hat{\mbf x}}}$ such that $\mcl P=\mcl P^*\ge \epsilon I$ and
\[
\mcl A^*\mcl P \mcl T+\mcl T^*\mcl P \mcl A \le -\delta \mcl T^*\mcl T.
\] Then the GPDE model defined by $\{n,\mbf G_{\mathrm{o}}, \mbf G_{\mathrm b}, \mbf G_{\mathrm{p}}\}$ is exponentially stable in the $\R L_2^{n_x, n_{\hat{\mbf x}}}$ norm.
\end{theorem}
\begin{proof}
    This can proved by using a Lyapunov function candidate $V(\mbf x) = \ip{\mcl T\ubar{\mbf x}}{\mcl P\mcl T\ubar{\mbf x}}_{\R L_2}$. Then, by the assumptions of the Theorem $\dot{V}(\ubar{\mbf x}(t))\le -2\alpha V(\ubar{\mbf x}(t))$. Thus, one can use the coercivity of $V$ to show that any solution $\mbf x=\mcl T\ubar{\mbf x}$ to the GPDE model is exponentially stable in $\R L_2$-norm.
\end{proof}
Lastly, we note that other LPI tests for properties of the GPDE in terms of the associated PIE include $L_2$-gain~\cite{shivakumar_CDC2019}, $H_\infty$-optimal estimator design~\cite{das_2019ACC}, and $H_\infty$-optimal full-state feedback controller synthesis~\cite{shivakumar_2020CDC}.

\subsection{A Side Note on PIETOOLS}\label{sec:PIETOOLS}
The generality of the class of models (GPDE) considered requires the identification of a large number of system parameters --- most of which are typically zero or sparse. Furthermore, construction of the associated PIE system using the formulae in~\Cref{fig:Gb_definitions,fig:PIE_subsystem_equation} can be cumbersome. This complicated process of identification of parameters and application of formulae may thus be an impediment to the practical application of the results in this paper. For this reason, PIETOOLS versions include a Graphical User Interface and symbolic parser for the construction of GPDE models, which do not require the user to understand the notational system defined in this paper. 
Additional details can be found in the PIETOOLS user manual~\cite{PIETOOLS2021b}.
%
%
%
%
In addition to the GUI, PIETOOLS includes many tools for the analysis, control, estimation and simulation of PIE systems in the context of: simple PDE models, advanced GPDE models and Delay Differential Equations.	
%
%
In the following section, we apply this GUI to several GPDE models. We will also include results generated by the analysis, control, and simulation tools in PIETOOLS when relevant.


	\section{Examples of the PIE representation}\label{sec:numerical}
	Now, we illustrate the PIE representation of four GPDE models. We can use the PIETOOLS GUI, as described in Section~\ref{sec:PIETOOLS}, to construct the associated PIE system. More examples can be found in Appendix~\ref{app:more_examples} or~\cite{shivakumar_2022GPDE,PIETOOLS2021b}.
	\subsection{Damped wave equation with delay and motor dynamics}\label{subsec:wave_example}
	First, we revisit the GPDE model studied in~\Cref{ill:Datko}. Since we have already identified the parameters of the GPDE model and we can construct the associated PIE system using the formulae in~\cref{fig:PIE_subsystem_equation}. For simplicity, we choose $\mu(s) = 1$ which yields the following nonzero PIE system parameters.\vspace{-2mm}

	{\small
	\begin{align*}
	&\mcl T = \fourpiFull{1}{\bmat{0&0&0}}{\bmat{0\\-2\\-2s}}{0_{3}}{\bmat{1&0&0\\0&0&0\\0&-s\theta&-\theta}}{\bmat{0&0&0\\0&-1&0\\0&-s\theta&-s}},\\
	&\mcl A \hspace{-0.5mm}=\hspace{-0.5mm} \fourpiFull{-\frac{5}{2}}{0_{1, 3}}{\bmat{0\\0\\0}}{\bmat{0&\hspace*{-1.5mm}0&\hspace*{-1.5mm}1\\0&\hspace*{-1.5mm}\frac{1}{0.2}&\hspace*{-1.5mm}0\\0&\hspace*{-1.5mm}0&\hspace*{-1.5mm}0}}{0_{3}}{0_{3}},\, \mcl T_w \hspace{-0.5mm}=\hspace{-0.5mm} \fourpi{\emptyset}{\emptyset}{\bmat{0\\0\\s}}{\emptyset},\\
 &\mcl B_2 = \fourpi{1}{\emptyset}{0_{3, 1}}{\emptyset},\; \mcl C_2 = \fourpi{0}{\bmat{1~0~0}}{\emptyset}{\emptyset}\\
	&\mcl C_1 = \fourpi{\bmat{-1\\0}}{\bmat{0&-0.5s&-0.5s^2-s\\0&0&0}}{\emptyset}{\emptyset},\\
	&\mcl D_{12} = \fourpi{\bmat{0\\1}}{\emptyset}{\emptyset}{\emptyset},\;\mcl D_{11} = \fourpi{\bmat{0.5\\0}}{\emptyset}{\emptyset}{\emptyset}.
	\end{align*}}
\begin{figure}
    \centering
    \includegraphics[width=0.4\textwidth]{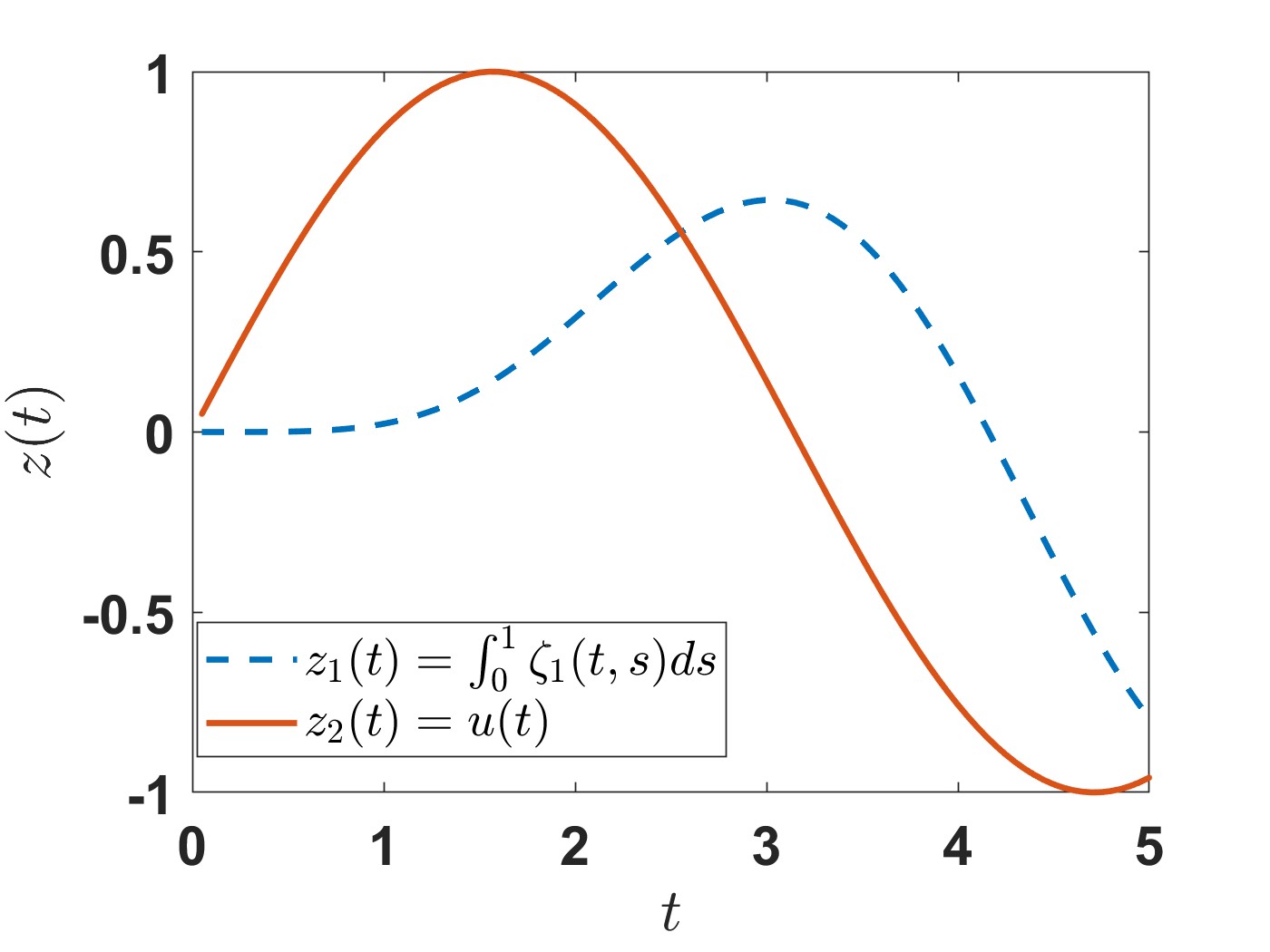}
    \caption{Output response of the wave equation driven by a motor with delay via boundary control. We see in the plot that, for $u(t)=sin(t)$ and $w=0$, the effect of the input on $z_1$ is delayed.}
    \label{fig:wave_output}
\end{figure}

Having obtained a PIE representation, we simulate the response of the system to control input with delay by using a single mode control excitation $u(t)=sin(t)$ and no exogenous inputs. In \Cref{fig:wave_output}, we see that output $z_1$ lags behind the control input $z_2$, which is expected since the maximum delay value is $\tau = 0.5$ and the history of the control input is set to zero. 

\subsection{A $4^{th}$ order PDE: Timoshenko Beam Equation}\label{sec:timoshenko}
In this example, we find the PIE system associated with a GPDE model with a $4^{th}$ order spatial derivative. While the dynamics of the Timoshenko beam~\cite{timoshenko} are often modeled as two coupled $2^{nd}$ order PDEs, if the beam is elastic, isotropic and homogeneous with constant cross-section then these equations can be combined to obtain a $4^{th}$ order GPDE representation.
\begin{align*}
&\rho A\ddot{w}(t,s) - \left(\rho I+\frac{EI\rho}{\kappa G}\right) \partial_s^2\ddot{w}(t,s)+\frac{\rho^2I}{\kappa G} \ddddot{w}(t,s)\\
&= -EI\partial_s^4w(t,s)+d(t),
\end{align*}
where $\rho$ is the density of the beam material, $A$ is the cross section area, $I$ is the second moment of area, $\kappa$ is the Timoshenko beam constant, $E$ is the elastic modulus, $G$ is the shear modulus, and $d$ is some distributed exogenous disturbance. For simplicity, we take $\rho=A=I=\kappa=G=E=1$. The BCs are given by $w(t,0)=0$, $\partial_sw(t,0)=0$, $\partial_s^2w(t,1)=w(t,1)$, and $\partial_s^3w(t,1)=\partial_sw(t,1)$.

We define the state variables as $w,\dot w, \ddot w$ and $\dddot w$ and based on the generator and BCs, we partition the state variables as $\hat{\mbf x}_0 = \text{col}(\dot{w},\dddot{w})$, $\hat{\mbf x}_1=\ddot{w}$, and $\hat{\mbf x}_4 = w$. The full state is then $\hat{\mbf x} =col(\hat{\mbf x}_0,\, \hat{\mbf x}_2,\, \hat{\mbf x}_4)$ so that the continuity condition is $n=\{2,0,1,0,1\}$, implying $n_{S_1}=2$, $n_{S_2}=2$, $n_{S_3}=1$, $n_{S_4}=1$, and hence $n_{S}=6$. Because we require $n_{BC}=n_{S}$, we need two additional BCs. To get these new BCs, we differentiate $w(t,0)=0$ and $\partial_sw(t,0)=0$ twice in time to obtain $\ddot{w}(t,0)=0$ and $\partial_s\ddot{w}(t,0)=0$. We now use the PIETOOLS GUI to calculate the PIE representation as
\begin{align}\label{eqn:Timoshenko_beam}
\fourpi{\emptyset}{\emptyset}{\emptyset}{G_i} \dot{\ubar{\hat{\mbf x}}} \hspace*{1.75mm} (t,s)&= \fourpi{\emptyset}{\emptyset}{\emptyset}{A_i}\ubar{\hat{\mbf x}}(t,s)\\&\qquad+\fourpi{\emptyset}{\emptyset}{\bmat{0~0~0~1}^T}{\emptyset}d(t)\notag
\end{align}
where $\ubar{\hat{\mbf x}} = \text{col}(\dot{w},\dddot{w},\partial_s^2\ddot{w},\partial_s^4w)$ and
\begin{align*}
&A_0(s) = \bmat{0&0&0&0\\0&0&1&-1\\0&1&0&0\\1&0&0&0},A_1(s,\theta) = \bmat{0&0&s-\theta&0\\0&0&\theta-s&0\\0&0&0&0\\0&0&0&0},\\
&A_2(s,\theta) = 0_4,f_0(s,\theta) = -\frac{1}{39}s^3\theta^3 + \frac{s^2\theta^2}{26}\left(3s - \theta   - 2\right),\\
&G_0(s) = \bmat{I_{2}\\&0_{2}}, G_2(s,\theta) = \bmat{0_{3}&\\&\hspace*{-2mm}f_0(s,\theta)- \frac{1}{6}s^2(s + 3\theta)},\\
&G_1(s,\theta) = \bmat{0_{2}&\\&\bmat{s-\theta&0\\0&f_0(s,\theta)+ \frac{1}{6}\theta^2(3s - \theta)}}.
\end{align*}

\begin{figure}
    \centering
    \includegraphics[width=0.4\textwidth]{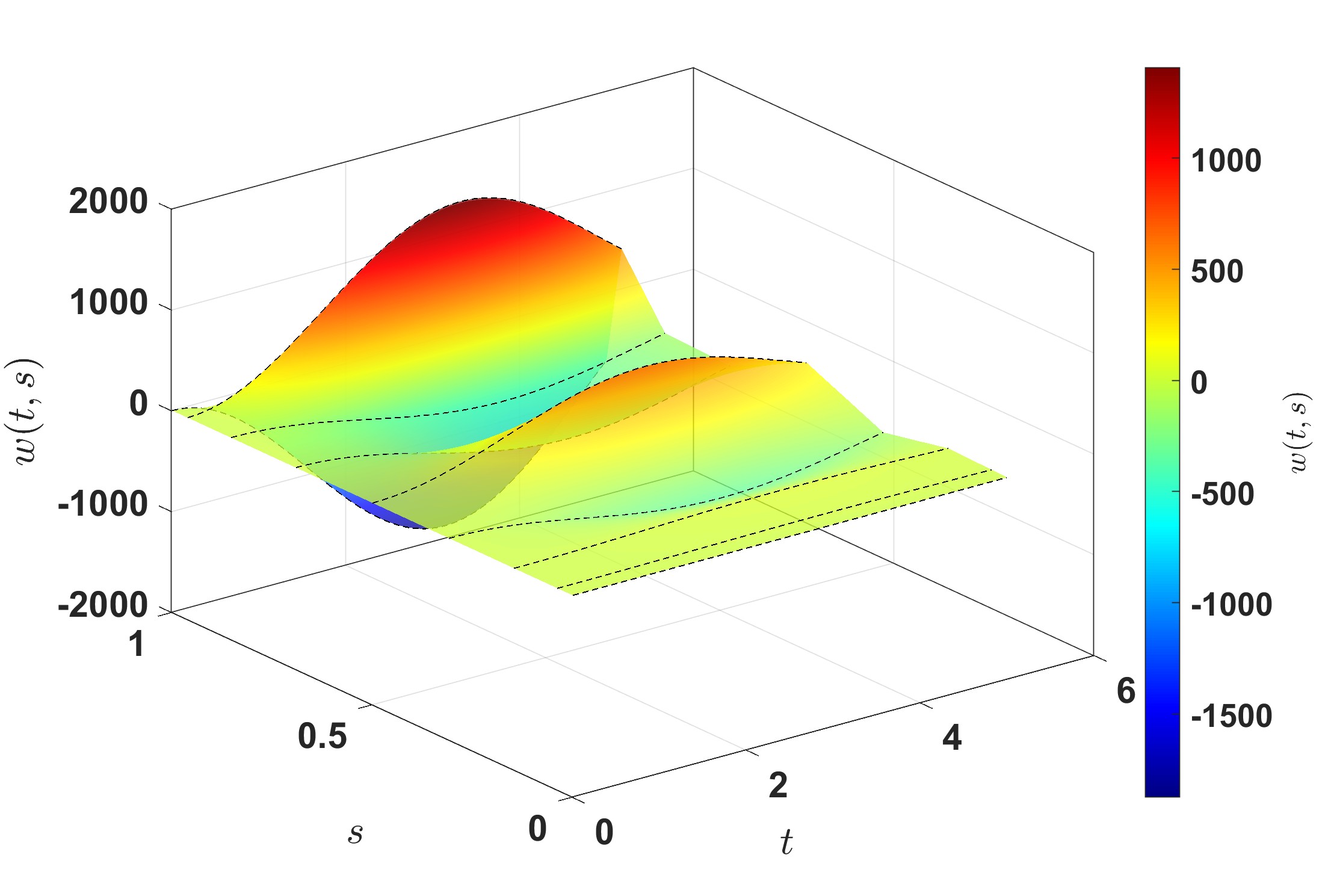}
    \caption{In this figure, we plot the lateral displacement of the Timoshenko Beam equations Eq.~\eqref{eqn:Timoshenko_beam} under a sinosuidal excitation $d(t)=sin(t)$ using the PIE representation given in \Cref{sec:timoshenko}.}
    \label{fig:timo_beam}
\end{figure}

Using the PIE representation given above, we can simulate the system response under external excitation $d$ for any initial conditions. In \Cref{fig:timo_beam}, we simulate the PIE system using PIESIM module of the PIETOOLS with zero initial conditions and a sinusoidal input $d(t)=sin(t)$.  

\subsection{Simulation of diffusion PDE with nonlocal boundary conditions}\label{subsec:entropy}
In this example, we consider a simple linear PDE model of heat-conduction in an elastic material  \cite{day2013heat} which is modeled as a reaction-diffusion equation with some non-local effects that are modeled as volumetric constraints at the boundary. Normalizing the domain, we obtain the GPDE model
\begin{align}\label{eqn:example3_PDE}
\dot{\mbf x}(t,s) &= \mbf x(t,s)+\partial_s^2 \mbf x(t,s)+w_a(t),\quad \mbf x(0,s)=0,\notag\\
\mbf x(t,0) &= \int_0^1 \mbf x(t,s) ds + w_b(t),\qquad \mbf x(t,1)=0,\notag\\
w_b(t) &= \frac{\sin(t)}{t},\quad w_a(t) = cos(2t)e^{-t}.
\end{align}

\begin{figure}[t]
\centering
\includegraphics*[width=0.4\textwidth]{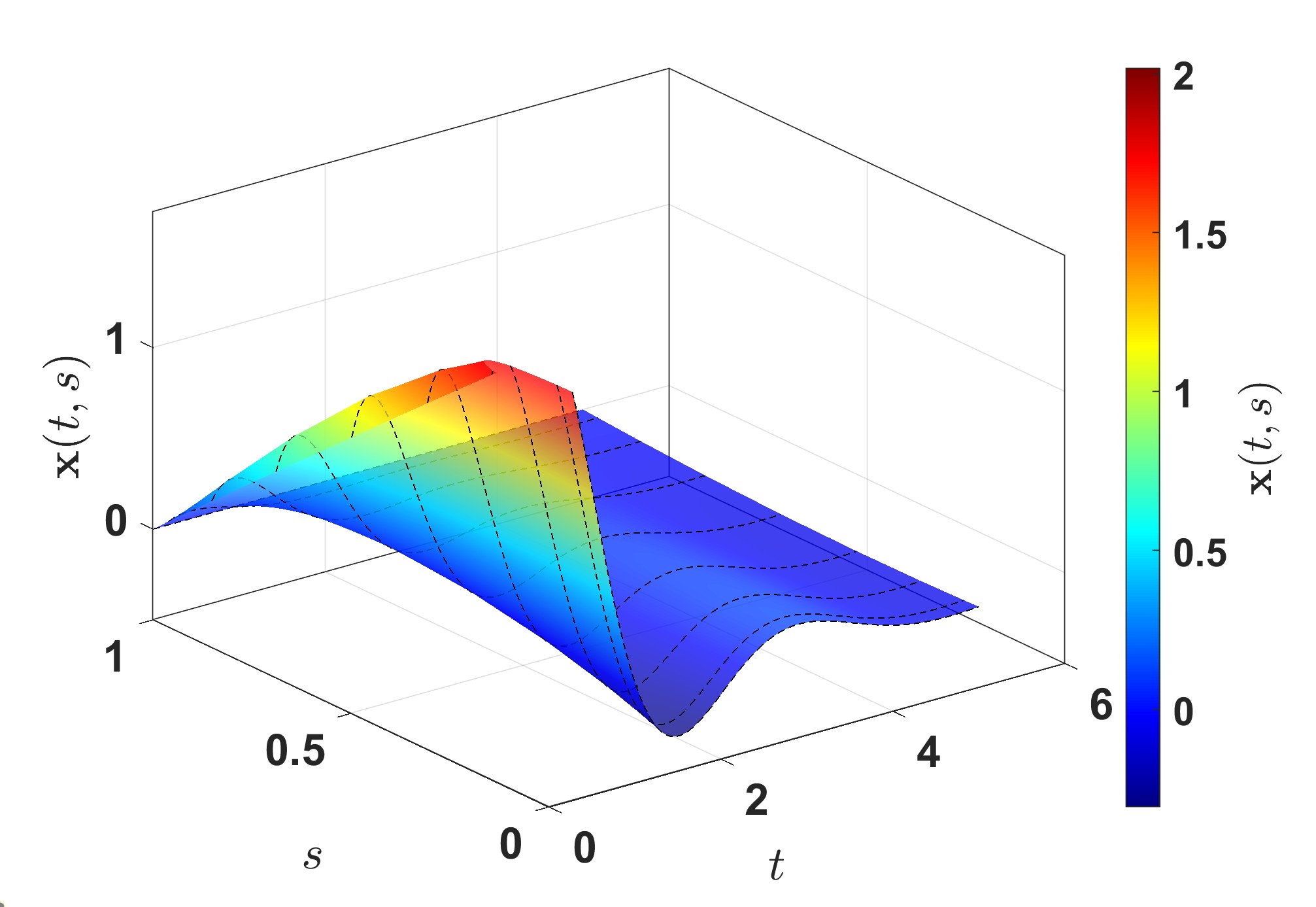}
\caption{A simulation of the open-loop GPDE model defined by Eq.~\eqref{eqn:example3_PDE} in~\Cref{subsec:entropy} via PIESIM module of PIETOOLS and the PIE representation of a GPDE. The initial conditions are set to zero. However, non-zero decaying sinusoidal disturbances, $w_a=sinc(t)$ and $w_b(t) = cos(2t)e^{-t}$, are applied at $t=0$.}
\label{fig:plot2}
\end{figure}

In simulation of such PDEs, using discretization or Galerkin methods, one requires a coupled ODE model such as
\[
\dot{x}_{1:N}(t) = Ax_{1:N}(t)+ B w_a(t),\quad C_b x_{1:N}(t) = C_w w_b(t),
\]
where $x_{1:N}$ is the finite-dimensional discretization (or projection) of the PDE state, $A$ and $B$ are sparse matrices from PDE dynamics, and $C_b$ and $C_w$ are \emph{non-sparse} matrices from the BC conditions. However, typical schemes for numerical simulation of PDEs require a suitable choice of basis functions and often suffers from high complexity and numerical instability~\cite{martin}. The existence of a PIE representation, however, simplifies the process of simulation using the methods proposed in, e.g.~\cite{ypeet}. To illustrate, in~\Cref{fig:plot2}, we use the PDE modeling tool and PIESIM module of PIETOOLS to simulate the proposed PDE model with minimal effort.

\subsection{Control of Reaction-Diffusion PDE} \label{subsec:react_diff} 			
Consider a reaction-diffusion PDE model with an ODE-based controller acting at the boundary.
\begin{align}
&\dot{x}(t) = -x(t)+ u(t), \; \dot{\mbf x}(t,s) = 10\mbf  x(t,s) +\mbf  \partial_s^2 \mbf x(t,s)+w(t),\notag\\
&z(t) = \bmat{\mbf x(t,1)\\ u(t)}, \quad s\in(0,1), \quad  t\ge 0\notag\\
&\mbf x(t,0) =0, \;\mbf  x(t,1) = x(t),\; x(0)=0, \;\mbf x(0,s) = \sin(\pi s),\label{eqn:react_diff_control}
\end{align}
where $x$ is the state of the dynamic boundary controller, $\mbf{x}$ is the distributed state, $z$ is the regulated output and $w$ is a disturbance. The control input, $u(t)$, enters the system through an ODE which is then coupled with the PDE state $\mbf x$ at the boundary. Using the PIETOOLS GUI to define the GPDE, we construct the associated PIE system. The open loop GPDE model is unstable, however, the PIETOOLS tool for stabilizing state-feedback controller synthesis (based on~\cite{shivakumar_2020CDC}) provides the following state-feedback controller.
\begin{align*}
u(t) &= -13.45x(t)+\int_0^1 k(s) \partial_s^2\mbf x(s,t) ds, \qquad \text{where}\\
k(s) &= -9.39s^{10} + 19.7s^9 + 34.7s^8 - 124s^7 + 83.5s^6 \\
&\quad+48.2s^5 - 78.9s^4 + 25.4s^3+ 3.98s^2- 8.73s + 6.61.
\end{align*}

\begin{figure}[t]
\centering
\begin{subfigure}[b]{0.24\textwidth}
\centering
\includegraphics*[width=\textwidth]{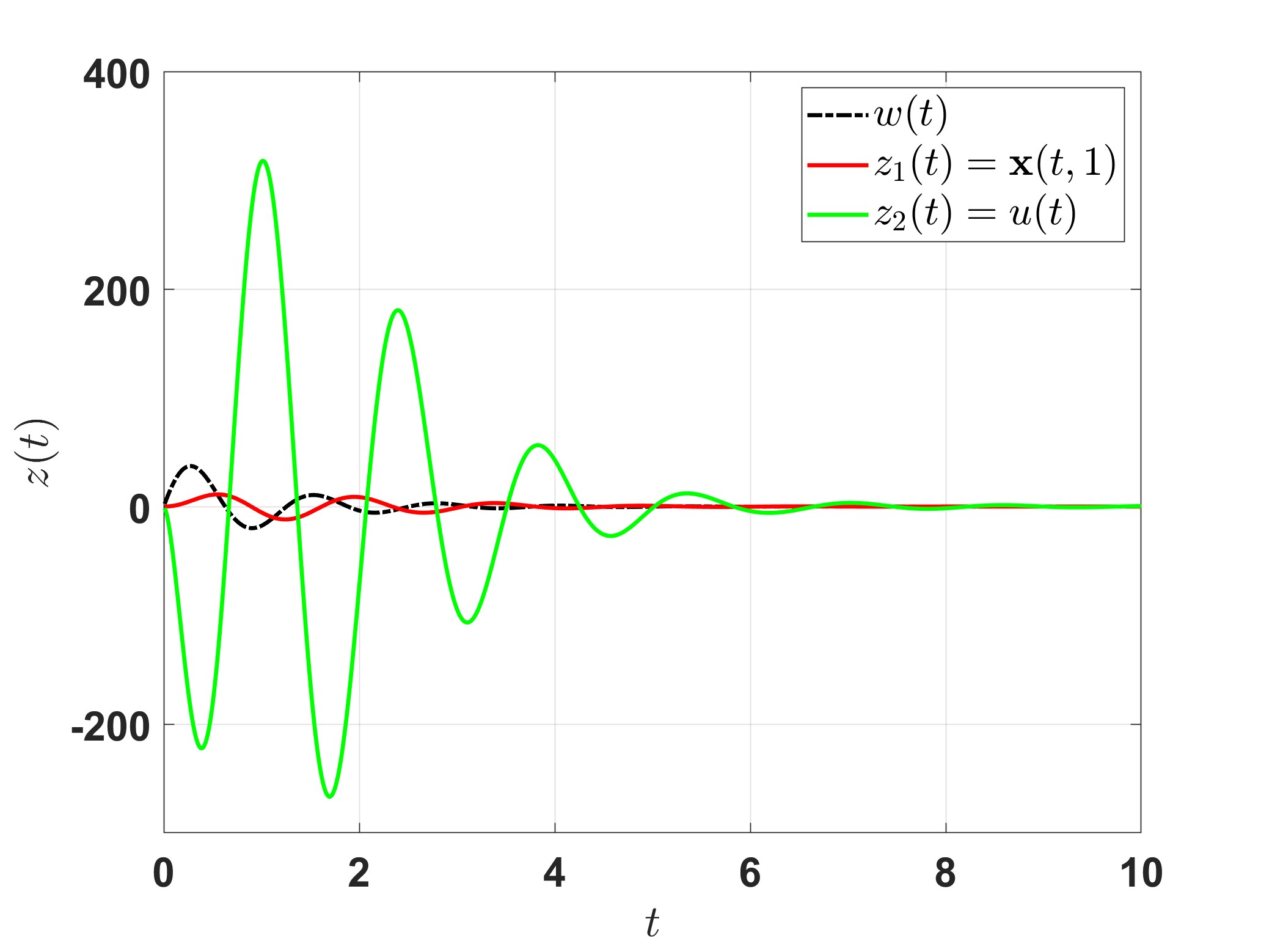}
\caption{Variation of $z(t)$ vs $t$}
\label{fig:output_react_diff}
\end{subfigure}
\begin{subfigure}[b]{0.24\textwidth}
\centering
\includegraphics[width=\textwidth]{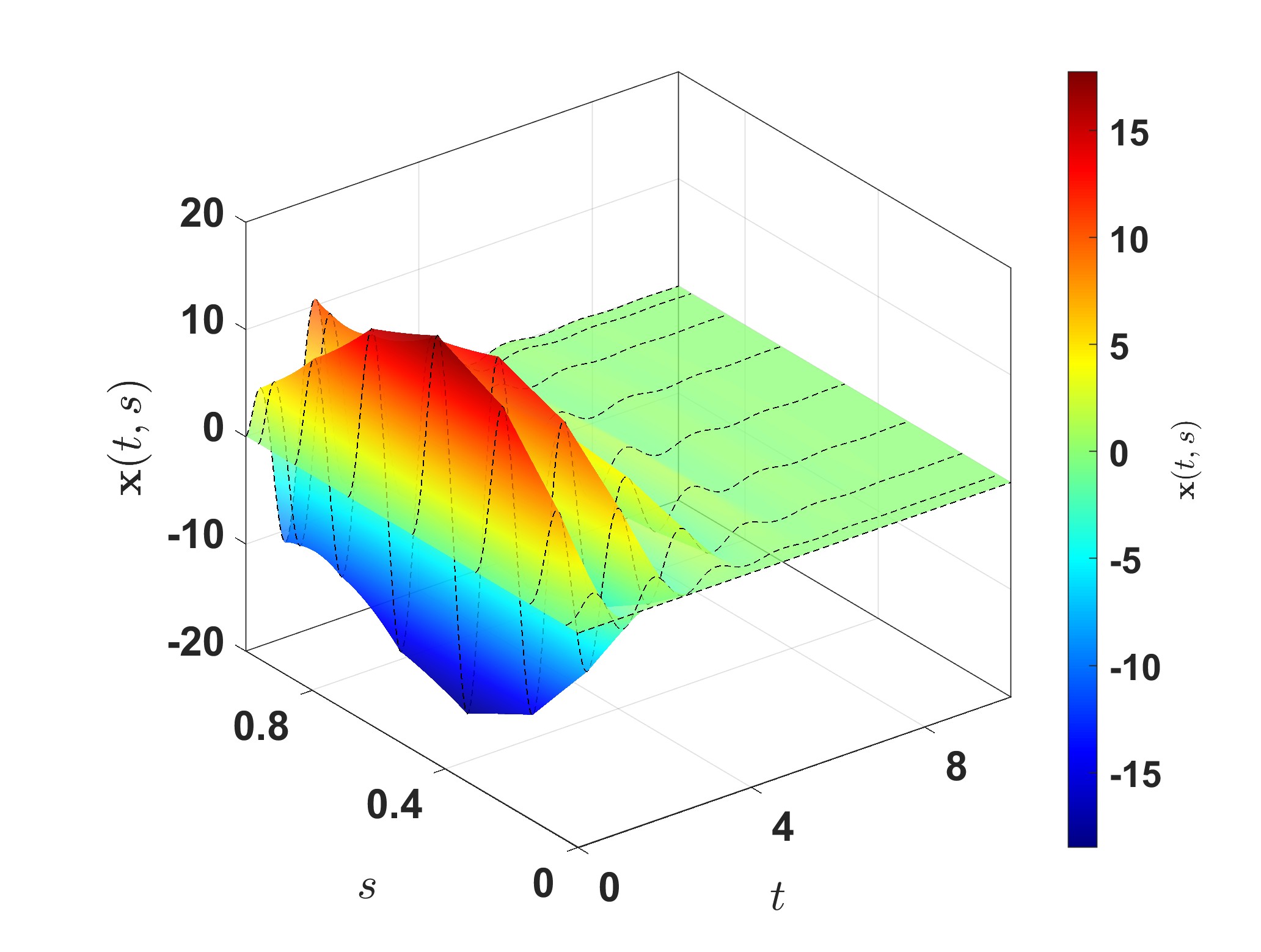}
\caption{Variation of $\mbf x(t)$ vs $t$}
\label{fig:state_react_diff}
\end{subfigure}
\vspace{2mm}
\caption{A simulation of the closed-loop GPDE model defined by Eq.~\eqref{eqn:react_diff_control} in~\Cref{subsec:react_diff} using the stabilizing state-feedback controller generated using PIETOOLS. The initial conditions and decaying sinusoidal disturbance, $w$, are as defined in \cref{subsec:react_diff}.}
\label{fig:plot1}
\end{figure}

We now use the PIESIM package in PIETOOLS to simulate the closed-loop PIE system and reconstruct the GPDE solution where the disturbance is $w(t) = \frac{sin(10t)}{10t+10^{-5}}$. Both the output and control input are shown in~\Cref{fig:plot1}(a) -- verifying that the proposed controller stabilizes the system.

	\section{Conclusion}\label{sec:conclusion}
	A general parametrization was proposed for a class of linear coupled ODE-PDE models (GPDEs) with $N$th-order spatial derivatives that allows for inputs and outputs that enter through: the limit values of the GPDE model, the in-domain dynamics of the PDE subsystem, and a coupled ODE. This parametrization also allows integrals of the PDE state to appear in the ODE-PDE dynamics, the boundary conditions, and the system's outputs --- unifying several existing classes of PDE models in a single parameterized framework.

	Having parameterized a broad class of coupled ODE-PDE models, we proposed a test for the PIE-compatibility of a given well-posed GPDE model. We showed that such compatibility implies the existence of an associated Partial Integral Equation (PIE) representation of the GPDE model with a unitary PI map from the state of the PIE system to the state of the GPDE model.
    Finally, we have shown that many properties of the GPDE model and associated PIE system are equivalent -- including the existence of solutions, input-output properties, internal stability, and controllability.

	To aid in the practical application of the proposed GPDE models and PIE conversion formulae, we have described efficient open-source software (PIETOOLS) for the construction of the GPDE model, conversion to PIE system, simulation of the GPDE/PIE, and analysis/control of the GPDE/PIE -- features demonstrated on several example problems.

	\bibliographystyle{plain}
	\bibliography{references}
	\vspace{-10 mm}	
	\begin{IEEEbiography}[\vspace*{-6mm}{\includegraphics[width=1in,height=1.15in,clip,keepaspectratio]{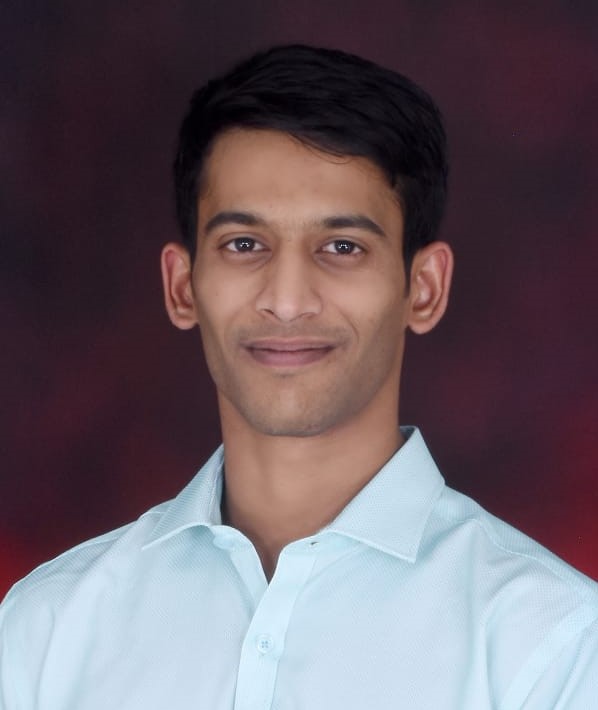}}]{Sachin Shivakumar} received his B.Tech (2015) from the Indian Institute of Technology, Kharagpur, and M.S (2018) Arizona State University (ASU), both in Mechanical Engineering. Since 2018, he has been part of the Cybernetic Systems and Controls Lab (CSCL) at ASU, working on developing mathematical and computational tools to analyze and control PDE systems.
	\end{IEEEbiography}
\vspace{-20 mm}
	\begin{IEEEbiography}[\vspace*{-6mm}{\includegraphics[width=1in,height=1.15in,clip,keepaspectratio]{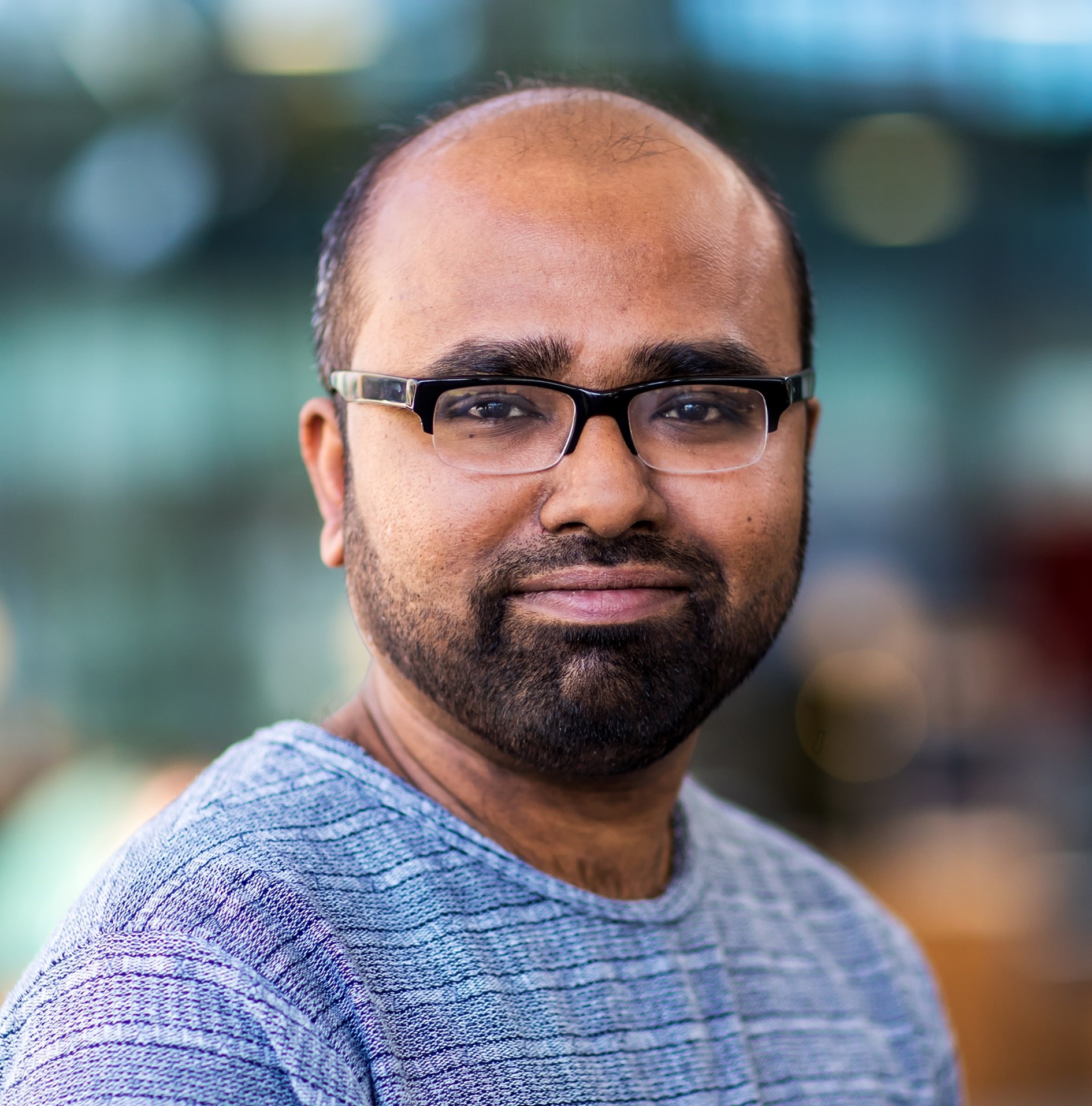}}]{Amritam Das} completed his Bachelors in Mechatronics Engineering in India and joined Eindhoven University of Technology where he finished his masters in 2016 and PhD in 2020. His research interests include thermodynamical systems, nonlinear systems, optimal control, and model reduction.
	\end{IEEEbiography}	
\vspace{-20 mm}
	\begin{IEEEbiography}[\vspace*{-6mm}{\includegraphics[width=1in,height=1.15in,clip,keepaspectratio]{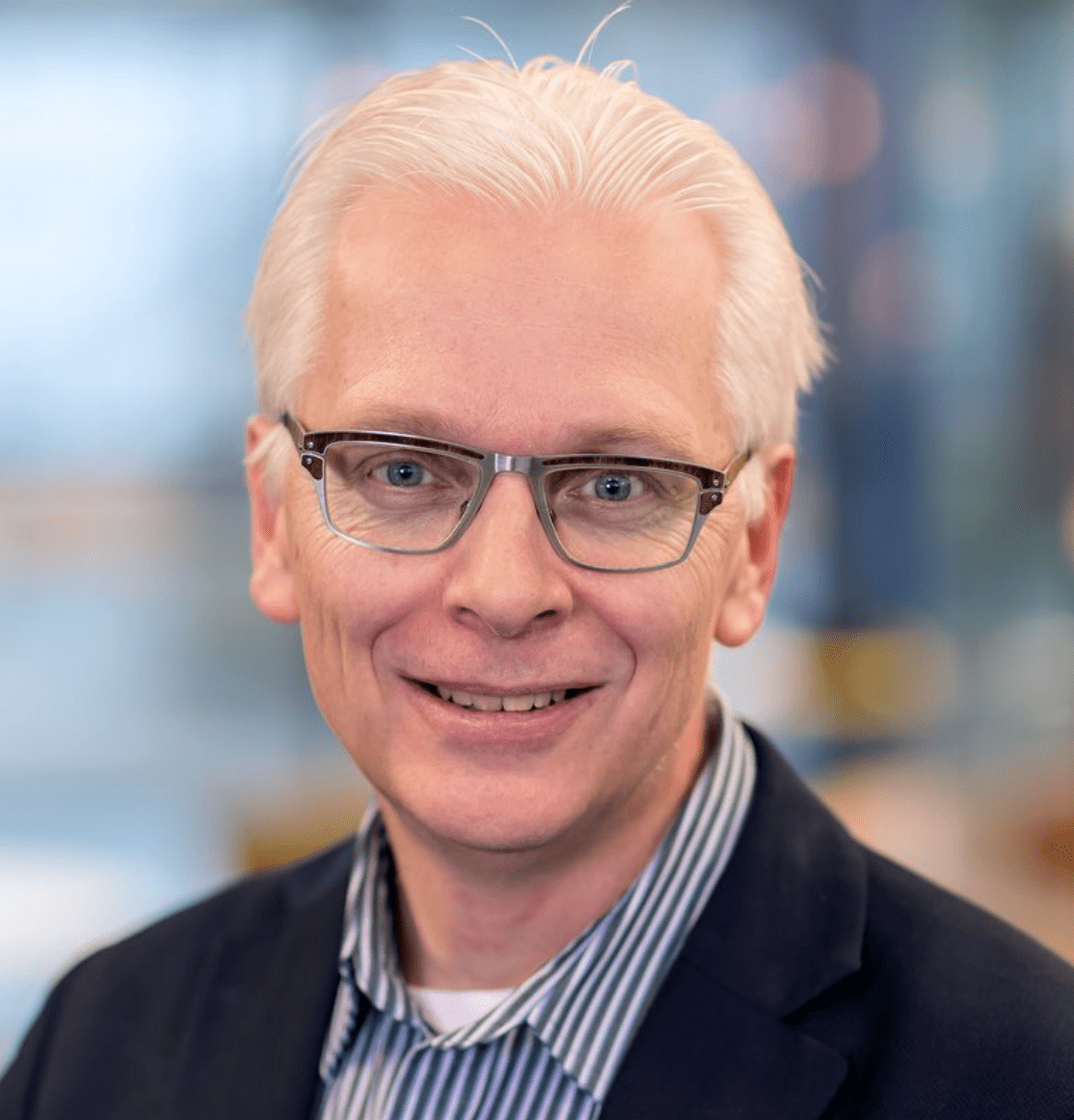}}]{Siep Weiland} received both his MSc (1986) and PhD (1991) degrees in mathematics from the University of Groningen. He worked as a postdoctoral research associate at the Department of Electrical Engineering and Computer Engineering of Rice University, USA (1991-92). Since 1992, he has been affiliated with the Eindhoven University of Technology and was appointed as a full-time professor in Spatial-temporal Systems for Control in 2010.
	\end{IEEEbiography}
\vspace{-16 mm}
	\begin{IEEEbiography}[\vspace*{-6mm}{\includegraphics[width=1in,height=1.15in,clip,keepaspectratio]{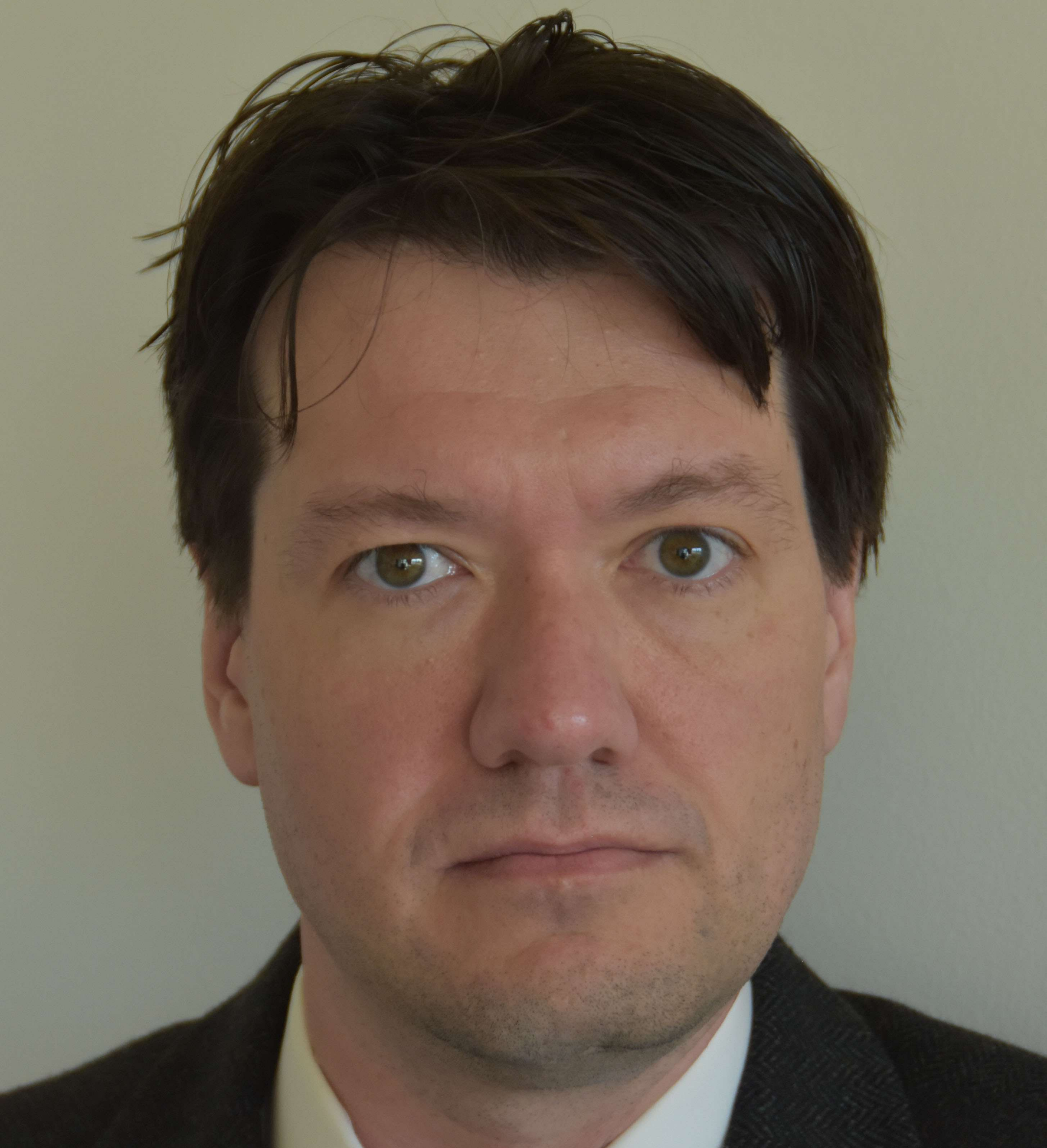}}]{Matthew Peet}
	received M.S. and Ph.D. degrees in aeronautics and astronautics
	from Stanford University (2000-2006). He was a postdoc at INRIA (2006-2008) and Asst. Professor at the Illinois Institute of Technology (2008-2012). Currently, he
	is Associate Professor of Aerospace Engineering
	at Arizona State University.
	\end{IEEEbiography}

	\onecolumn

	\appendix
	\section{Proofs of Theorems}
	In these appendices, we provide proofs for all theorems, lemmas, and corollaries stated in the paper along with several intermediate lemmas. In~\Cref{app:T_map}, the goal is to prove~\Cref{thm:T_map} - the map between the domain of the PIE subsystem and the domain of the PDE subsystem. In~\Cref{app:equivalence}, we prove equivalence of solutions for the PIE subsystem and PDE subsystem. In~\Cref{app:T_map_GPDE} we construct the map between the domain of the GPDE and associated PIE representation. In~\Cref{app:equivalence_2}, we prove equivalence of solutions of the GPDE and associated PIE system. In~\Cref{app:unitary_T,app:norm_equivalence,app:stability_equivalence}, we prove that map from PIE to GPDE state is unitary and that internal stability of PIE and GPDE model is equivalent. Finally, in~\Cref{app:algebra_pi,app:pi_concat} we show that the Partial Integral (PI) operators form a $*$-algebra and provide formulae for composition, adjoint, and concatenation of PI operators.
	\subsection{The map between PIE and PDE states: Proof of \Cref{thm:T_map}}\label{app:T_map}
	To find a map between the \textit{fundamental} state (state of the PIE) and the \textit{primal} state (state of the PDE subsystem), we will use the Fundamental Theorem of Calculus (FTC) and the BCs. First, we recall the FTC and extend it to vector-valued functions on the interval $[a,b]$ as shown below.
	\CauchyLem*
	\begin{proof}
	We prove this using the principle of induction. Suppose the lemma is true for some $N$ and $\mbf x\in W_{N+1}^n[a,b]$. Because the lemma is true for $N$, we have
	\begin{align}\label{eq:induction_step}
	&\mbf x(s) = \mbf x(a) + \sum\limits_{j=1}^{N-1}\frac{(s-a)^j}{j!}(\partial_s^j\mbf x)(a)+\int_{a}^{s}\frac{(s-\theta)^{N-1}}{(N-1)!}(\partial_s^N\mbf x)(\theta)d\theta.
	\end{align}
	Now, by the FTC, we have
	\begin{align*}
	\partial_s^N \mbf x(s) = (\partial_s^N \mbf x)(a)+\int_{a}^{s}(\partial_s^{(N+1)} \mbf x)(\theta)d\theta.
	\end{align*}
	Next, we substitute the above identity into \Cref{eq:induction_step}, and using the integral identity
	\[
	\int_{a}^b \int_{a}^\theta f(\theta,\eta)d \eta d\theta
	=\int_{a}^b \int_{\eta}^b f(\theta,\eta)d \theta d\eta
\]
we have
\begin{align*}
&\mbf x(s) =\mbf x(a) + \sum\limits_{j=1}^{N-1}\frac{(s-a)^j}{j!}\partial_s^j\mbf x(a)
+\int_{a}^{s}\frac{(s-\theta)^{N-1}}{(N-1)!}\left(\partial_s^N \mbf x(a)+\int_{a}^{\theta}\partial_s^{(N+1)} \mbf x(\eta)d\eta\right)d\theta.
\end{align*}
While the first two terms are close to the required form, the last term (the integral term) is not and can be simplified by using integration by parts. Then, from the integral term we get
\begin{align*}
&\int_{a}^{s}\frac{(s-\theta)^{N-1}}{(N-1)!}\left(\partial_s^N \mbf x(a)+\int_{a}^{\theta}\partial_s^{(N+1)} \mbf x(\eta)d\eta\right)d\theta\\
&=\left(\int_{a}^{s}\frac{(s-\theta)^{N-1}}{(N-1)!}d\theta\right) \partial_s^N \mbf x(a)
+\int_{a}^{s}\int_{a}^{\theta}\frac{(s-\theta)^{N-1}}{(N-1)!}\partial_s^{(N+1)} \mbf x(\eta)d\eta \, d\theta\\
&= \frac{(s-a)^N}{N!}\partial_s^N\mbf x(a)
+\int_{a}^{s}\left(\int_{\eta}^{s}\frac{(s-\theta)^{N-1}}{(N-1)!}d\theta\right)\partial_s^{(N+1)} \mbf x(\eta)d\eta\\
&= \frac{(s-a)^N}{N!}\partial_s^N\mbf x(a)+\int_{a}^{s}\frac{(s-\eta)^{N}}{N!}\partial_s^{(N+1)} \mbf x(\eta)d\eta.
\end{align*}
Finally, by substituting the above terms back into the equation we get,
\begin{align*}
\mbf x(s) =\mbf x(a) + \sum\limits_{j=1}^{N}\frac{(s-a)^j}{j!}\partial_s^j\mbf x(a)+\int_{a}^{s}\frac{(s-\eta)^{N}}{N!}\partial_s^{(N+1)} \mbf x(\eta)d\eta.
\end{align*}
Therefore, if the statement of the lemma is true for $N$, then it is also true for $N+1$. Clearly, the lemma is true for $N=1$.
\end{proof}

We can extend Lemma~\ref{lem:cauchy} to obtain an expression for the derivatives of $\mbf{x}\in C_N^{n}$ in terms of $\partial_s^N\mbf x$ and of a given set of \textit{core} boundary values of $\mbf x$.

\begin{lemma}\label{lem:ftc_derivatives1}
Suppose $\mbf{x}\in W_{N}^n$. Then, for any $i < N$, we have
\begin{align*}
(\partial_s^{i}\mbf x)(s) = \sum\limits_{j=i}^{N-1}\tau_{j-i}(s-a)(\partial_s^{j}\mbf x)(a)+\int_{a}^{s}\tau_{N-i-1}(s-\theta)(\partial_s^N\mbf x)(\theta)d\theta
\end{align*}
where $\tau_i(s)=\frac{s^i}{i!}$.
\end{lemma}
\begin{proof}
First note that $\tau_i(0)=0$ for any $i > 0$, $\partial_s \tau_0(s)=0$ and
\[
\tau_i(s) = \frac{s^{i}}{i!}\qquad  \rightarrow \qquad  \partial_s\tau_i(s) = i\frac{s^{i-1}}{i!}=\tau_{i-1}(s)
\]
and suppose the formula holds for $i-1\ge 0$. Then
\[
(\partial_s^{i-1}\mbf x)(s) = \sum\limits_{j=i-1}^{N-1}\tau_{j-i+1}(s-a)(\partial_s^{j}\mbf x)(a)+\int_{a}^{s}\tau_{N-i}(s-\theta)(\partial_s^N\mbf x)(\theta)d\theta.
\]
and hence, since $\partial_s \tau_0(s)=0$, we have
{\small
\begin{align*}
&	(\partial_s^{i}\mbf x)(s)=		\partial_s(\partial_s^{i-1}\mbf x)(s)  \\
&= \sum\limits_{j=i-1}^{N-1}\left(\partial_s\tau_{j-i+1}(s-a)\right)(\partial_s^{j}\mbf x)(a)+\tau_{N-i}(0)(\partial_s^N\mbf x)(s)+\int_{a}^{s}\left(\partial_s\tau_{N-i}(s-\theta)\right)(\partial_s^N\mbf x)(\theta)d\theta.\\
&=  \sum\limits_{j=i}^{N-1}\left(\partial_s\tau_{j-i+1}(s-a)\right)(\partial_s^{j}\mbf x)(a)+\tau_{N-i}(0)(\partial_s^N\mbf x)(s)+\int_{a}^{s}\left(\partial_s\tau_{N-i}(s-\theta)\right)(\partial_s^N\mbf x)(\theta)d\theta.\\
&=  \sum\limits_{j=i}^{N-1}\tau_{j-i}(s-a)(\partial_s^{j}\mbf x)(a)+\int_{a}^{s}\tau_{N-i-1}(s-\theta)(\partial_s^N\mbf x)(\theta)d\theta.
%
\end{align*}
}
By \cref{lem:cauchy}, the result holds for $i=0$, which completes the proof.
\end{proof}
We now propose a mixed-order version of Lemma~\ref{lem:cauchy}

\begin{corollary}\label{cor:ftc_primary}
Suppose $\mbf{x}\in \prod_{i=0}^NW_i^{n_i}$ and define
\begin{align*}
&J_{i,j} =\bmat{0_{n_{i:j-1}\times n_{j:N}}\\I_{n_{j:N}}}\in \R^{n_{i:N} \times n_{j:N}},\qquad \tau_i(s) = \frac{s^{i}}{i!},\qquad T_{i,j}(s) =\tau_{j-i}(s)J_{i,j} \quad j\ge i,\qquad\mcl C\mbf{x}:= \bmat{S\mbf x\\ \partial_s S^2\mbf x \\ \vdots \\ \partial_s^{N-1}S^{N}\mbf x},
\end{align*}
we have
\begin{align}
\bmat{\mbf{x}_1(s)\\\vdots\\\mbf x_N(s)} &= \bmat{T_{1,1}(s-a)&T_{1,2}(s-a)&\cdots&T_{1,N}(s-a)}(\mcl C\mbf{x})(a) \notag\\
&\quad+ \int_{a}^{s}\bmat{\tau_0(s-\theta)I_{n_1}&&\\&\ddots&\\&&\tau_{N-1}(s-\theta)I_{n_k}}\bmat{\partial_{\theta}\mbf{x}_1(\theta)\\\vdots\\\partial_{\theta}^N \mbf x_N(\theta)}d\theta\notag\\
&=T_1(s-a)(\mcl C\mbf{x})(a)+\int_{a}^{s}Q_1(s-\theta)\bmat{\mbf{x}_0(\theta)\\\vdots\\\partial_{\theta}^N \mbf x_N(\theta)}d\theta.
\end{align}
\end{corollary}

\begin{proof}
For convenience, let us denote $P_{i,j}\in \R^{n_i \times n_{S_{j}}}$ to be the uniquely defined 0-1 matrix so that $\mbf x_i(s)=P_{i,j}S^j \mbf x(s)$ and which is given by
\[
P_{i,j}:=\bmat{0_{n_i\times(n_{S_j}-n_{S_{i-1}})}&I_{n_{i}}&0_{n_i\times n_{S_{i+1}}}}=\bmat{0_{n_i\times n_{j:i-1}}&I_{n_{i}}&0_{n_i\times n_{i+1:N}}}.
\]

We now use $P_{i,j}$  and the identity from \Cref{lem:cauchy} to write $\mbf x_i$ in terms of $(\mcl C\mbf x)(a)$ and $\partial_s^i\mbf x_i$. Specifically, if $\mbf x_k\in C_k^{n_k}[a,b]$, then		{\small
\begin{align*}
\mbf x_i(s) &= \sum\limits_{j=0}^{i-1}\tau_j(s-a)\partial_s^j\mbf x_i(a)+\int_{a}^{s}\tau_{i-1}(s-\theta)\partial_s^i\mbf x_i(\theta)d\theta\\
&= \sum\limits_{j=0}^{i-1}\tau_j(s-a) P_{i,j+1}\partial_s^j S^{j+1} \mbf x(a)+\int_{a}^{s}\tau_{i-1}(s-\theta)\partial_s^i\mbf x_i(\theta)d\theta\\
&=\bmat{ \tau_0(s-a) P_{i,1} &\cdots & \tau_{i-1}(s-a) P_{i,i}}\bmat{S\mbf x(a)\\ \vdots \\ \partial_s^{i-1} S^{i} \mbf x(a)}+\int_{a}^{s}\tau_{i-1}(s-\theta)\partial_s^i\mbf x_i(\theta)d\theta\\
&=\bmat{ \tau_0(s-a) P_{i,1} &\cdots & \tau_{i-1}(s-a) P_{i,i} & 0_{n_i\times n_{S_{i+1:N}}}}\bmat{S\mbf x(a)\\ \vdots \\ \partial_s^{N-1} S^{N} \mbf x(a)}+\int_{a}^{s}\tau_{i-1}(s-\theta)\partial_s^i\mbf x_i(\theta)d\theta\\
&=\bmat{ \tau_0(s-a) P_{i,1} &\cdots & \tau_{i-1}(s-a) P_{i,i} & 0_{n_i\times n_{S_{i+1:N}}}}(\mcl C\mbf x)(a)+\int_{a}^{s}\tau_{i-1}(s-\theta)\partial_s^i\mbf x_i(\theta)d\theta.
\end{align*}
}	
Now, we can concatenate the $\mbf x_i$'s to get,
\begin{align*}
\bmat{\mbf{x}_1(s)\\\vdots\\\mbf x_N(s)} &=\bmat{\tau_0(s-a)P_{1,1} &0&0\\\vdots&\ddots&0\\ \tau_0(s-a)P_{N,1}&\cdots&\tau_{N-1}(s-a)P_{N,N} }(\mcl C\mbf
x)(a)+\int_{a}^{s}\bmat{\tau_0(s-\theta)\partial_s\mbf x_1(\theta)\\\vdots\\\tau_{N-1}(s-\theta)\partial_s^N\mbf x_N(\theta)}d\theta\\
&=\bmat{\tau_0(s-a)\bmat{P_{1,1}\\P_{2,1}\\\vdots\\P_{N,1}} &\tau_1(s-a)\bmat{0\\P_{2,2}\\\vdots\\P_{N,2}}&\cdots&\tau_{N-1}(s-a)\bmat{0\\0\\\vdots\\P_{N,N}}}(\mcl C\mbf x)(a)\\
&\qquad +\int_{a}^{s}\bmat{\tau_0(s-\theta)\partial_s\mbf x_1(\theta)\\\vdots\\\tau_{N-1}(s-\theta)\partial_s^N\mbf x_N(\theta)}d\theta\\
&=\bmat{\tau_0(s-a)J_{1,1}&\tau_1(s-a)J_{1,2}&\cdots&\tau_{N-1}(s-a)J_{1,N}}(\mcl C\mbf x)(a)\\
&\qquad+ \int_{a}^{s}\bmat{\tau_0(s-\theta)&&\\&\ddots&\\&&\tau_{N-1}(s-\theta)}\bmat{\partial_{\theta}\mbf{x}_1(\theta)\\\vdots\\\partial_{\theta}^N \mbf x_N(\theta)}d\theta\\
&= \bmat{T_{1,1}(s-a)&\cdots&T_{1,N}(s-a)}(\mcl C\mbf{x})(a) \\
&\qquad+ \int_{a}^{s}\bmat{\tau_0(s-\theta)&&\\&\ddots&\\&&\tau_{N-1}(s-\theta)}\bmat{\partial_{\theta}\mbf{x}_1(\theta)\\\vdots\\\partial_{\theta}^N \mbf x_N(\theta)}d\theta
\end{align*}
where we have used the fact that for any $i$
\begin{align*}
\bmat{P_{i,i}\\P_{i+1,i}\\\vdots\\P_{N,i}} &= \bmat{\bmat{0_{n_i\times n_{i:i-1}}&I_{n_{i}}&0_{n_i\times n_{i+1:N}}}\\\bmat{0_{n_{i+1}\times n_{i:i}}&I_{n_{i+1}}&0_{n_{i+1}\times n_{i+2:N}}}\\\vdots\\\bmat{0_{n_N\times n_{i:N-1}}&I_{n_{N}}&0_{n_N\times n_{N+1:N}}}} = \bmat{I_{n_{i}}&&&\\&I_{n_{i+1}}&&\\&&\ddots&\\&&&I_{n_N}}= I_{n_{i:N}}
\end{align*}
and hence
\begin{align*}
\bmat{0_{n_{1:i-1}\times n_{i:N}}\\P_{i,i}\\\vdots\\P_{N,i}} &= \bmat{0_{n_{1:i-1}\times n_{i:N}}\\I_{n_{i:N}}} = J_{1,i}.
\end{align*}
\end{proof}

We conclude that it is possible to express any function $\mbf x\in \prod_{i=1}^N W_i^{n_i}$ using left boundary values (at $s=a$) of the continuous partial derivatives $(\mcl C\mbf x)$ and the fundamental state $\ubar{\hat{\mbf x}} = \text{col}(\mbf x_0,\cdots, \partial_s^N\mbf x_N)$. Since we require a map from the fundamental state $\ubar{\hat{\mbf x}}$ to the primal state $\mbf x$, we need to eliminate the left boundary values $(\mcl C\mbf x)(a)$. The first step in this direction is to express $\mcl C\mbf x$ in terms of $\ubar{\hat{\mbf x}}$ and $(\mcl C\mbf x)(a)$.

\begin{corollary}\label{cor:x_c}
Suppose $\mbf{x}\in \prod_{i=0}^NW_i^{n_i}$. Then, for $T$ and $Q$ are as defined in \Cref{fig:PIE_subsystem_equation}, and
\begin{align*}
(\mcl C\mbf x):= \bmat{S\mbf x\\ \partial_s S^2\mbf x \\ \vdots \\ \partial_s^{N-1}S^{N}\mbf x}, \qquad \ubar{\hat{\mbf x}} := \bmat{\mbf{x}_0\\\partial_s^1\mbf{x}_{1}\\\vdots\\\partial_s^N\mbf{x}_N},
\end{align*}
we have
\[
(\mcl C\mbf x)(s)= T(s-a)(\mcl C\mbf x)(a)+\int_{a}^{s}Q(s-\theta)\ubar{\hat{\mbf x}}(\theta)d\theta.
\]	\end{corollary}

\begin{proof}
We will use the identity from \cref{cor:ftc_primary,lem:ftc_derivatives1} to find $\partial_s^{i-1}S^i\mbf x$ for all $1\le i\le N$ and concatenate them vertically to obtain $(\mcl C\mbf x)$. First, we need to find an expression for $\partial_s^{i-1}S^i\mbf x$. By definition, we have
\[
\partial_s^{i-1}S^i \mbf x(s)=  \bmat{\partial_s^{i-1}\mbf x_i(s)\\\partial_s^{i-1} \mbf x_{i+1}(s)\\\vdots \\\partial_s^{i-1}\mbf x_N(s)}.
\]
By \cref{lem:ftc_derivatives1},
\begin{align*}
(\partial_s^{i}\mbf x)(s) = \sum\limits_{j=i}^{N-1}\tau_{j-i}(s-a)(\partial_s^{j}\mbf x)(a)+\int_{a}^{s}\tau_{N-i-1}(s-\theta)(\partial_s^N\mbf x)(\theta)d\theta
\end{align*}i
which can be generalized for $\mbf x_k \in \mcl C_{k}$ with $k< N$ as
\begin{align*}
(\partial_s^{i}\mbf x_k)(s) = \sum\limits_{j=i}^{k-1}\tau_{j-i}(s-a)(\partial_s^{j}\mbf x_k)(a)+\int_{a}^{s}\tau_{k-i-1}(s-\theta)(\partial_s^k\mbf x_k)(\theta)d\theta.
\end{align*}
To find the $(i-1)^{th}-$ derivative for each component of the vector we just perform concatenation to get
{\small	\begin{align*}
&\partial_s^{i-1}S^i \mbf x(s)=\bmat{\partial_s^{i-1}\mbf x_i(s)\\\partial_s^{i-1} \mbf x_{i+1}(s)\\\vdots \\\partial_s^{i-1}\mbf x_N(s)}\\
&\quad = \bmat{(\partial_s^{i-1}\mbf x_i)(a)\\(\partial_s^{i-1}\mbf x_{i+1})(a)+\tau_1(s-a)(\partial_s^i \mbf x_{i+1})(a)\\\vdots\\(\partial_s^{i-1}\mbf x_N)(a)+\sum\limits_{j=i}^{N-1}\tau_{j-i+1}(s-a)(\partial_s^j\mbf x_N)(a)}+\int_a^s \bmat{(\partial_s^{i} \mbf x_i)(\theta)\\\tau_1(s-\theta)(\partial_s^{i+1}\mbf x_{i+1})(\theta)\\\vdots\\\tau_{(N-i)}(s-\theta)(\partial_s^{N} \mbf x_N)(\theta)}d\theta.
\end{align*}
}
The matrices $J_{i,j}$ for $j>i$ are used to select the elements from $(\partial_s^{j-1}S^{j}\mbf x)(a)\in \R^{n_{j:N}}$ ($j$th part of $(\mcl C\mbf x)(a)$) which appear in the $(j-i)$th to $(N-i)$th components of $(\partial_s^{i-1}S^{i}\mbf x)(s)\in \R^{n_{i:N}}$ ($i$th part of $(\mcl C\mbf x)(s)$). Specifically, for $j\ge i$, we will see that $\partial_s^{i-1}S^{i}\mbf x$ is the combination of terms of the form
\begin{align*}
&\bmat{0_{n_{i:j-1}\times 1}\\(\partial_s^{j-1} \mbf x_{j})(s)\\\vdots\\(\partial_s^{j-1}\mbf x_N)(s)}=\bmat{0_{n_{i:j-1}\times n_{j:N}}\\I_{n_{j:N}}}\bmat{(\partial_s^{j-1} \mbf x_{j})(s)\\\vdots\\(\partial_s^{j-1}\mbf x_N)(s)}=\bmat{0_{n_{i:j-1}\times n_{j:N}}\\I_{n_{j:N}}}(\partial_s^{j-1}S^{j}\mbf x)(s)=J_{i,j}(\partial_s^{j-1}S^{j}\mbf x)(s).
\end{align*}
By exploiting the $J_{i,j}$ notation, we are able to represent the first term in the expression for $\partial_s^{i-1}S^{i}\mbf x$ as
\begin{align*}
\bmat{(\partial_s^{i-1}\mbf x_i)(a)\\(\partial_s^{i-1}\mbf x_{i+1})(a)+\tau_1(s-a)(\partial_s^i \mbf x_{i+1})(a)\\\vdots\\(\partial_s^{i-1}\mbf x_N)(a)+\sum\limits_{j=i}^{N-1}\tau_{j-i+1}(s-a)(\partial_s^j\mbf x_N)(a)}
&= \tau_0(s-a)\bmat{\partial_s^{i-1}\mbf x_i(a)\\\partial_s^{i-1}\mbf x_{i+1}(a)\\\vdots\\\partial_s^{i-1}\mbf x_N(a)}
+\cdots
+\tau_{j-i}(s-a)\bmat{0_{n_{i:j-1}\times 1}\\(\partial_s^{j-1} \mbf x_{j})(a)\\\vdots\\(\partial_s^{j-1}\mbf x_N)(a)}+\cdots\\
&\quad+\tau_{N-i}(s-a)\bmat{0_{n_i}\\0_{n_{i+1}}\\\vdots\\\partial_s^{N-1}\mbf x_N(a)}\\
&= \tau_0(s-a)(\partial_s^{i-1}S^i\mbf x)(a)+\cdots+\tau_{j-i}(s-a)J_{i,j}(\partial_s^{j-1}S^{j}\mbf x)(a)+\cdots\\
&\quad +\tau_{N-i}(s-a)J_{i,N}(\partial_s^{N-1} S^{N}\mbf x)(a)\\
&= \bmat{\tau_0(s-a)J_{i,i}&\cdots &\tau_{N-i}(s-a)J_{i,N}}\bmat{\partial_s^{i-1} S^i \mbf x(a)\\\vdots\\\partial_s^{N-1}S^N\mbf x(a)}\\
&= \bmat{0_{n_{i:N}\times n_{S_{1:i-1}}}&\tau_0(s-a)J_{i,i}&\cdots &\tau_{N-i}(s-a)J_{i,N}}\bmat{S\mbf x(a)\\ \vdots \\\partial_s^{i-1} S^i \mbf x(a)\\\partial_s^{i-1}S^{i+1}\mbf x(a)\\\vdots\\\partial_s^{N-1}S^N\mbf x(a)}\\
&= T_{i}(s-a)(\mcl C\mbf x)(a).
\end{align*}
The second vector with the integral terms can be written as
\begin{align*}
\int_a^s \bmat{\partial_s^{i} \mbf x_i(\theta)\\\tau_1(s-\theta)\partial_s^{i+1}\mbf x_{i+1}(\theta)\\\vdots\\\tau_{N-i}(s-\theta)\partial_s^{N} \mbf x_N(\theta)}d\theta& =\int_a^s \bmat{0&I_{n_i}&&&\\0&&\tau_1(s-\theta)I_{n_{i+1}}&&\\\vdots&&&\ddots&\\0&&&&\tau_{N-i}(s-\theta)I_{n_N}} \bmat{\mbf x_0\\\partial_s \mbf x_1(\theta)\\\vdots\\\partial_s^N \mbf x_N(\theta)}d\theta\\
&= \int_a^s Q_i(s-\theta)\ubar{\hat{\mbf x}}(\theta)d\theta
\end{align*}
Therefore, combining both the terms we get
\begin{align*}
\partial_{s}^{i-1}S^i\mbf{x}(s) &=  T_{i}(s-a)(\mcl C\mbf x)(a)+\int_{a}^{s}Q_i(s-\theta)\ubar{\hat{\mbf x}}(\theta)d\theta.
\end{align*}
Since $\partial_{s}^{i-1}S^i\mbf{x}(s)$ can be uniquely determined using $(\mcl C\mbf x)(a)$ and $\ubar{\hat{\mbf x}}$, we can now generalize this to all of $(\mcl C\mbf x)$ by using concatenation of $\partial_{s}^{i-1}S^i\mbf{x}(s)$ over all $i\in n$ as
\begin{align*}
(\mcl C\mbf x)(s) &= \bmat{S\mbf x(s)\\ \vdots \\ \partial_s^{N-1}S^{N}\mbf x(s)}= \bmat{T_{1}(s-a)\\T_{2}(s-a)\\\vdots\\T_{N}(s-a)}\bmat{S\mbf x(a)\\ \partial_s S^2\mbf x(a) \\ \vdots \\ \partial_s^{N-1}S^{N}\mbf x(a)}+\int_{a}^{s}\bmat{Q_1(s-\theta)\\\vdots\\Q_N(s-\theta)}\ubar{\hat{\mbf x}}(\theta)d\theta\\
&= T(s-a)(\mcl C\mbf x)(a)+\int_{a}^{s}Q(s-\theta)\ubar{\hat{\mbf x}}(\theta)d\theta.
\end{align*}
\end{proof}

Next, we use the map from $(\mcl C\mbf x)(a)$ and $\ubar{\hat{\mbf x}}$ to $(\mcl C\mbf x)$ to obtain the following list of identities.

\begin{corollary}\label{cor:x_D}
Suppose $\hat{\mbf{x}}\in W^n$ for some $v\in \R^q$. Define
\begin{align*}
(\mcl F\hat{\mbf x}):= \bmat{\hat{\mbf x}\\ \partial_s S\hat{\mbf x} \\ \vdots \\ \partial_s^{N}S^{N}\hat{\mbf x}}, \qquad (\mcl C\hat{\mbf x}):= \bmat{S\hat{\mbf x}\\ \partial_s S^2\hat{\mbf x} \\ \vdots \\ \partial_s^{N-1}S^{N}\hat{\mbf x}}, \qquad \ubar{\hat{\mbf x}} := \bmat{\hat{\mbf x}_0\\\partial_s^1\hat{\mbf x}_{1}\\\vdots\\\partial_s^N\hat{\mbf x}_N}.
\end{align*}
Then we have the following.
\begin{enumerate}[label=(\alph*)]
\item For $U_i$, $T$, and $Q$ as defined in \Cref{fig:PIE_subsystem_equation},
we have
\begin{align*}
(\mcl F\hat{\mbf x})(s)&=U_1\ubar{\hat{\mbf x}}(s) + U_2(\mcl C\hat{\mbf x})(s) \\
&= U_2 T(s-a)(\mcl C\hat{\mbf x})(a)+  U_1\ubar{\hat{\mbf x}}(s)+\int_{a}^{s}U_2 Q(s-\theta)\ubar{\hat{\mbf x}}(\theta)d\theta.
\end{align*}
\item Given a set of parameters $\mbf G_{\mathrm b}$, if $\{n,\mbf G_{\mathrm b}\}$ is PIE-compatible, then for $B_Q$ as defined in \cref{fig:PIE_subsystem_equation}
\begin{align*}
B\bmat{(\mcl C\hat{\mbf x})(a)\\ (\mcl C\hat{\mbf x})(b)}-\int_{a}^{b}B_{I}(s)(\mcl F\hat{\mbf x})(s)ds&= B_T(\mcl C\hat{\mbf x})(a)-\int_{a}^{b}B_TB_Q(s)\ubar{\hat{\mbf x}}(s)ds.
\end{align*}
\item Given a set of parameters $\mbf G_{\mathrm b}$, if $\hat{\mbf x} \in X_v$ and $\{n,\mbf G_{\mathrm b}\}$ is PIE-compatible, then
\[
(\mcl C\hat{\mbf x})(a)= \int_{a}^{b}B_Q(\theta)\ubar{\hat{\mbf x}}(\theta)d\theta+B_T^{-1}B_v v.
\]
\item Given a set of parameters $\mbf G_{\mathrm b}$, if $\hat{\mbf x} \in X_v$ and $\{n,\mbf G_{\mathrm b}\}$ is PIE-compatible, then
\[
\bmat{\bmat{v\\(\mcl B \hat{\mbf x})}\\(\mcl F\hat{\mbf x})(\cdot)}= \fourpi{\bmat{I_{n_v}\\B_T^{-1}B_v\\T(b-a)B_T^{-1}B_v}}{\bmat{0_{n_r \times n_x}\\B_Q(s)\\T(b-a)B_Q(s)+Q(b-s)}}{U_2T(s-a)B_T^{-1}B_v}{U_1,R_{D,1},R_{D,2}}\bmat{v\\\ubar{\hat{\mbf x}}}.
\]
\end{enumerate}
\end{corollary}
\begin{proof}
Let $\hat{\mbf x}\in X_v$ for some $v\in \R^{q}$.

For (a), we examine the terms $\partial_s^i S^i\hat{\mbf x}$ in the vector $(\mcl F\hat{\mbf x})$. These terms may be divided into those from $\ubar{\hat{\mbf x}}$ and those from $(\mcl C\hat{\mbf x})$. Specifically, we define the permutation matrix $U=\bmat{U_1 & U_2}$  so that
\[
\underbrace{\bmat{\hat{\mbf x}(s)\\ \partial_s S \hat{\mbf x}(s)\\ \vdots \\ (\partial_s^N S^N\hat{\mbf x})(s)}}_{(\mcl F\hat{\mbf x})}=U \bmat{\hat{\mbf x}(s)\\  \bmat{S\hat{\mbf x}(s)\\ \vdots \\ \partial_s^{N-1}S^{N}\hat{\mbf x}(s)}}=\bmat{U_1 & U_2} \bmat{\ubar{\hat{\mbf x}}(s)\\ (\mcl C\hat{\mbf x})(s)}.
\]
To justify our expression for the permutation matrix, $U$, first note that
\begin{align*}
\partial_s^iS^i\hat{\mbf x} &= \bmat{\partial_s^i \hat{\mbf x}_i\\\partial_s^i \hat{\mbf x}_{i+1}\\ \vdots \\\partial s^i \hat{\mbf x}_N}=\bmat{\ubar{\hat{\mbf x}}_{i}\\ \partial_s^{i} S^{i+1} \hat{\mbf x}}=\bmat{\ubar{\hat{\mbf x}}_{i}\\ (\mcl C\hat{\mbf x})_{i+1}}\\
&=\bmat{\bmat{I_{n_i} \\ 0_{ n_{i+1:N}\times n_i}} & \bmat{0_{n_i\times n_{i+1:N} }\\  I_{n_{i+1:N}}}}\bmat{\ubar{\hat{\mbf x}}_{i}\\ (\mcl C\hat{\mbf x})_{i+1}}=\bmat{U_{1,i}&U_{2,i}}\bmat{\ubar{\hat{\mbf x}}_{i}\\ (\mcl C\hat{\mbf x})_{i+1}}
\end{align*}
which holds for $i<N$. For $i=N$, we simply have
\[
\partial_s^NS^N\hat{\mbf x}:=\partial_s^N\hat{\mbf x}_N=\ubar{\hat{\mbf x}}_{N} = \underbrace{I_{n_N}}_{U_{1,N}}\ubar{\hat{\mbf x}}_{N}.
\]
Clearly, then
{\small
\begin{align*}
&(\mcl F\hat{\mbf x})\\
&=\bmat{\partial_s^0S^0\hat{\mbf x}\\\vdots\\\partial_s^NS^N\hat{\mbf x}} = \underbrace{\bmat{U_{1,0}&&\\&\ddots&\\&&U_{1,N}}}_{U_1}\underbrace{\bmat{\ubar{\hat{\mbf x}}_{0}\\ \vdots\\ \ubar{\hat{\mbf x}}_{N}}}_{\ubar{\hat{\mbf x}}} + \underbrace{\bmat{U_{2,0}&&\\&\ddots&\\&&U_{2,(N-1)}&\\0_{n_N\times n_{1:N}}&\cdots&0_{n_N\times n_{N:N}}}}_{U_2}\underbrace{\bmat{(\mcl C\hat{\mbf x})_{1}\\\vdots\\ (\mcl C\hat{\mbf x})_{N}}}_{(\mcl C\hat{\mbf x})}
\\&=U_1\ubar{\hat{\mbf x}} + U_2(\mcl C\hat{\mbf x}).
\end{align*}
}
Finally, by Corollary~\ref{cor:x_c}, we write
\begin{align*}
(\mcl F\hat{\mbf x})(s)&=U_1\ubar{\hat{\mbf x}}(s) + U_2(\mcl C\hat{\mbf x})(s)= U_1\ubar{\hat{\mbf x}}(s) + U_2 T(s-a)(\mcl C\hat{\mbf x})(a)+\int_{a}^{s}U_2 Q(s-\theta)\ubar{\hat{\mbf x}}(\theta)d\theta\\
&=U_2 T(s-a)(\mcl C\hat{\mbf x})(a)+  U_1\ubar{\hat{\mbf x}}(s)+\int_{a}^{s}U_2 Q(s-\theta)\ubar{\hat{\mbf x}}(\theta)d\theta
\end{align*}

Now, suppose we are given $\mbf G_{\mathrm b} = \{B, B_I, B_v\}$ such that $B_T$ is invertible where
\[
B_T = B\bmat{I_{n_S}\\T(b-a)}-\int_{a}^{b}B_I(s)U_2 T(s-a)ds.
\]
For (b), 
we temporarily partition $B$ as $B= \bmat{B_l& B_r}$ where both $B_l$ and $B_r$ have equal number of columns. Then, we look at the expression $B\bmat{(\mcl C\hat{\mbf x})(a)\\ (\mcl C\hat{\mbf x})(b)}-\int_{a}^{b} B_{I}(s)(\mcl F\hat{\mbf x})(s) ds$. Clearly, we need an expression for $(\mcl C\hat{\mbf x})(b)$ which can be obtained from \cref{cor:x_c} (by substituting $s=b$) as
\[
(\mcl C\hat{\mbf x})(b) = T(b-a)(\mcl C\hat{\mbf x})(a)+\int_{a}^{b}Q(b-s)\ubar{\hat{\mbf x}}(s)ds.
\]
Replacing $(\mcl C\hat{\mbf x})(b)$ and $(\mcl F\hat{\mbf x})$ in the expression for $B\bmat{(\mcl C\hat{\mbf x})(a)\\ (\mcl C\hat{\mbf x})(b)}-\int_{a}^{b} B_{I}(s)(\mcl F\hat{\mbf x})(s) ds$, we get
{\small
\begin{align*}
&\underbrace{\bmat{B_l & B_r}}_{B}\bmat{(\mcl C\hat{\mbf x})(a)\\ (\mcl C\hat{\mbf x})(b)}-\int_{a}^{b} B_{I}(s)(\mcl F\hat{\mbf x})(s) ds \\
&= (B_l+B_rT(b-a))(\mcl C\hat{\mbf x})(a)+\int_{a}^{b}B_r Q(b-s)\ubar{\hat{\mbf x}}(s)ds\\
&-\left(\int_{a}^{b}B_I(s)U_2 T(s-a)ds\right) (\mcl C\hat{\mbf x})(a)-  \int_{a}^{b}\left( B_I(s)U_1+\int_{s}^{b}B_I(\theta) U_2 Q(\theta-s)d\theta\right)\ubar{\hat{\mbf x}}(s)ds\\
&= \underbrace{\left(B_l+B_rT(b-a)-\int_{a}^{b}B_I(s)U_2 T(s-a)ds\right)}_{B_T}(\mcl C\hat{\mbf x})(a)\\
&+\int_{a}^{b}\underbrace{\left(B_r Q(b-s)-B_I(s)U_1-\int_{s}^{b}B_I(\theta) U_2 Q(\theta-s)d\theta\right)}_{B_TB_Q(\theta)}\ubar{\hat{\mbf x}}(s)ds\\
&= B_T(\mcl C\hat{\mbf x})(a)-\int_{a}^{b}B_TB_Q(s)\ubar{\hat{\mbf x}}(s)ds,
\end{align*}
}%
which proves the second statement of the corollary.

For (c), we have the additional constraint that $\hat{\mbf x}\in X_v$. Then, we know that
\[
B\bmat{(\mcl C\hat{\mbf x})(a)\\ (\mcl C\hat{\mbf x})(b)}-\int_{a}^{b} B_{I}(s)(\mcl F\hat{\mbf x})(s) ds - B_v v = 0.
\]
Therefore, from second statement of the corollary, we have
\[
B_T(\mcl C\hat{\mbf x})(a)-\int_{a}^{b}B_TB_Q(s)\ubar{\hat{\mbf x}}(s)ds - B_v v = 0,
\]
and since $B_T$ is invertible, we can conclude that
\[
(\mcl C\hat{\mbf x})(a)=\int_{a}^{b}B_Q(s)\ubar{\hat{\mbf x}}(s)ds+B_T^{-1}B_v v.
\]

For (d), we know that $(\mcl F\hat{\mbf x})$ and $(\mcl C\hat{\mbf x})(a)$ (from steps (a) and (b)) can be expressed as
\begin{align*}
(\mcl F\hat{\mbf x})(s)&=U_2 T(s-a)(\mcl C\hat{\mbf x})(a)+  U_1\ubar{\hat{\mbf x}}(s)+\int_{a}^{s}U_2 Q(s-\theta)\ubar{\hat{\mbf x}}(\theta)d\theta,\\
(\mcl C\hat{\mbf x})(a)&=\int_{a}^{b}B_Q(s)\ubar{\hat{\mbf x}}(s)ds+B_T^{-1}B_v v.
\end{align*}
Thus, by substituting $(\mcl C\hat{\mbf x})(a)$ in the expression for $(\mcl F\hat{\mbf x})$, we get
{\small
\begin{align*}
(\mcl F\hat{\mbf x})(s) &= U_2 T(s-a)(\mcl C\hat{\mbf x})(a)+  U_1\ubar{\hat{\mbf x}}(s)+\int_{a}^{s}U_2 Q(s-\theta)\ubar{\hat{\mbf x}}(\theta)d\theta\\
&= U_2 T(s-a)\left(\int_{a}^{b}B_Q(s)\ubar{\hat{\mbf x}}(s)ds+B_T^{-1}B_v v\right)+  U_1\ubar{\hat{\mbf x}}(s)+\int_{a}^{s}U_2 Q(s-\theta)\ubar{\hat{\mbf x}}(\theta)d\theta\\
&=\fourpiFull{\emptyset}{\emptyset}{U_2T(s-a)B_T^{-1}B_v}{U_1}{U_2(T(s-a)B_Q(\theta)+Q(s-\theta))}{U_2T(s-a)B_Q(\theta)}\bmat{v\\\ubar{\hat{\mbf x}}(\cdot)}\\
&= \left(\fourpiFull{\emptyset}{\emptyset}{U_2T(s-a)B_T^{-1}B_v}{U_1}{R_{D,1}}{R_{D,2}}\bmat{v\\\ubar{\hat{\mbf{x}}}}\right)(s).
\end{align*}}
where we define the variables
\[
R_{D,1}(s,\theta) = R_{D,2}(s,\theta)+U_2Q(s-\theta), \quad R_{D,2}(,s\theta) = U_2T(s-a)B_Q(\theta).
\]

Now, since $\hat{\mbf x} \in X_{v}$ for all $t\ge 0$, by Corollary~\ref{cor:x_c}, we have
\[
(\mcl C\hat{\mbf x})(s)= T(s-a)(\mcl C\hat{\mbf x})(a)+\int_{a}^{s}Q(s-\theta)\ubar{\hat{\mbf x}}(\theta)d\theta.
\]
Furthermore, since $B_T$ is invertible, from Corollary~\ref{cor:x_D}, we have
\begin{align*}
(\mcl C\hat{\mbf x})(a) &= \int_{a}^{b}B_Q(\theta)\ubar{\hat{\mbf x}}(\theta)d\theta+B_T^{-1}B_v v,
\end{align*}
and hence we can express $(\mcl B \hat{\mbf x})$ in terms of $\hat{\mbf x}$ and $v$ as
\begin{align*}
&(\mcl B \hat{\mbf x}) = \bmat{(\mcl C\hat{\mbf x})(a)\\(\mcl C\hat{\mbf x})(b)} = \bmat{I\\T(b-a)}(\mcl C\hat{\mbf x})(a) + \int_a^b \bmat{0\\Q(b-\theta)}\ubar{\hat{\mbf x}}(\theta)d\theta\\
&= \bmat{I\\T(b-a)}\left(\int_{a}^{b}B_Q(\theta)\ubar{\hat{\mbf x}}(\theta)d\theta+B_T^{-1}B_v v\right) + \int_a^b \bmat{0\\Q(b-\theta)}\ubar{\hat{\mbf x}}(\theta)d\theta\\	&=\fourpi{\bmat{I\\T(b-a)}B_T^{-1}B_v}{\bmat{I\\T(b-a)}B_Q+\bmat{0\\Q}}{\emptyset}{\emptyset}\bmat{v\\\ubar{\hat{\mbf x}}}.
\end{align*}
To get an expression for the combined $v$, $(\mcl B \hat{\mbf x})$ and $(\mcl F\hat{\mbf x})$, we can just concatenate them vertically to get
\[
\bmat{\bmat{v\\(\mcl B \hat{\mbf x})}\\(\mcl F\hat{\mbf x})(\cdot)}= \fourpi{\bmat{I_{n_v}\\B_T^{-1}B_v\\T(b-a)B_T^{-1}B_v}}{\bmat{0_{n_r \times n_x}\\B_Q(s)\\T(b-a)B_Q(s)+Q(b-s)}}{U_2T(s-a)B_T^{-1}B_v}{U_1,R_{D,1},R_{D,2}}\bmat{v\\\ubar{\hat{\mbf x}}}.
\]
\end{proof}

Now, from a), we have a map from $\{(\mcl C\hat{\mbf x})(a),~\ubar{\hat{\mbf x}}\}$ to the vector of all well-defined terms, $\mcl F\hat{\mbf x}$. Furthermore, from c), when the BCs are PIE-compatible we have a map from $\{\ubar{\hat{\mbf x}},~v\}$ to $(\mcl C\hat{\mbf x})(a)$. This allows us to express the left boundary values, $(\mcl C\hat{\mbf x})(a)$ in terms of $\{\ubar{\hat{\mbf x}},~v\}$ -- yielding a map from $\ubar{\hat{\mbf x}}$ to $(\mcl F\hat{\mbf x})$. Extending this result, we can use~\cref{cor:ftc_primary} to obtain a map from $\{\ubar{\hat{\mbf x}},~v\}$ to $\hat{\mbf x}$.

\Tmap*

\begin{proof} \textbf{Proof of Part 1.} Let $\hat{\mbf x}\in X_v$ for some $v\in \R^{q}$. Clearly, by definition of $X_v$, $\partial_s^i\hat{\mbf x}_i\in L_2^{n_i}$. Therefore, $\mcl D\hat{\mbf x}\in L_2^{n_{\hat{\mbf x}}}$. Next we need to express $\hat{\mbf x}$ in terms of $\ubar{\hat{\mbf x}}:=\mcl D\hat{\mbf x}$ and $v$. For that, we will first express $(\mcl C\hat{\mbf x})(a)$ solely in terms of $\ubar{\hat{\mbf x}}$ and $v$. From \cref{cor:x_D}, we know that if $\{n,\mbf G_{\mathrm b}\}$ is PIE-compatible, then
\[
(\mcl C\hat{\mbf x})(a)= \int_{a}^{b}B_Q(\theta)\ubar{\hat{\mbf x}}(\theta)d\theta+B_T^{-1}B_v v.
\]
Now that we have an expression for $(\mcl C\hat{\mbf x})(a)$, we simply substitute this into the expression for $\hat{\mbf x}$ from Corollary~\ref{cor:ftc_primary} to obtain
\begin{align*}
&\hat{\mbf x}_{1:N}(s)=T_1(s-a)(\mcl C\hat{\mbf x})(a)+\int_{a}^{s}Q_1(s-\theta)\ubar{\hat{\mbf x}}(\theta)d\theta\\
&=\int_{a}^{b}T_1(s-a)B_Q(\theta)\ubar{\hat{\mbf x}}(\theta)d\theta+ \int_{a}^{s}Q_1(s-\theta)\ubar{\hat{\mbf x}}(\theta)d\theta+T_1(s-a)B_T^{-1}B_v v.
\end{align*}
Adding on the somewhat incongruous $\hat{\mbf x}_0$ term, we obtain
\begin{align}
\hat{\mbf x}(s)=\bmat{\hat{\mbf x}_0(s)\\\hat{\mbf x}_{1:N}(s)}&=\underbrace{\bmat{I_{n_0}&0\\0&0_{n_x-n_0}}}_{G_0}\ubar{\hat{\mbf x}}(s)+\int_{a}^{b}\underbrace{\bmat{0_{n_0 \times n_{\hat{\mbf x}}}\\T_1(s-a)B_Q(\theta)}}_{G_2(s,\theta)}\ubar{\hat{\mbf x}}(\theta)d\theta+ \int_{a}^{s}\underbrace{\bmat{0_{n_0 \times n_{\hat{\mbf x}}}\\Q_1(s-\theta)}}_{G_1(s,\theta)-G_2(s,\theta)}\ubar{\hat{\mbf x}}(\theta)d\theta\label{eqn:ThmA6_p1a}\\
&\quad +\underbrace{\bmat{0_{n_0 \times n_v}\\T_1(s-a)B_T^{-1}B_v}}_{G_v(s)} v\notag\\
&=G_0\ubar{\hat{\mbf x}}(s)
+\int_{s}^{b}G_2(s,\theta)\ubar{\hat{\mbf x}}(\theta)d\theta+ \int_{a}^{s}G_1(s,\theta)\ubar{\hat{\mbf x}}(\theta)d\theta+G_v(s) v= (\hat{\mcl T}\ubar{\hat{\mbf x}})(s) + (\mcl T_v v) (s).\notag
\end{align}
\end{proof}
\begin{proof} \textbf{Proof of Part 2.}
Let $v\in \R^q$ and $\ubar{\hat{\mbf x}} \in L_2^{n_{\hat{\mbf x}}}$ be arbitrary. Our first step is to show that $\mcl D\left(\hat{\mcl{T}}  \ubar{\hat{\mbf{x}}}+\mcl{T}_{v}v\right)=\ubar{\hat{\mbf x}} \in L_2^{n_{\hat{\mbf x}}}$. By the definition of $\hat{\mcl T}$ and $\mcl T_v$, Eqn~\eqref{eqn:ThmA6_p1a} at the end of the proof of Part 1 shows that for any $\ubar{\hat{\mbf x}} \in L_2^{n_{\hat{\mbf x}}}$ and $v \in \R^{q}$,
\begin{align*}
&(\hat{\mcl T}\ubar{\hat{\mbf x}})(s) + (\mcl T_v v) (s)\\
&=\bmat{I&0\\0&0}\ubar{\hat{\mbf x}}(s)
+\int_{a}^{b}\bmat{0\\T_1(s-a)B_Q(\theta)}\ubar{\hat{\mbf x}}(\theta)d\theta+ \int_{a}^{s}\bmat{0\\Q_1(s-\theta)}\ubar{\hat{\mbf x}}(\theta)d\theta\\
&\quad+\bmat{0\\T_1(s-a)B_T^{-1}B_v} v.
\end{align*}
Thus, we may group the terms with $T_1(s-a)$ together and apply the $\mcl D$ operator to obtain
\begin{align*}
\mcl D\left(\hat{\mcl T}\ubar{\hat{\mbf x}} + \mcl T_v v\right)(s)&=\mcl D\bmat{0\\T_1(s-a)}\left(\int_{a}^{b} B_Q(\theta)\ubar{\hat{\mbf x}}(\theta)d\theta+B_T^{-1}B_v v\right)\\
&\quad+ \mcl D\left(\bmat{I_{n_0}&0\\0&0}\ubar{\hat{\mbf x}}(s)+\int_{a}^{s}\bmat{0_{n_0 \times n_{\hat{\mbf x}}}\\Q_1(s-\theta)}\ubar{\hat{\mbf x}}(\theta)d\theta\right).
\end{align*}
Now, examining the first term we have
\begin{align*}
\mcl D\bmat{0\\T_1(s-a)}&=\bmat{I_{n_0}&&\\&\ddots&\\&&\partial_s^N I_{n_N}}\bmat{0\\T_1(s-a)}=\bmat{0\\\left(\bmat{\partial_sI_{n_1}&&\\&\ddots&\\&&\partial_s^N I_{n_N}}T_1(s-a)\right)}
\end{align*}
and since $\partial_s^i \tau_j(s)=0$ for any $j>i$, we have
{\small
\begin{align*}
\left(\bmat{\partial_sI_{n_1}&&\\&\ddots&\\&&\partial_s^N I_{n_N}}T_1(s-a)\right)&=\bmat{\partial_sI_{n_1}&&\\&\ddots&\\&&\partial_s^N I_{n_N}}\bmat{\tau_0(s-a)J_{1,1}&\tau_1(s-a)J_{1,2}&\cdots&\tau_{N-1}(s-a)J_{1,N}}\\
&=\bmat{\partial_sI_{n_1}&&\\&\ddots&\\&&\partial_s^N I_{n_N}}\bmat{\tau_0(s-a)I_{n_{1:N}}&\bmat{0\\\tau_1(s-a)I_{n_{2:N}}}&\cdots&\bmat{0\\\tau_{N-1}(s-a)I_{n_{N}}}}\\
&=0.
\end{align*}}
Hence the first term in our expression for $\mcl D\left(\hat{\mcl T}\ubar{\hat{\mbf x}} + \mcl T_v v\right)$ is zero. Now, consider the second term in the expression for $\mcl D \left(\hat{\mcl T}\ubar{\hat{\mbf x}} + \mcl T_v v\right)$,
\begin{align*}
&\mcl D\left(\bmat{I_{n_0}&0\\0&0}\ubar{\hat{\mbf x}}(s)+\int_{a}^{s}\bmat{0_{n_0 \times n_x}\\Q_1(s-\theta)}\ubar{\hat{\mbf x}}(\theta)d\theta\right)\\
&=\mcl D\left(\bmat{I_{n_0}&0\\0&0}\ubar{\hat{\mbf x}}(s)+\int_{a}^{s}\bmat{0_{n_0}&0&&&0\\0&\tau_0(s)I_{n_1}&&&\\ \vdots&&&\ddots&\\0&&&&\tau_{N-1}(s)I_{n_N}}\ubar{\hat{\mbf x}}(\theta)d\theta\right).
\end{align*}
For this term, we use an inductive approach. Specifically, we factor $\mcl D$ into first-order derivative operators as
\[
\mcl D=\bmat{I&&&\\&\partial_s I&&\\&&\ddots&\\&&&\partial_s^N I}=\prod_{i=1}^N \underbrace{\bmat{I_{n_{0:i-1}}&0\\0&I_{n_{i:N}}\partial_s}}_{\mcl D_i}.
\]
Now, since $\partial_s \tau_i(s)=\tau_{i-1}(s)$ for $i\ge 1$, $\tau_i(0)=0$ for $i>0$ and $\tau_0(0)=1$, we have that for any $i<N$,
\begin{align*}
&\mcl D_i\bmat{0_{n_{0:i-1} \times n_x}\\Q_i(s-\theta)}=\bmat{I_{n_{0:i-1}}&0\\0&I_{n_{i:N}}\partial_s}\bmat{0_{n_{0:i-1}\times n_{0:i-1}}&&&\\&\tau_0(s)I_{n_i}&&\\ &&\ddots&\\&&&\tau_{N-i}(s)I_{n_N}}\\
&=\bmat{0_{n_{0:i}\times n_{0:i}}&&&\\&\tau_0(s)I_{n_{i+1}}&&\\ &&\ddots&\\&&&\tau_{N-i-1}(s)I_{n_N}}=\bmat{0_{n_{0:i} \times n_x}\\Q_{i+1}(s-\theta)}
\end{align*}
and $\mcl D_N\bmat{0_{n_{0:N-1} \times n_x}\\Q_N(s-\theta)}=0$. Additionally,
for $i\ge 0$, we have
\[
\bmat{0_{n_{0:i-1} \times n_x}\\Q_i(0)}=\bmat{0_{n_{0:i-1}\times n_{0:i-1}}&&\\&I_{n_i}&\\&&0_{n_{i+1:N}\times n_{i+1:N}}}.
\]
We conclude that
\begin{align*}
&\mcl D_i\left(\bmat{I_{n_{0:i-1}}&\\&0}\ubar{\hat{\mbf x}}(s)+\int_{a}^{s}\bmat{0_{n_{0:i-1} \times n_x}\\Q_i(s-\theta)}\ubar{\hat{\mbf x}}(\theta)d\theta\right)\\
&=\bmat{I_{n_{0:i-1}}&\\&0}\ubar{\hat{\mbf x}}(s)+\bmat{0_{n_{0:i-1} \times n_x}\\Q_i(0)}\ubar{\hat{\mbf x}}(s)+\int_{a}^{s}\bmat{0_{n_{0:i} \times n_x}\\Q_{i+1}(s-\theta)}\ubar{\hat{\mbf x}}(\theta)d\theta\\
&=\bmat{I_{n_{0:i}}&\\&0}\ubar{\hat{\mbf x}}(s)+\int_{a}^{s}\bmat{0_{n_{0:i} \times n_x}\\Q_{i+1}(s-\theta)}\ubar{\hat{\mbf x}}(\theta)d\theta.
\end{align*}
Applying this inductive step to each of the $\mcl D_i$ operators in $\mcl D\hat{\mbf x}$, we have
\begin{align*}
&\mcl D\left(\bmat{I_{n_0}&0\\0&0}\ubar{\hat{\mbf x}}(s)+\int_{a}^{s}\bmat{0\\Q_1(s-\theta)}\ubar{\hat{\mbf x}}(\theta)d\theta\right)\\
&= \mcl D_{N}\cdots \mcl D_{1}\left(\bmat{I_{n_0}&0\\0&0}\ubar{\hat{\mbf x}}(s)+\int_{a}^{s}\bmat{0\\Q_1(s-\theta)}\ubar{\hat{\mbf x}}(\theta)d\theta\right)=\ubar{\hat{\mbf x}}(s).
\end{align*}
Combining these results, we conclude that for any $\ubar{\hat{\mbf x}} \in L_{2}^{n_{\hat{\mbf x}}}$
{\small
\begin{align*}
\left(\mcl D\left(\hat{\mcl T}\mbf{\hat x}+\mcl T_v v\right)\right)(s)=
&=\mcl D\left(\int_{a}^{b} G_2(s,\theta)\ubar{\hat{\mbf x}}(\theta)d\theta+G_v(s) v\right)+ \mcl D\left(\bmat{I_{n_0}&0\\0&0}\ubar{\hat{\mbf x}}(s)+\int_{a}^{s}\bmat{0\\Q_1(s-\theta)}\ubar{\hat{\mbf x}}(\theta)d\theta\right)=\ubar{\hat{\mbf x}}(s).
\end{align*}
}%

Finally, we need to show that for any $\ubar{\hat{\mbf x}} \in L_{2}^{n_{\hat{\mbf x}}}$, $\hat{\mcl{T}}  \ubar{\hat{\mbf{x}}}+\mcl{T}_{v}v \in X_v$. Let $\hat{\mbf x}:=\hat{\mcl{T}}  \ubar{\hat{\mbf{x}}}+\mcl{T}_{v}v$. Clearly, since $\mcl D\hat{\mbf x}=\ubar{\hat{\mbf x}} \in L_2^{n_{\hat{\mbf x}}}$, we have $\hat{\mbf x}\in W^n$. To show that $\hat{\mbf x}\in X_v$, however, we must now show that the BCs are satisfied. For this part, we have that if $\hat{\mbf x}=\hat{\mcl{T}}  \ubar{\hat{\mbf{x}}}+\mcl{T}_{v}v\in W^n$, then by Corollary~\ref{cor:x_D},
{\small
\begin{align*}
&B\bmat{(\mcl C\hat{\mbf x})(a)\\ (\mcl C\hat{\mbf x})(b)}-\int_{a}^{b} B_{I}(s)(\mcl F\hat{\mbf x})(s) ds -B_v v= B_T\left((\mcl C\hat{\mbf x})(a)-\int_{a}^{b}B_Q(s)\ubar{\hat{\mbf{x}}}(s)ds-B_T^{-1}B_v v\right).
\end{align*}
}
Since $\hat{\mbf x}\in W^n$ and $B_T$ is invertible, we have that $\hat{\mbf x} \in X_v$ if and only if
\[
(\mcl C\hat{\mbf x})(a)-\int_{a}^{b}B_Q(s)\ubar{\hat{\mbf{x}}}(s)ds-B_T^{-1}B_v v=0.
\]
Recall from the beginning of the proof of Part 2 that
\begin{align*}
\hat{\mbf x}(s):=\bmat{\hat{\mbf x}_0(s)\\\hat{\mbf x}_{1:N}(s)}&=\bmat{I&0\\0&0}\ubar{\hat{\mbf x}}(s)+\int_{a}^{b}\bmat{0\\T_1(s-a)B_Q(\theta)}\ubar{\hat{\mbf x}}(\theta)d\theta+ \int_{a}^{s}\bmat{0\\Q_1(s-\theta)}\ubar{\hat{\mbf x}}(\theta)d\theta\\
&\quad +\bmat{0\\T_1(s-a)B_T^{-1}B_v} v
\end{align*}
and hence
\begin{align*}
&\hat{\mbf x}_{1:N}(s)=\int_{a}^{b}T_1(s-a)B_Q(\theta)\ubar{\hat{\mbf x}}(\theta)d\theta+ \int_{a}^{s}Q_1(s-\theta)\ubar{\hat{\mbf x}}(\theta)d\theta+T_1(s-a)B_T^{-1}B_v v.
\end{align*}
In addition, from Corollary A.3, we have
\[
\hat{\mbf x}_{1:N}(s)=T_1(s-a)(\mcl C\hat{\mbf x})(a)+\int_{a}^{s}Q_1(s-\theta)\ubar{\hat{\mbf x}}(\theta)d\theta.
\]
Substituting this identity in the previous equation, we get
\begin{align*}
&T_1(s-a)(\mcl C\hat{\mbf x})(a)+\int_{a}^{s}Q_1(s-\theta)\ubar{\hat{\mbf x}}(\theta)d\theta\\
&=\int_{a}^{b}T_1(s-a)B_Q(\theta)\ubar{\hat{\mbf x}}(\theta)d\theta+ \int_{a}^{s}Q_1(s-\theta)\ubar{\hat{\mbf x}}(\theta)d\theta+T_1(s-a)B_T^{-1}B_v v,
\end{align*}
which implies
\begin{align*}
T_1(s-a)\left((\mcl C\hat{\mbf x})(a)-\int_{a}^{b}B_Q(\theta)\ubar{\hat{\mbf x}}(\theta)d\theta-B_T^{-1}B_v v\right)=0.
\end{align*}
We will use induction to show that the above equality holds when $T_1$ is replaced by $T_i$. First, suppose
\begin{align*}
T_i(s-a)\left((\mcl C\hat{\mbf x})(a)-\int_{a}^{b}B_Q(\theta)\ubar{\hat{\mbf x}}(\theta)d\theta-B_T^{-1}B_v v\right)=0.
\end{align*}
Then, since
\begin{align*}
&\partial_s \bmat{0_{n_{i+1:N}\times n_i}&I_{n_{i+1:N}}} T_i(s-a) = T_{i+1}(s-a), \\
\end{align*}
we have the relation
\begin{align*}
&\partial_s \bmat{0_{n_{i+1:N}\times n_i}&I_{n_{i+1:N}}} T_i(s-a)\left((\mcl C\hat{\mbf x})(a)-\int_{a}^{b}B_Q(\theta)\ubar{\hat{\mbf x}}(\theta)d\theta-B_T^{-1}B_v v\right)\\
&=T_{i+1}(s-a)\left((\mcl C\hat{\mbf x})(a)-\int_{a}^{b}B_Q(\theta)\ubar{\hat{\mbf x}}(\theta)d\theta-B_T^{-1}B_v v\right)=0.
\end{align*}
Since the equality is true for $i=1$, by induction we can conclude, for any $i\ge 1$,
\begin{align*}
T_i(s-a)\left((\mcl C\hat{\mbf x})(a)-\int_{a}^{b}B_Q(\theta)\ubar{\hat{\mbf x}}(\theta)d\theta-B_T^{-1}B_v v\right)=0.
\end{align*}
By stacking all $T_i$'s and using $T = \text{col}(T_1,\cdots,T_N)$, for any $s\in[a,b]$ we have
\begin{align*}
T(s-a)\left((\mcl C\hat{\mbf x})(a)-\int_{a}^{b}B_Q(\theta)\ubar{\hat{\mbf x}}(\theta)d\theta-B_T^{-1}B_v v\right)=0.
\end{align*}
and since $T(0)=I_{n_S}$, we have that
\begin{align*}
(\mcl C\hat{\mbf x})(a)-\int_{a}^{b}B_Q(\theta)\ubar{\hat{\mbf x}}(\theta)d\theta-B_T^{-1}B_v v=0.
\end{align*}
which completes the proof.
\end{proof}

\subsection{Equivalence of PIE and PDE subsystems: Proof of \Cref{thm:equivalence}}\label{app:equivalence}
Now that we have established a PI map from $L_2$ to $X_v$, we will obtain the PIE associated with a PDE subsystem by replacing $\hat{\mbf x}$ in the PDE subsystem with $\hat{\mbf x}= \hat{\mcl T}\ubar{\hat{\mbf x}}+\mcl T_v v$. Because we have shown that this PI map is a bijection, we will then conclude that existence of a solution for the PIE subsystem guarantees the existence of a solution for the PDE subsystem. This proof is split into two parts.

\thmEquivalence*
\begin{proof}\textbf{Proof of 1).}
Suppose $\{\hat{\mbf x},~r\}$ satisfies the PDE~\Cref{eq:general_pde_subsystem} defined by $n\in \N^{N+1}$ and $\{\mbf G_{\mathrm b}, \mbf G_{\mathrm{p}}\}$ with initial conditions $\hat{\mbf x}^0$ and input $v$. Then by Definition~\ref{defn:PDE}: a)	$r\in L_{2e}^{n_r}[\R_+]$; b) $\hat{\mbf x}(t) \in X_{v(t)}$ for all $t\ge 0$; c) $\hat{\mbf x}$ is Frech\'et differentiable with respect to the $L_2$-norm  on $\R_+$; d) \Cref{eq:general_pde_subsystem} is satisfied for all $t\ge 0$; and e) $\hat{\mbf x}(0) = \hat{\mbf x}^0$.

Let $\ubar{\hat{\mbf x}} = \mcl{D}\hat{\mbf x}$, $\ubar{\hat{\mbf x}}^0 = \mcl{D}\hat{\mbf x}^0$, $n=n_{\hat{\mbf x}}$ and $m=0$. Our goal is to show that for $\mbf G_{\mathrm{PIE}}$ as defined above, $\{\ubar{\hat{\mbf x}},r,\emptyset\}$ satisfies the PIE defined by $\mbf G_{PIE}$ for initial condition $\ubar{\hat{\mbf x}}^0$ and input $\{v,\emptyset\}$. For this, we must show that: 1) $v\in L_{2e}^{n_v}[\R_+]$ and $(\mcl{T}_{v} v)(\cdot,s)\in W_{1e}^{n_{\hat{\mbf x}}}[\R_+]$ for all $s \in [a,b]$; 2) $\ubar{\hat{\mbf x}}:\R_+\to \R L_2^{0,n_{\hat{\mbf x}}}[a,b]$ and $r\in L_{2e}^{n_r}[\R_+]$; 3) $\ubar{\hat{\mbf x}}^0\in \R L_2^{0,n_{\hat{\mbf x}}}[a,b]$ and $\ubar{\hat{\mbf x}}(0) = \ubar{\hat{\mbf x}}^0$; 4) $\hat{\mcl T}\ubar{\hat{\mbf x}}$ is Frech\'et differentiable on $\R_+$; and 5) Equation~\eqref{eq:PIE_full} is satisfied for all $t\in\R_+$.

\indent For 1), $v\in L_{2e}^{n_v}[\R_+]$ from the theorem statement and by the definition of $\mcl T_w$, $B_v v\in W_{1e}^{2n_S}[\R_+]$ implies
\[
(\mcl{T}_{w} w)(s) = (\mcl T_v v)(s) = \bmat{0\\T_1(s-a)}B_T^{-1}B_v v\in W_{1e}^{n_{\hat{\mbf x}}}[\R_+].
\]

For 2), from \Cref{thm:T_map}a we have that $\hat{\mbf x}(t) \in X_{v(t)}$ implies $\ubar{\hat{\mbf x}}(t)=\mcl{D}\hat{\mbf x}(t) \in RL_2^{0,n_{\hat{\mbf x}}}=L_2^{n_{\hat{\mbf x}}}$ for all $t \ge 0$. Furthermore, from the definition of solution of the PDE, $r\in L_{2e}^{n_r}[\R_+]$.

For 3), from \Cref{thm:T_map}a we have that $\hat{\mbf x}^0 \in X_{v(0)}$ implies $\ubar{\hat{\mbf x}}^0=\mcl{D}\hat{\mbf x}^0 \in RL_2^{0,n_{\hat{\mbf x}}}=L_2^{n_{\hat{\mbf x}}}$. Furthermore, since $\ubar{\hat{\mbf x}}(t) = \mcl{D}\hat{\mbf x}(t)$ for all $t\ge 0$, we have $\ubar{\hat{\mbf x}}(0) = \mcl{D}\hat{\mbf x}(0)=\mcl D \hat{\mbf x}^0 = \ubar{\mbf x}^0$.

For 4), because $\hat{\mbf x}$ is Frech\'et differentiable  on $\R_+$, the limit of $\frac{\hat{\mbf x}(t+h)-\hat{\mbf x}(t)}{h}$ as $h\to0^+$ exists for any $t\ge 0$ when convergence is defined with respect to the $L_2$ norm. This, and the fact that $\mcl T_v v \in W_{1e}^{n_v}$ implies  that
\[
\lim_{h\rightarrow 0^+}\frac{\hat{\mcl T}\ubar{\hat{\mbf x}}(t+h)-\hat{\mcl T}\ubar{\hat{\mbf x}}(t) }{h}=\lim_{h\rightarrow 0^+}\frac{\hat{\mbf x}(t+h)-\hat{\mbf x}(t) }{h}-\lim_{h\rightarrow 0^+}\frac{\mcl T_v v(t+h)-\mcl T_v v(t) }{h}
\]
similarly exists for all $t\ge 0$. Thus, $\hat{\mcl T} \ubar{\hat{\mbf x}}$ is Frech\'et differentiable with respect to $L_2$-norm.

Lastly, for 5), since $\hat{\mbf x}(t)$ satisfies \eqref{eq:odepde_b}-\eqref{eq:general_pde_subsystem} for all $t\ge 0$, we have
\begin{flalign}
\bmat{\dot{\hat{\mbf x}}(t,s)\\ r(t)} &= \sum\limits_{i=0}^N \bmat{A_{0}(s) +\int\limits_{a}^{s}A_{1}(s,\cdot)+\int\limits_{s}^{b} A_{2}(s,\cdot)\\
\int_{a}^{b} C_{r}(\cdot)	}(\mcl F\hat{\mbf x})(t,\cdot) +\bmat{B_{xv}(s) & B_{xb}(s)\\ 0& D_{rb} }\bmat{v(t)\\ (\mcl B \hat{\mbf x})(t)}.\label{eq:general_pde_subsystem2}
\end{flalign}
Since $\hat{\mbf x}(t)\in X_{v(t)}$ and $\ubar{\hat{\mbf x}}(t)=\mcl D \hat{\mbf x}(t)$ for all $t\ge 0$, from \Cref{thm:T_map}, we have that
\begin{align*}
\hat{\mbf x}(t) &= \hat{\mcl{T}}\ubar{\hat{\mbf x}}(t)+\mcl{T}_{v}v(t)\qquad \text{which implies}\qquad  \dot{\hat{\mbf x}}(t) = \hat{\mcl{T}}\ubar{\dot{\hat{\mbf{x}}}}(t)+\mcl{T}_{v}\dot{v}(t).
\end{align*}
We can substitute this into Eq.~\ref{eq:general_pde_subsystem2} and re-write Eq.~\ref{eq:general_pde_subsystem2} using the PI operator notation to get the compact relation
\begin{align}\label{eq:pie_temp}
\bmat{r(t)\\\hat{\mcl{T}}\ubar{\dot{\hat{\mbf{x}}}}(t)+\mcl{T}_{v}\dot{v}(t)} &= \fourpi{\bmat{0&D_{rb}}}{C_r}{\bmat{B_{xv}&B_{xb}}}{A_i}\bmat{\bmat{v(t)\\(\mcl B \hat{\mbf x})(t)}\\ (\mcl F\hat{\mbf x})(t)}.
\end{align}
where we define
\begin{align*}
&(\mcl F\hat{\mbf x})(t):= U_1\ubar{\hat{\mbf x}}(t) + U_2(\mcl C\hat{\mbf x})(t), \qquad (\mcl B \hat{\mbf x})(t) := \bmat{(\mcl C\hat{\mbf x})(t,a)\\(\mcl C\hat{\mbf x})(t,b)}\qquad (\mcl C\hat{\mbf x})(t):= \bmat{(S\hat{\mbf x})(t)\\ (\partial_s S^2\hat{\mbf x})(t) \\ \vdots \\ (\partial_s^{N-1}S^{N}\hat{\mbf x})(t)}.
\end{align*}

We know from \Cref{cor:x_D}d that when $B_T$ is invertible, $\hat{\mbf x}(t)\in X_{v(t)}$ and $\ubar{\hat{\mbf x}}(t) = \mcl D\hat{\mbf x}(t)$, we have the relation
\[
\bmat{\bmat{v(t)\\(\mcl B \hat{\mbf x})(t)}\\(\mcl F\hat{\mbf x})(t)}=\fourpi{\bmat{I\\P_{b}}}{\bmat{0\\Q_{b}}}{U_2T(s-a)B_T^{-1}B_v}{\{U_1,R_{D,1},R_{D,2}\}}\bmat{v(t)\\\ubar{\hat{\mbf x}}(t)}.
\]
Using the above expression for $\bmat{v(t)\\(\mcl B \hat{\mbf x})(t)\\(\mcl F\hat{\mbf x})(t)}$, we now expand Eq.~\ref{eq:pie_temp} to obtain
\begin{align*}
&	\bmat{r(t)\\\hat{\mcl{T}}\ubar{\dot{\hat{\mbf{x}}}}(t)+\mcl{T}_{v}\dot{v}(t)} = \fourpi{\bmat{0&D_{rb}}}{C_r}{\bmat{B_{xv}&B_{xb}}}{A_i}\bmat{\bmat{v(t)\\(\mcl B \hat{\mbf x})(t)}\\ (\mcl F\hat{\mbf x})(t)}\\
&=\fourpi{\bmat{0&D_{rb}}}{C_r}{\bmat{B_{xv}&B_{xb}}}{A_i}\fourpi{\bmat{I\\P_{b}}}{\bmat{0\\Q_{b}}}{U_2T(s-a)B_T^{-1}B_v}{\{U_1,R_{D,1},R_{D,2}\}}\bmat{v(t)\\\ubar{\hat{\mbf x}}(t)}=\fourpi{ D_{rv}}{ C_{rx}}{ B_{xv}}{\hat A_i}\bmat{v(t)\\\ubar{\hat{\mbf x}}(t)}
\end{align*}
where
\begin{align*}
\bmat{D_{rv}& C_{rx}\\ B_{xv}&\hat A_i}=\picomp{\bmat{0&D_{rb}}}{C_r}{\bmat{B_{xv}&B_{xb}}}{A_i}{\bmat{I_{n_v}\\B_T^{-1}B_v\\T(b-a)B_T^{-1}B_v}}{\bmat{0_{n_r \times n_x}\\B_Q(s)\\T(b-a)B_Q(s)+Q(b-s)}}{U_2T(s-a)B_T^{-1}B_v}{U_1,R_{D,1},R_{D,2}}.
\end{align*}
which shows that $\{\ubar{\hat{\mbf x}},r,\emptyset\}$ satisfies the PIE defined by $\mbf G_{\mathrm{PIE}}$ for initial condition $\ubar{\hat{\mbf x}}^0$ and input $\{v,\emptyset\}$.
\end{proof}

\begin{proof} \textbf{Proof of 2).}
Suppose $\{\ubar{\hat{\mbf x}},r,\emptyset\}$ satisfies the PIE ~\Cref{eq:PIE_full} defined by the set of parameters $\mbf G_{\mathrm{PIE}}$ for initial conditions $\ubar{\hat{\mbf x}}^0$ and input $\{v,\emptyset\}$. Then we have: a) $r\in L_{2e}^{n_r}[\R_+]$; b) $\ubar{\hat{\mbf x}}(t,\cdot)\in \R L_2^{m,n}[a,b]$ for all $t\ge 0$; c) $\hat{\mcl T}\ubar{\hat{\mbf x}}$ is Frech\'et differentiable on $\R_+$; d) \Cref{eq:PIE_full} is satisfied for all $t\in\R_+$; and e) $\ubar{\hat{\mbf x}}(0,\cdot) = \ubar{\hat{\mbf x}}^0$.
Let
\begin{align*}
\hat{\mbf x}(t) = \hat{\mcl{T}}\ubar{\hat{\mbf x}}(t)+\mcl{T}_{v}v(t), \quad \hat{\mbf x}^0 = \hat{\mcl{T}}\ubar{\hat{\mbf x}}^0+\mcl{T}_v v(0).
\end{align*}
Then, our goal is to show that, $\{\hat{\mcl T}\ubar{\hat{\mbf x}}+\mcl T_v v,~r\}$ satisfies the PDE~\Cref{eq:general_pde_subsystem} defined by $n\in \N^{N+1}$ and $\{\mbf G_{\mathrm b}, \mbf G_{\mathrm{p}}\}$ with initial conditions $\hat{\mbf x}^0 = \hat{\mcl T}\ubar{\hat{\mbf x}}^0 + \mcl T_v v(0)$ and input $v$. For this, we must show: 1) $r\in L_{2e}^{n_r}[\R_+]$; 2) $\hat{\mbf x}(t) \in X_{v(t)}$ for all $t\ge 0$; 3) $\hat{\mbf x}^0 \in X_{v(0)}$ and $\hat{\mbf x}(0,\cdot) = \hat{\mbf x}^0$; 4) $\hat{\mbf x}$ is Frech\'et differentiable with respect to the $L_2$-norm  on $\R_+$; and 5) \Cref{eq:general_pde_subsystem} is satisfied for all $t\ge 0$.

For 1), $r\in L_{2e}^{n_r}[\R_+]$ holds immediately by the definition of solution of the PIE.

For 2), \Cref{thm:T_map}b states that for any $v(t)\in \R$ and $\ubar{\hat{\mbf x}}(t)\in L_2^{n_{\hat{\mbf x}}}$, we have $\hat{\mbf x}(t)=\hat{\mcl{T}}\ubar{\hat{\mbf x}}(t)+\mcl{T}_{v}v(t)\in X_{v(t)}$.

For 3), \Cref{thm:T_map}b states that for any $v(0)\in \R$ and $\ubar{\hat{\mbf x}}^0\in L_2^{n_{\hat{\mbf x}}}$, we have $\hat{\mbf x}^0 = \hat{\mcl T}\ubar{\hat{\mbf x}}^0 + \mcl T_v v(0)\in X_{v(0)}$. In addition, $\hat{\mbf x}(0) = \hat{\mcl{T}}\ubar{\hat{\mbf x}}(0)+\mcl{T}_{v}v(0)=\hat{\mcl{T}}\ubar{\hat{\mbf x}}^0+\mcl{T}_{v}v(0)\in X_{v(0)}$.

For 4), we know $\hat{\mcl T}\ubar{\hat{\mbf x}}$ is Frech\'et differentiable, which implies that $\lim_{h\rightarrow 0^+}\frac{\mcl T\ubar{\hat{\mbf x}}(t+h)-\mcl T\ubar{\hat{\mbf x}}(t) }{h}$ exists when the convergence is with respect to $L_2$-norm. Since $\mcl T_v v \in W_{1e}^{n_v}$, we conclude that
\[
\lim_{h\rightarrow 0^+}\frac{\hat{\mbf x}(t+h)-\hat{\mbf x}(t) }{h}=\lim_{h\rightarrow 0^+}\frac{\hat{\mcl T}\ubar{\hat{\mbf x}}(t+h)-\hat{\mcl T}\ubar{\hat{\mbf x}}(t) }{h}+\lim_{h\rightarrow 0^+}\frac{\mcl T_v v(t+h)-\mcl T_v v(t) }{h}
\]
exists for all $t\ge 0$. Thus, $\hat{\mbf x}$ is Frech\'et differentiable with respect to $L_2$-norm.

For 5), since $\hat{\mbf x}$ is Frech\'et differentiable and $\ubar{\hat{\mbf x}}$ satisfies the PIE, we have
\begin{align*}
\dot{\hat{\mbf x}}(t) &= \hat{\mcl{T}}\ubar{\dot{\hat{\mbf x}}}(t)+\mcl{T}_{v}\dot v(t)= \hat{\mcl{A}}\ubar{\hat{\mbf x}}(t)+\mcl{B}_{v}v(t)
\end{align*}
and furthermore, $r(t) =  \mcl C_{r}\ubar{\hat{\mbf x}}(t)+\mcl{D}_{rv}v(t)$. Combining these expressions, we obtain
\begin{align*}
&\bmat{r(t)\\\dot{\hat{\mbf x}}(t)} = \bmat{\mcl{D}_{rv}&\mcl C_{r}\\\mcl{B}_{v}&\hat{\mcl{A}}}\bmat{v(t)\\\ubar{\hat{\mbf x}}(t)}=\fourpi{D_{rv}}{C_{rx}}{B_{xv}}{\hat A_i}\bmat{v(t)\\\ubar{\hat{\mbf x}}(t)}.
\end{align*}
Now, we use the relation from~\cref{fig:PIE_subsystem_equation}
{\small
\begin{align*}
\fourpi{D_{rv}}{C_{rx}}{B_{xv}}{\hat A_i}&=\mcl P\left[\picomp{\bmat{0&D_{rb}}}{C_r}{\bmat{B_{xv}&B_{xb}}}{A_i}{\bmat{I_{n_v}\\B_T^{-1}B_v\\T(b-a)B_T^{-1}B_v}}{\bmat{0_{n_r \times n_x}\\B_Q(s)\\T(b-a)B_Q(s)+Q(b-s)}}{U_2T(s-a)B_T^{-1}B_v}{U_1,R_{D,1},R_{D,2}}\right] \\
&=\fourpi{\bmat{0&D_{rb}}}{C_r}{\bmat{B_{xv}&B_{xb}}}{A_i}\fourpi{\bmat{I_{n_v}\\B_T^{-1}B_v\\T(b-a)B_T^{-1}B_v}}{\bmat{0_{n_r \times n_x}\\B_Q(s)\\T(b-a)B_Q(s)+Q(b-s)}}{U_2T(s-a)B_T^{-1}B_v}{U_1,R_{D,1},R_{D,2}},
\end{align*}}
to obtain
\begin{align*}
&\bmat{r(t)\\\dot{\hat{\mbf x}}(t)}=\fourpi{D_{rv}}{C_{rx}}{B_{xv}}{\hat A_i}\bmat{v(t)\\\ubar{\hat{\mbf x}}(t)}=\\
& \fourpi{\bmat{0&D_{rb}}}{C_r}{\bmat{B_{xv}&B_{xb}}}{A_i}\fourpi{\bmat{I_{n_v}\\B_T^{-1}B_v\\T(b-a)B_T^{-1}B_v}}{\bmat{0_{n_r \times n_x}\\B_Q(s)\\T(b-a)B_Q(s)+Q(b-s)}}{U_2T(s-a)B_T^{-1}B_v}{U_1,R_{D,1},R_{D,2}}\bmat{v(t)\\\ubar{\hat{\mbf x}}(t)}.
\end{align*}
We need to eliminate $\ubar{\hat{\mbf x}}(t)$ from the right hand side to get an expression solely in terms of $\hat{\mbf x}$. For this purpose, we use \Cref{thm:T_map}b, which gives us the relation $\ubar{\hat{\mbf x}}(t) = \mcl D \hat{\mbf x}(t)$. Defining now
\begin{align*}
(\mcl F\hat{\mbf x})(t):= \bmat{\hat{\mbf x}(t)\\ \partial_s S\hat{\mbf x}(t) \\ \vdots \\ \partial_s^{N}S^{N}\hat{\mbf x}(t)}, \qquad (\mcl C\hat{\mbf x})(t):= \bmat{S\hat{\mbf x}(t)\\ \partial_s S^2\hat{\mbf x}(t) \\ \vdots \\ \partial_s^{N-1}S^{N}\hat{\mbf x}(t)} \qquad (\mcl B \hat{\mbf x}) = \bmat{(\mcl C\hat{\mbf x})(t,a)\\(\mcl C\hat{\mbf x})(t,b)}.
\end{align*}
Using Corollaries~\ref{cor:x_c} and~\ref{cor:x_D}d, these definitions now imply
\begin{align*}
&\bmat{\bmat{v(t)\\(\mcl B \hat{\mbf x})(t)}\\(\mcl F\hat{\mbf x})(t,\cdot)}= \fourpi{\bmat{I_{n_v}\\B_T^{-1}B_v\\T(b-a)B_T^{-1}B_v}}{\bmat{0_{n_r \times n_x}\\B_Q(s)\\T(b-a)B_Q(s)+Q(b-s)}}{U_2T(s-a)B_T^{-1}B_v}{U_1,R_{D,1},R_{D,2}}\bmat{v(t)\\\ubar{\hat{\mbf x}}(t)}.
\end{align*}

Then, we can re-write the expressions for $r$ and $\dot{\hat{\mbf x}}$ as
\begin{align*}
&\bmat{r(t)\\\dot{\hat{\mbf x}}(t)}
= \fourpi{\bmat{0&D_{rb}}}{C_r}{\bmat{B_{xv}&B_{xb}}}{A_i}\bmat{\bmat{v(t)\\(\mcl B \hat{\mbf x})(t)}\\(\mcl F\hat{\mbf x})(t,\cdot)}\\
&= \sum\limits_{i=0}^N \bmat{
\int_{a}^{b} C_{r}(\cdot)	\\A_{0}(s) +\int\limits_{a}^{s}A_{1}(s,\cdot)+\int\limits_{s}^{b} A_{2}(s,\cdot)}(\mcl F\hat{\mbf x})(t,\cdot) +\bmat{ 0& D_{rb}\\B_{xv}(s) & B_{xb}(s) }\bmat{v(t)\\ (\mcl B \hat{\mbf x})(t)}.
\end{align*}

Thus we conclude that $\{\hat{\mbf x}, r\}$ satisfies the PDE \Cref{eq:general_pde_subsystem} with initial condition $\hat{\mbf x}^0$ and input $v$.
\end{proof}

\subsection{Bijective map between PIE and GPDE states: Proof of \Cref{cor:T_map_GPDE}}\label{app:T_map_GPDE}
We now construct the map between the domain of the GPDE and associated PIE representation and show this is a bijection.

\TmapGPDE*

\begin{proof} \textbf{Proof of Part 1.} Let $\bmat{x\\\hat{\mbf x}}\in \mcl X_{w,u}$ for some $w\in \R^{p}, u\in \R^q$. Clearly, by definition of $\mcl X_{w,u}$, $\hat{\mbf x}\in X_v$ with $v= C_vx+D_{vw}w+D_{vu}u$ for arbitrary matrices $C_v$, $D_{vw}$, and $D_{vu}$. Therefore, from \cref{thm:T_map}a, $\mcl D\hat{\mbf x}\in L_2^{n_{\hat{\mbf x}}}$ and hence $\{x,\mcl D\hat{\mbf x}\} \in \R L_2^{n_x,n_{\hat{\mbf x}}}$. Furthermore, for $\hat{\mcl T}$ and $\mcl T_v$ as defined in \cref{fig:Gb_definitions}, we have
\[
\hat{\mbf x} = \hat{\mcl T}\mcl D\hat{\mbf x}+ \mcl T_v v = \bmat{\mcl T_v C_v&\hat{\mcl T}}\bmat{x\\\mcl D \hat{\mbf x}}+ \mcl T_v \bmat{D_{vw}&D_{vu}}\bmat{w\\u}.
\]
Then, by concatenating $x$ and $\hat{\mbf x}$ and by using the definitions of $\mcl T$, $\mcl T_w$, $\mcl T_u$, we have
\begin{align*}
\bmat{x\\\hat{\mbf x}}&= \bmat{x \\\bmat{\mcl T_v C_v&\hat{\mcl T}}\bmat{x\\\mcl D \hat{\mbf x}}+ \mcl T_v \bmat{D_{vw}&D_{vu}}\bmat{w\\u}}\\
&= \bmat{x\\\bmat{\mcl T_v C_v&\hat{\mcl T}}\bmat{x\\\mcl D \hat{\mbf x}}} + \bmat{0\\\mcl T_v D_{vw}}w +\bmat{0\\\mcl T_vD_{vu}}u\\
&= \bmat{I&0\\\mcl T_v C_v&\hat{\mcl T}}\bmat{x\\\mcl D \hat{\mbf x}} + \bmat{0\\\mcl T_v D_{vw}}w +\bmat{0\\\mcl T_vD_{vu}}u = \mcl T \bmat{x\\\mcl D\hat{\mbf x}}+\mcl T_w w+\mcl T_u u.
\end{align*}
\end{proof}
\begin{proof} \textbf{Proof of Part 2.}
Let $w\in \R^p$, $u\in \R^q$ and $\ubar{\mbf x} \in \R L_2^{n_x, n_{\hat{\mbf x}}}$ be arbitrary.

Let $\bmat{x\\\ubar{\hat{\mbf x}}} := \ubar{\mbf x}$ where $x \in \R^{n_x}$ and $\ubar{\hat{\mbf x}}\in L_2^{n_{\hat{\mbf x}}}$.

By substituting the definitions of $\mcl T$, $\mcl T_w$ and $\mcl T_u$,
\begin{align*}
\mcl T \bmat{x\\\ubar{\hat{\mbf x}}} +\mcl T_w w+\mcl T_u u& = \bmat{I&0\\\mcl T_v C_v&\hat{\mcl T}}\bmat{x\\\ubar{\hat{\mbf x}}} + \bmat{0\\\mcl T_v D_{vw}}w +\bmat{0\\\mcl T_vD_{vu}}u =\underbrace{\bmat{x \\\hat{\mcl T}\ubar{\hat{\mbf x}}+ \mcl T_v \left(\bmat{C_v&D_{vw}&D_{vu}}\bmat{x\\w\\u}\right)}}_{\bmat{x\\\hat{\mbf x}}:=}.
\end{align*}
Clearly, from \cref{thm:T_map}b, defining $\hat{\mbf x}$ as $\hat{\mbf x} = \hat{\mcl T}\ubar{\hat{\mbf x}}+ \mcl T_v \left(\bmat{C_v&D_{vw}&D_{vu}}\bmat{x\\w\\u}\right)$ implies that $\hat{\mbf x}\in X_v$ with $v = \bmat{C_v&D_{vw}&D_{vu}}\bmat{x\\w\\u}$. Therefore, by definition of $\mcl X_{w,u}$, $\bmat{x\\\hat{\mbf x}}\in \mcl X_{w,u}$.

Our next step is to show that $\bmat{I&0\\0&\mcl D}\left(\mcl T  \ubar{\mbf{x}}+\mcl{T}_{w}w+\mcl T_u u\right)=\ubar{\mbf x} \in \R L_2^{n_x, n_{\hat{\mbf x}}}$. Earlier, we defined $\bmat{x\\\ubar{\hat{\mbf x}}}:= \ubar{\mbf x}$ and showed that if we define $\left(\mcl T  \ubar{\mbf{x}}+\mcl{T}_{w}w+\mcl T_u u\right) = \bmat{x\\\hat{\mbf x}}$ for some $\hat{\mbf x}\in X_v$ with $v = C_v x+D_{vw}w +D_{vu}u$. Thus, from \cref{thm:T_map}b, we have
\[
\mcl D\hat{\mbf x} = \mcl D\left(\hat{\mcl T}\ubar{\hat{\mbf x}}+\mcl T_v v\right) = \ubar{\hat{\mbf x}}.
\]
Therefore,
\[
\bmat{I&0\\0&\mcl D}\left(\mcl T  \ubar{\mbf{x}}+\mcl{T}_{w}w+\mcl T_u u\right) = \bmat{I&0\\0&\mcl D}\bmat{x\\\hat{\mbf x}} = \bmat{x\\\mcl D\hat{\mbf x}} = \bmat{x\\\ubar{\hat{\mbf x}}} = \ubar{\mbf x}.
\]

\end{proof}

\subsection{Equivalence of PIE and GPDE: Proof of \Cref{thm:equivalence_2_reverse}}\label{app:equivalence_2}
The equivalence of solutions between a GPDE model and associated PIE is a straightforward extension of~\Cref{thm:equivalence}). This proof is split into two parts.

\corEquivalenceReverse*
\begin{proof}\textbf{Proof of Part 1}
Suppose $\{x,$ $\hat{\mbf x},$ $z,$ $y,$ $v,$ $r\}$ satisfies the GPDE defined by $\{\mbf G_{\mathrm{o}},$ $\mbf G_{\mathrm b},$ $\mbf G_{\mathrm{p}}\}$ with initial condition $\{x^0,$ $\hat{\mbf x}^0\}$ and input $\{w,u\}$. Then, we have: a) $x\in W_{1e}^{n_x}[\R_+]$, $z\in L_{2e}^{n_z}[\R_+]$, $y\in L_{2e}^{n_y}[\R_+]$, $v\in L_{2e}^{n_v}[\R_+]$, $r\in L_{2e}^{n_r}[\R_+]$; b) $\hat{\mbf x}(t) \in X_{v(t)}$ for all $t\ge 0$; c) $x$ is differentiable  on $\R_+$, $\hat{\mbf x}$ is Frech\'et differentiable with respect to the $L_2$-norm  on $\R_+$; d) \Cref{eq:odepde_general,eq:general_pde_subsystem} are satisfied for all $t\ge 0$; and e) $x(0)=x^0$, $\hat{\mbf x}(0) = \hat{\mbf x}^0$ and $\hat{\mbf x}^0 \in X_{v(0)}$.

Now, from above points, since $\hat{\mbf x}^0 \in X_{v(0)}$ and $\hat{\mbf x}(0) = \hat{\mbf x}^0$, $r\in L_{2e}^{n_r}[\R_+]$, $\hat{\mbf x}(t) \in X_{v(t)}$ for all $t\ge 0$, $\hat{\mbf x}$ is Frech\'et differentiable with respect to the $L_2$-norm  on $\R_+$, and \Cref{eq:general_pde_subsystem} is satisfied for all $t\ge 0$, we have that $\{\hat{\mbf x},$ $r\}$ satisfies the PDE defined by $n\in \N^{N+1}$ and $\{\mbf G_{\mathrm b}, \mbf G_{\mathrm{p}}\}$ with initial conditions $\hat{\mbf x}^0$ and input $v$. Furthermore, since
\[
v(t) = C_vx(t)+D_{vw}w(t)+D_{vu}u(t),
\]
we have that $v\in L_{2e}^{n_v}[\R_+]$ with $B_v v\in W_{1e}^{2n_S}[\R_+]$. Thus, by \Cref{thm:equivalence}, $\{\mcl D\hat{\mbf x}, ~r\}$ is a solution to the PIE defined
\[
\mbf G_{\mathrm{PIE}_s} = \left\{\hat{\mcl T},\mcl T_v, \emptyset, \hat{\mcl A},\mcl B_v,\emptyset,\mcl C_r,\emptyset,\mcl D_{rv},\emptyset,\emptyset,\emptyset\right\}
\]
with initial condition  $\mcl D\hat{\mbf x}^0 \in L_2^{n_{\hat{\mbf x}}}$. Therefore, if we define $\ubar{\hat{\mbf x}}(t)=\mcl D \hat{\mbf x}(t)$ and $\ubar{\hat{\mbf x}}^0=\mcl D \hat{\mbf x}^0$, we have that: f) $v\in L_{2e}^{n_v}[\R_+]$ and $(\mcl{T}_{v} v)(\cdot,s)\in W_{1e}^{n_{\hat{\mbf x}}}[\R_+]$ for all $s \in [a,b]$; g) $\ubar{\hat{\mbf x}}:\R_+\to \R L_2^{0,n_{\hat{\mbf x}}}[a,b]$ and $r\in L_{2e}^{n_r}[\R_+]$; h) $\ubar{\hat{\mbf x}}^0\in \R L_2^{0,n_{\hat{\mbf x}}}[a,b]$ and $\ubar{\hat{\mbf x}}(0) = \ubar{\hat{\mbf x}}^0$; i) $\hat{\mcl T}\ubar{\hat{\mbf x}}$ is Frech\'et differentiable on $\R_+$; and j) Equation~\eqref{eq:PIE_full} (defined by $\mbf G_{\mathrm{PIE}_s}$) is satisfied for all $t\in\R_+$, i.e.
\begin{align*}
&	\bmat{r(t)\\\hat{\mcl{T}}\ubar{\dot{\hat{\mbf{x}}}}(t)+\mcl{T}_{v}\dot{v}(t)} = \bmat{\mcl D_{rv}&\mcl C_{r}\\\mcl B_{v}&\hat{\mcl A}}\bmat{v(t)\\\ubar{\hat{\mbf x}}(t)}.
\end{align*}

Now, let $\ubar{\mbf x}^0=\bmat{x^0\\\mcl D\hat{\mbf x}^0}=\bmat{x^0\\\ubar{\hat{\mbf x}}^0}$ and $\ubar{\mbf x}(t) = \bmat{x(t)\\\mcl D\hat{\mbf x}(t)}= \bmat{x(t)\\\ubar{\hat{\mbf x}}(t)}$ for all $t\ge 0$. Our goal is to show that $\{\ubar{\mbf x}, z, y\}$ satisfies the PIE defined by $\mbf G_{\mathrm{PIE}}$ with initial condition $\ubar{\mbf x}^0$ and input $\{w, u\}$, which means we need to show that: 1) $(\mcl T_w w)(\cdot,s), (\mcl T_u u)(\cdot,s)\in W_{1e}^{n_x+n_{\hat{\mbf x}}}$ for all $s\in [a,b]$; 2) $\ubar{\mbf x}(t)\in \R L_2^{n_x,n_{\hat{\mbf x}}}[a,b]$ for all $t\ge 0$; 3) $\ubar{\mbf x}(0) = \ubar{\mbf x}^0$ and $\ubar{\mbf x}^0\in \R L_2^{n_x,n_{\hat{\mbf x}}}$; 4) $\mcl T\ubar{\mbf x}$ is Frech\'et differentiable on $\R_+$; and 5) \Cref{eq:PIE_full} (defined by $\mbf G_{\mathrm{PIE}}$) is satisfied for all $t\in\R_+$.

\noindent For 1), $B_vD_{vw} w\in W_{1e}^{2n_S}[\R_+]$ and hence by the definition of $\mcl T_w$, we have
\[
(\mcl{T}_{w} w(\cdot))(s) = \bmat{0\\T_1(s-a)}B_T^{-1}B_v D_{vw} w(\cdot)\in W_{1e}^{n_x+n_{\hat{\mbf x}}}[\R_+].
\]
Likewise, $B_vD_{vu} u\in W_{1e}^{2n_S}[\R_+]$ implies
\[
(\mcl{T}_{u} u(\cdot))(s) = \bmat{0\\T_1(s-a)}B_T^{-1}B_v D_{vu} u(\cdot)\in W_{1e}^{n_x+n_{\hat{\mbf x}}}[\R_+].
\]

\noindent For 2), since $\ubar{\hat{\mbf x}}(t) \in L_2^{0,n_{\hat{\mbf x}}}[a,b]$ and $x(t)\in \R^{n_x}$, we have $\ubar{\mbf x}(t) =\bmat{x(t)\\\ubar{\hat{\mbf x}}(t)} \in \R L_2^{n_x, n_{\hat{\mbf x}}}$ for all $t\ge 0$.

\noindent For 3), since $\ubar{\hat{\mbf x}}^0 \in L_2^{0,n_{\hat{\mbf x}}}[a,b]$ and $x^0\in \R^{n_x}$, we have $\ubar{\mbf x}^0 =\bmat{x^0\\\ubar{\hat{\mbf x}}^0} \in \R L_2^{n_x, n_{\hat{\mbf x}}}$. Furthermore, $
\ubar{\mbf x}(0) =\bmat{x(0)\\\ubar{\hat{\mbf x}}(0)}=\bmat{x^0\\\ubar{\hat{\mbf x}}^0} = \ubar{\mbf x}^0$.

\noindent For 4), by definitions of $\mcl T$ and $\ubar{\mbf x}$, there exists a $k>0$ such that
\begin{align*}
\norm{\mcl T\ubar{\mbf x}(t)}_{L_2} &= \norm{\bmat{I&0\\G_vC_v&\hat{\mcl T}}\bmat{x(t)\\\ubar{\hat{\mbf x}}(t)}}_{L_2}=\norm{\bmat{x(t)\\G_vC_vx(t)+\hat{\mcl T}\ubar{\hat{\mbf x}}(t)}}_{L_2}\\
&=\norm{x(t)}_{L_2}+\norm{G_vC_vx(t)+\hat{\mcl T}\ubar{\hat{\mbf x}}(t)}_{L_2} \le k\norm{x(t)}_{\R}+\norm{\hat{\mcl T}\ubar{\hat{\mbf x}}(t)}_{L_2}.
\end{align*}
Since $\hat{\mcl T}\ubar{\hat{\mbf x}}(t)$ is Frech\'et differentiable and $x \in W_{1e}^{n_x}$ is differentiable, we have that $\mcl T\ubar{\mbf x}(t)$ is Frech\'et differentiable.

Finally, for 5), we need to show that
\[
\bmat{\mcl{T}\ubar{\dot{\mbf x}}(t)\\z(t)\\y(t)}=\bmat{\mcl{A}&\mcl{B}_1&\mcl{B}_2\\\mcl{C}_1&\mcl{D}_{11}&\mcl{D}_{12}\\\mcl{C}_2&\mcl{D}_{21}&\mcl{D}_{22}}\bmat{\ubar{\mbf x}(t)\\w(t)\\u(t)}-\bmat{\mcl{T}_{w}\dot{w}(t)+\mcl{T}_{u}\dot{u}(t)\\0\\0}
\]
is satisfied for all $t\ge 0$.

Since $x,z,y,v$ satisfy the GPDE, we have
\begin{equation}
\bmat{\dot{x}(t)\\\hline z(t)\\y(t)\\v(t)} = \bmat{A&\vlines&B_{xw}&B_{xu} & B_{xr}\\\hline C_z&\vlines&D_{zw}&D_{zu}&D_{zr}\\C_y&\vlines&D_{yw}&D_{yu}&D_{yr}\\C_v&\vlines&D_{vw}&D_{vu}&0}\bmat{x(t)\\\hline w(t)\\u(t)\\r(t) }.\label{eqn:subsys_ODE}
\end{equation}
Furthermore, as stated above,
\begin{equation}
\bmat{\hat{\mcl{T}}\ubar{\dot{\hat{\mbf{x}}}}(t)+\mcl{T}_{v}\dot{v}(t)\\ r(t)} = \bmat{\hat{\mcl A}&\mcl B_{v}\\\mcl C_{r}&\mcl D_{rv}}\bmat{\ubar{\hat{\mbf x}}(t)\\v(t)}.\label{eqn:subsys_temp}
\end{equation}
These two identities are all that are required to conclude the proof. Specifically, extracting expression for $v,r$ and $\mcl T_v\dot v$, we obtain
\begin{align*}
v(t)&=\bmat{C_v&D_{vw}&D_{vu}}\bmat{x(t)\\w(t)\\u(t)},\\
r(t) &= \mcl C_r \ubar{\hat{\mbf x}}(t)  + D_{rv}v(t)= \bmat{\mcl D_{rv}C_v&\mcl C_r&\mcl D_{rv}D_{vw}&\mcl D_{rv}D_{vu}}\bmat{x(t)\\ \ubar{\hat{\mbf x}}(t)\\w(t)\\u(t)},\\
(\mcl T_v\dot v(t))(s)&=G_v(s)\bmat{C_v&D_{vw}&D_{vu}}\bmat{\dot x(t)\\\dot w(t)\\\dot u(t)}\\&=G_v(s)C_v\dot x(t)+G_v(s) D_{vw}\dot w(t)+G_v(s) D_{vu} \dot u(t).
\end{align*}
Substituting these expressions back into Eq.~\eqref{eqn:subsys_temp} yields
\begin{align*}
&\hat{\mcl{T}}\ubar{\dot{\hat{\mbf{x}}}}(t)+G_v(s)C_v\dot x(t)+G_v(s) D_{vw}\dot w(t)+G_v(s) D_{vu} \dot u(t) = \hat{\mcl A}\ubar{\hat{\mbf x}}(t)+\mcl B_{v}\bmat{C_v&D_{vw}&D_{vu}}\bmat{x(t)\\w(t)\\u(t)}
\end{align*}
or
\begin{align*}
&\bmat{G_v(s)C_v&\hat{\mcl{T}}}\bmat{\dot x(t)\\\ubar{\dot{\hat{\mbf{x}}}}(t)}+G_v(s) D_{vw}\dot w(t)+G_v(s) D_{vu} \dot u(t) = \bmat{\mcl B_{v}C_v&\hat{\mcl A}&\mcl B_{v}D_{vw}&\mcl B_{v}D_{vu}}\bmat{x(t)\\\ubar{\hat{\mbf x}}(t)\\w(t)\\u(t)}.
\end{align*}
Appending the above equation to the system of equations in Eq.~\eqref{eqn:subsys_ODE} and omitting the equation for $v$ yields
\begin{flalign*}
&\bmat{\bmat{I&0\\G_v(s)C_v&\hat{\mcl{T}}}\bmat{\dot{x}(t)\\ \ubar{\dot{\hat{\mbf{x}}}}}+\mcl T_w \dot w(t)+\mcl T_u \dot u(t)\\\hline z(t)\\y(t)}=
\bmat{A&0&\vlines&B_{xw}&B_{xu} \\
\mcl B_{v}C_v&\hat{\mcl A}& \vlines &\mcl B_{v}D_{vw} & \mcl B_{v}D_{vu}\\
\hline C_z&0&\vlines&D_{zw}&D_{zu}\\
C_y&0&\vlines&D_{yw}&D_{yu}}
\bmat{x(t)\\ \ubar{\hat{\mbf x}}(t)\\\hline w(t)\\u(t) }+\bmat{\bmat{B_{xr}\\0}\\D_{zr}\\D_{yr}}r(t)\notag\\
&\qquad=
\bmat{A&0&\vlines&B_{xw}&B_{xu} \\
\mcl B_{v}C_v&\hat{\mcl A}& \vlines &\mcl B_{v}D_{vw} & \mcl B_{v}D_{vu}\\
\hline C_z&0&\vlines&D_{zw}&D_{zu}\\
C_y&0&\vlines&D_{yw}&D_{yu}}
\bmat{x(t)\\ \ubar{\hat{\mbf x}}(t)\\\hline w(t)\\u(t) }+\bmat{\bmat{B_{xr}\\0}\\D_{zr}\\D_{yr}}\bmat{\mcl D_{rv}C_v&\mcl C_r&\mcl D_{rv}D_{vw}&\mcl D_{rv}D_{vu}}\bmat{x(t)\\ \ubar{\hat{\mbf x}}(t)\\w(t)\\u(t)}\notag\\
&\qquad=
\bmat{\mcl A&\vlines&\mcl B_{1}&\mcl B_{2} \\
\hline \mcl C_1&\vlines&\mcl D_{11}&\mcl D_{12}\\
\mcl C_2&\vlines&\mcl D_{21}&\mcl D_{22}}
\bmat{\ubar{\mbf x}(t)\\\hline w(t)\\u(t) }.
\end{flalign*}
We conclude that
\[
\bmat{\mcl{T}\ubar{\dot{\mbf{x}}}(t)\\z(t)\\y(t)}=\bmat{\mcl{A}&\mcl{B}_1&\mcl{B}_2\\\mcl{C}_1&\mcl{D}_{11}&\mcl{D}_{12}\\\mcl{C}_2&\mcl{D}_{21}&\mcl{D}_{22}}\bmat{\ubar{\mbf x}(t)\\w(t)\\u(t)}-\bmat{\mcl{T}_{w}\dot{w}(t)+\mcl{T}_{u}\dot{u}(t)\\0\\0}
\]
which implies that $\{\bmat{x\\\mcl D\hat{\mbf x}},\, z,\,y\}$ satisfies the PIE defined by
$\mbf G_{\mathrm{PIE}}$
with initial condition $\bmat{x^0\\\mcl D\hat{\mbf x}} \in \R L_2^{n_x, n_{\hat{\mbf x}}}$ and input $\{w,u\}$.
\end{proof}

\begin{proof}\textbf{Proof of Part 2} In this proof, we will use definitions in~\Cref{fig:PIE_subsystem_equation} using the parameters contained in $\{\mbf G_{\mathrm{o}},$ $\mbf G_{\mathrm b},$ $\mbf G_{\mathrm{p}}\}$.

Now, suppose $\{\ubar{\mbf x},~z,~y\}$ satisfies the PIE defined by $\mbf G_{\mathrm{PIE}}$ with initial condition $\ubar{\mbf x}^0$ and input $\{w,u\}$. Then, by definition of solution of a PIE: a) $z\in L_{2e}^{n_z}[\R_+]$, $y\in L_{2e}^{n_y}[\R_+]$; b) $\ubar{\mbf x}(t)\in \R L_2^{n_x,n_{\hat{\mbf x}}}[a,b]$ for all $t\ge 0$; c) $\mcl T\ubar{\mbf x}$ is Frech\'et differentiable on $\R_+$; d) $\ubar{\mbf x}(0) = \ubar{\mbf x}^0$; and e) The equation
\begin{align}\label{eq:PIE_full_temp}
\bmat{\mcl{T}\ubar{\dot{\mbf{x}}}(t)+\mcl{T}_{w}\dot{w}(t)+\mcl{T}_{u}\dot{u}(t)\\z(t)\\y(t)}=\bmat{\mcl{A}&\mcl{B}_1&\mcl{B}_2\\\mcl{C}_1&\mcl{D}_{11}&\mcl{D}_{12}\\\mcl{C}_2&\mcl{D}_{21}&\mcl{D}_{22}}\bmat{\ubar{\mbf x}(t)\\w(t)\\u(t)}
\end{align} is satisfied for all $t\in\R_+$.

For $\ubar{\mbf x}(t)\in \R L_2^{n_x,n_{\hat{\mbf x}}}$ we define $\hat x(t)\in \R^{n_x}$ and $\ubar{\hat{\mbf x}}(t)\in L_2^{n_{\hat{\mbf x}}}$ by $\bmat{\hat{x}(t)\\\ubar{\hat{\mbf x}}(t)}:=\ubar{\mbf x}(t)$. Similarly, we define the elements $\bmat{\hat x^0\\\ubar{\hat{\mbf x}}^0}:=\ubar{\mbf x}^0$. Now, by the definitions of $\mcl T$, $\mcl T_w$ and $\mcl T_u$, we have
\[
\bmat{x(t)\\\hat{\mbf x}(t)} = \bmat{I&0\\G_v C_v&\hat{\mcl T}}\bmat{\hat x(t)\\\ubar{\hat{\mbf x}}(t)}+\bmat{0\\G_v D_{vw}}w(t)+\bmat{0\\G_v D_{vu}}u(t)
\]
and hence $x(t)=\hat x(t)$. Similarly, $x^0=\hat x^0$. Hence we have $\ubar{\mbf x}(t)=\bmat{x(t)\\\ubar{\hat{\mbf x}}(t)}$ and $\ubar{\mbf x}^0=\bmat{x^0\\\ubar{\hat{\mbf x}}^0}$.

Now, using the definitions of $r$ and $v$ and examining the right hand side of Eq.~\eqref{eq:PIE_full_temp}, we have
{\small\begin{flalign*}
\bmat{\mcl A&\vlines&\mcl B_{1}&\mcl B_{2} \\
\hline \mcl C_1&\vlines&\mcl D_{11}&\mcl D_{12}\\
\mcl C_2&\vlines&\mcl D_{21}&\mcl D_{22}}
\bmat{\ubar{\mbf x}(t)\\\hline w(t)\\u(t) }
&=
\bmat{A&0&\vlines&B_{xw}&B_{xu} \\
\hat B_{xv}C_v&\hat{\mcl A}& \vlines &\hat B_{xv}D_{vw} & \hat B_{xv}D_{vu}\\
\hline C_z&0&\vlines&D_{zw}&D_{zu}\\
C_y&0&\vlines&D_{yw}&D_{yu}}
\bmat{x(t)\\ \ubar{\hat{\mbf{x}}}(t)\\\hline w(t)\\u(t) }+\bmat{\bmat{B_{xr}\\0}\\D_{zr}\\D_{yr}}\bmat{\mcl D_{rv}C_v&\mcl C_r&\mcl D_{rv}D_{vw}&\mcl D_{rv}D_{vu}}\bmat{x(t)\\ \ubar{\hat{\mbf x}}(t)\\w(t)\\u(t)}\notag\\
&=
\bmat{A&0&\vlines&B_{xw}&B_{xu} \\
\mcl B_vC_v&\hat{\mcl A}& \vlines &\mcl B_vD_{vw} & \mcl B_vD_{vu}\\
\hline C_z&0&\vlines&D_{zw}&D_{zu}\\
C_y&0&\vlines&D_{yw}&D_{yu}}
\bmat{x(t)\\ \ubar{\hat{\mbf{x}}}(t)\\\hline w(t)\\u(t) }+\bmat{\bmat{B_{xr}\\0}\\D_{zr}\\D_{yr}}r(t)\notag\\
&=
\bmat{A&0&\vlines&B_{xw}&B_{xu} & B_{xr}\\
\mcl B_{v}C_v&\hat{\mcl A}& \vlines &\mcl B_vD_{vw} & \mcl B_vD_{vu}&0\\
\hline C_z&0&\vlines&D_{zw}&D_{zu}&D_{zr}\\
C_y&0&\vlines&D_{yw}&D_{yu}&D_{yr}}
\bmat{x(t)\\ \ubar{\hat{\mbf x}}(t)\\\hline w(t)\\u(t)\\r(t)}.&
\end{flalign*}}
Likewise, if we substitute the definitions of the PI operators $\mcl T$, $\mcl T_w$, and $\mcl T_u$ in the left hand side of Eq.~\eqref{eq:PIE_full_temp}, we get
\begin{align*}
\bmat{\mcl{T}\ubar{\dot{\mbf{x}}}(t)+\mcl{T}_{w}\dot{w}(t)+\mcl{T}_{u}\dot{u}(t)\\z(t)\\y(t)}&=\bmat{ \bmat{I&0\\G_v(s) C_v&\hat{\mcl T}}\bmat{\dot x(t)\\\ubar{\dot{\hat{\mbf{x}}}}(t)}+\bmat{0\\G_v(s) D_{vw}}\dot{w}(t)+\bmat{0\\G_v(s) D_{vu}}\dot{u}(t)\\z(t)\\y(t)}\\
&=\bmat{\bmat{\dot x(t)\\\hat{\mcl T} \ubar{\dot{\hat{\mbf{x}}}}(t)}+\bmat{0\\\mcl T_v\dot{v}(t)}\\z(t)\\y(t)}.
\end{align*}
Adding the definition of $v$, we conclude that
\[
\bmat{\bmat{\dot x(t)\\\hat{\mcl T} \ubar{\dot{\hat{\mbf{x}}}}(t)}+\bmat{0\\\mcl T_v\dot{v}(t)}\\z(t)\\y(t)\\v(t)}
=\bmat{A&0&\vlines&B_{xw}&B_{xu} & B_{xr}\\
\mcl B_vC_v&\hat{\mcl A}& \vlines &\mcl B_vD_{vw} & \mcl B_vD_{vu}&0\\
\hline C_z&0&\vlines&D_{zw}&D_{zu}&D_{zr}\\
C_y&0&\vlines&D_{yw}&D_{yu}&D_{yr}\\
C_{v} &0&\vlines&D_{vw}&D_{vu}&0}
\bmat{x(t)\\ \ubar{\hat{\mbf x}}(t)\\\hline w(t)\\u(t)\\r(t) }.
\]
Therefore,
\begin{equation*}
\bmat{\dot{x}(t)\\\hline z(t)\\y(t)\\v(t)} = \bmat{A&\vlines&B_{xw}&B_{xu} & B_{xr}\\\hline C_z&\vlines&D_{zw}&D_{zu}&D_{zr}\\C_y&\vlines&D_{yw}&D_{yu}&D_{yr}\\C_v&\vlines&D_{vw}&D_{vu}&0}\bmat{x(t)\\\hline w(t)\\u(t)\\r(t) }
\end{equation*}
and
\begin{equation*}
\bmat{\hat{\mcl{T}}\ubar{\dot{\hat{\mbf{x}}}}(t)+\mcl{T}_{v}\dot{v}(t)\\ r(t)} = \bmat{\hat{\mcl A}&\mcl B_{v}\\\mcl C_{r}&\mcl D_{rv}}\bmat{\ubar{\hat{\mbf x}}(t)\\v(t)}.
\end{equation*}

Thus, we conclude: f) Since $\ubar{\mbf x}$ is Frech\'et differentiable, $x$ and $\hat{\mcl T}\ubar{\hat{\mbf x}}$ are Frech\'et differentiable; g) Since $x$ is Frech\'et differentiable, $ x\in W_{1e}^{n_x}$ and $w\in L_{2e}^{n_w}$, $u\in L_{2e}^{n_u}$ from the theorem statement, thus $v\in L_{2e}^{n_v}$; h) Since $\hat{\mcl T}\ubar{\hat{\mbf x}}$ is Frech\'et differentiable and $v\in L_{2e}^{n_v}$, we have $r\in L_{2e}^{n_r}$; i) Since $\ubar{\mbf x}^0\in \R L_2^{n_x,n_{\hat{\mbf x}}}$, we have $ x^0\in \R^{n_x}$, $\ubar{\hat{\mbf x}}^0\in L_2^{n_{\hat{\mbf x}}}$; and j) For all $t\ge 0$
\[
\bmat{r(t)\\\hat{\mcl T}\ubar{\dot{\hat{\mbf{x}}}}(t) + \mcl T_v \dot v(t)} = \bmat{\mcl D_{rv}&\mcl C_v\\\mcl B_v&\hat{\mcl A}}\bmat{v(t)\\\ubar{\hat{\mbf x}}(t)}.
\]

Thus, $\{\ubar{\hat{\mbf x}},r\}$ (as defined above), satisfies the PIE defined by
\[
\mbf G_{\mathrm{PIE}_s}= \{\hat{\mcl{T}},\mcl{T}_{v},\emptyset,\hat{\mcl A},\mcl{B}_{v},\emptyset,\mcl{C}_v, \emptyset,\mcl D_{rv},\emptyset,\emptyset,\emptyset\}
\]
for initial condition $\ubar{\hat{\mbf x}}^0$ and input $\{v,\emptyset\}$.
Thus, from \cref{thm:equivalence}, $\{\hat{\mcl T}\ubar{\hat{\mbf x}}+\mcl T_v v,~r\}$ satisfies the PDE defined by $n$ and $\{\mbf G_{\mathrm b}, \mbf G_{\mathrm{p}}\}$ with initial condition $\hat{\mcl T}\ubar{\hat{\mbf x}}^0 + \mcl T_v v(0)$ and input $v$. Since
\[
\bmat{x(t)\\\hat{\mbf x}(t)} = \mcl{T}\ubar{\mbf x}(t)+\mcl{T}_{w}w(t)+\mcl T_uu(t)=\bmat{ x(t)\\\hat{\mcl T} \ubar{\hat{\mbf{x}}}(t)}+\bmat{0\\\mcl T_v{v}(t)}
\]
and since similarly $\hat{\mbf{x}}^0=\hat{\mcl T} \ubar{\hat{\mbf{x}}}^0+\mcl T_v{v}(0)$, by the definition of solution of a PDE in~\ref{defn:PDE}, we have: k) $\hat{\mbf x}(t) \in X_{v(t)}$ for all $t\ge 0$; l) $\hat{\mbf x}$ is Frech\'et differentiable with respect to the $L_2$-norm on $\R_+$; m) \Cref{eq:general_pde_subsystem} is satisfied for all $t\ge 0$; and n) $\hat{\mbf x}(0) = \hat{\mbf x}^0$.

Reviewing all the above steps, we conclude that: 1) $z\in L_{2e}^{n_z}[\R_+]$ and $y\in L_{2e}^{n_y}[\R_+]$ by definition of solution of the PIE defined by $\mbf G_{\mathrm{PIE}}$; 2) $v\in L_{2e}^{n_v}[\R_+]$ and $r\in L_{2e}^{n_r}[\R_+]$ since $r,v$ satisfy the PDE; 3) $\hat{\mbf x}(t) \in X_{v(t)}$ for all $t\ge 0$ since $\hat{\mbf x}$ satisfies the PDE; 4) $x \in W_{1e}^{n_{x}}$ since $\hat{\mbf x}$ is Frech\'et differentiable; 5) $\hat{\mbf x}$ is Frech\'et differentiable with respect to the $L_2$-norm since $\hat{\mbf x}$ satisfies the PDE; $\bmat{x(0)\\\hat{\mbf x(0)}} = \mcl{T}\ubar{\mbf x}(0)+\mcl{T}_{w}w(0)+\mcl T_uu(0) = \mcl{T}\ubar{\mbf x}^0+\mcl{T}_{w}w(0)+\mcl T_uu(0)$ by definition of $x$ and $\hat{\mbf x}$; and 6) \Cref{eq:odepde_general,eq:general_pde_subsystem} are satisfied for all $t\ge 0$ as shown above.

We conclude that $\{x,$ $\hat{\mbf x},$ $z,$ $y,$ $v,$ $r\}$ satisfies the GPDE defined by $n$ and $\{\mbf G_{\mathrm{o}},$ $\mbf G_{\mathrm b},$ $\mbf G_{\mathrm{p}}\}$ with initial condition $\bmat{x^0\\\hat{\mbf x}^0}$ and input $\{w,u\}$.
\end{proof}
\subsection{Proof of \Cref{thm:unitary_T}}\label{app:unitary_T}
Having proven the equivalence between solutions of GPDE model and associated PIE, we now prove that these models have the same internal stability properties. Specifically, when $u=w=0$, the solution to associated PIE is stable if and only if the solution to GPDE model is internally stable. We do this in three parts. First, we show that the map $\ubar{\mbf x} \rightarrow \mcl T \ubar{\mbf x}+\mcl T_w w+\mcl T_u u$ is an isometric map between inner product spaces $L_2$ and $X^n$. Next, we show that the $W^n$ and $X^n$ (defined in \Cref{eq:ip_xv}) norms are equivalent. Finally, we show equivalence of internal stability in the respective norms.

\thmUnitary*

\begin{proof}
Let $\ubar{\hat{\mbf{x}}}, \ubar{\hat{\mbf y}}\in L_2^{n_{\hat{\mbf x}}}$  and $v_1, v_2\in \R^p$. Then, from \Cref{thm:T_map}, we have
\begin{align*}
\hat{\mcl{T}}\ubar{\hat{\mbf x}} +\mcl{T}_{v}v_1 \in X_{v_1},\quad \hat{\mcl{T}}\ubar{\hat{\mbf y}} +\mcl{T}_{v}v_2 \in X_{v_2}.
\end{align*}
Therefore, by definition \Cref{eq:ip_xv} and the result in \Cref{thm:T_map}b,
\begin{align*}
\ip{\left(\hat{\mcl{T}}\ubar{\hat{\mbf x}}+\mcl{T}_v v_1 \right)}{\left(\hat{\mcl{T}}\ubar{\hat{\mbf y}}+\mcl{T}_v v_2\right)}_{X^n} &= \ip{\mcl D\left(\hat{\mcl{T}}\ubar{\hat{\mbf x}}+\mcl{T}_v v_1 \right)}{\mcl D\left(\hat{\mcl{T}}\ubar{\hat{\mbf y}}+\mcl{T}_v v_2 \right)}_{L_2^{n_{\hat{\mbf x}}}} = \ip{\ubar{\hat{\mbf x}}}{\ubar{\hat{\mbf y}}}_{L_2^{n_{\hat{\mbf x}}}}.
\end{align*}

For b), let $\ubar{\mbf x}, \ubar{\mbf y} \in \R L_2^{n_x,n_{\hat{\mbf x}}}$ and $w_1, w_2\in \R^p$, $u_1, u_2\in \R^q$. Then, from \Cref{cor:T_map_GPDE}, we have
\begin{align*}
\mcl T\ubar{\mbf x} +\mcl{T}_{w}w_1+\mcl T_u u_1 \in \mcl X_{w_1,u_1},\quad \mcl T \ubar{\mbf y} +\mcl{T}_{w}w_2 +\mcl T_u u_2\in \mcl X_{w_2, u_2}.
\end{align*}
Since $\R^{n_x}\times X^n$ inner product is just sum of $\R$ and $X^n$ inner products, using definitions of $\mcl T$, $\mcl T_w$, and $\mcl T_u$ and the result in \Cref{cor:T_map_GPDE}b, we have
\begin{align*}
&\ip{\left(\mcl T\ubar{\mbf x} +\mcl{T}_{w}w_1+\mcl T_u u_1\right)}{\left(\mcl T\ubar{\mbf y} +\mcl{T}_{w}w_2+\mcl T_u u_2\right)}_{\R^{n_x} \times X^n} \\
&= \ip{\bmat{I&0\\0&\mcl D}\left(\mcl T\ubar{\mbf x} +\mcl{T}_{w}w_1+\mcl T_u u_1\right)}{\bmat{I&0\\0&\mcl D}\left(\mcl T\ubar{\mbf y} +\mcl{T}_{w}w_2+\mcl T_u u_2\right)}_{\R L_2^{n_x, n_{\hat{\mbf x}}}} =\ip{\ubar{\mbf x}}{\ubar{\mbf y}}_{\R L_2^{n_x, n_{\hat{\mbf x}}}}.
\end{align*}
\end{proof}
Using this result, we conclude that when $v=0$, the PI map ($\hat{\mcl T}$) is unitary. Since $\hat{\mcl T}$ is a unitary map from $L_2$ to $X_v$, the space $X_v$ is complete under the $X$-norm because $L_2$ is complete.
\subsection{Proof of \Cref{lem:norm_equivalence}}\label{app:norm_equivalence}
Next, we prove that the $\R X$ norm is equivalent to the $W^n$-norm on the subspace $\R\times X$.

\lemNormEquivalence*
\begin{proof}
Suppose $\mcl X_{0,0}$ is as defined in \Cref{eq:GPDE_domain}. Then, for any $\bmat{x\\\hat{\mbf x}}\in \mcl X_{0,0}$, we have, $x\in \R^{n_x}$ and $\hat{\mbf x}\in X_{C_v x}$ for some matrix $C_v$ and hence, from \Cref{thm:T_map}, we have $\hat{\mbf x} = \hat{\mcl T}\mcl D\hat{\mbf x}+\mcl T_v C_v x$ where $\hat{\mcl T}$ and $\mcl T_v$ are as defined in \Cref{fig:Gb_definitions}.

Let the space $\mcl X_{0,0}$ be equipped with two different inner products $\R^{n_x} \times X^n$ and $\R^{n_x}\times W^n$. Then
\begin{flalign*}
\norm{\bmat{x\\\hat{\mbf x}}}_{\R^{n_x}\times W^n}^2 &= \norm{x}_{\R}^2+\sum_{i=0}^{N}\sum_{j=0}^{i} \norm{\partial_s^j\hat{\mbf x}_i}_{L_2}^2 = \norm{x}_{\R}^2+ \sum_{i=0}^N\norm{\partial_s^i\hat{\mbf x}_i}_{L_2}^2 +\sum_{i=0}^{N}\sum_{j=0}^{i-1} \norm{\partial_s^j\hat{\mbf x}_i}_{L_2}^2  \\
&\ge \norm{x}_{\R}^2+ \sum_{i=0}^N\norm{\partial_s^i\hat{\mbf x}_i}_{L_2}^2  = \norm{x}_{\R}^2+\norm{\hat{\mbf x}}^2_{X^n}\ge \norm{\bmat{x\\\hat{\mbf x}}}_{\R^{n_x} \times X^n}^2. &
\end{flalign*}
For the reverse inequality we try to find an upper bound on $\norm{\cdot}_{\R^{n_x}\times W^n}$ as follows.
\begin{flalign*}
\norm{\bmat{x\\\hat{\mbf x}}}_{\R^{n_x}\times W^n}^2 &= \norm{x}_{\R}^2+\sum_{i=0}^{N}\sum_{j=0}^{i} \norm{\partial_s^j\hat{\mbf x}_i}_{L_2}^2 = \norm{x}_{\R}^2+\norm{(\mcl F\hat{\mbf x})}_{L_2}^2,&
\end{flalign*}
where $(\mcl F\hat{\mbf x}) = \text{col}(\partial_s^0S^0\hat{\mbf x}, \cdots, \partial_s^{N} S^{N}\hat{\mbf x})$.
Recall from Corollary \ref{cor:x_D},
\begin{align*}
(\mcl F\hat{\mbf x}) &= \mcl T_{D}\bmat{v\\\mcl D\hat{\mbf x}}, ~\text{for any}~ \hat{\mbf x}\in X_v~\text{and}~ \{n,\mbf G_{\mathrm b}\}~\text{PIE-compatible}
\end{align*}
where $\mcl T_{D} = \fourpi{\emptyset}{\emptyset}{U_2T(s-a)B_T^{-1}B_v}{U_1,R_{D,1},R_{D,2}}$ is a bounded PI operator. Substituting $v = C_v x$ specifically, we have
\begin{flalign*}
\norm{(\mcl F\hat{\mbf x})}_{L_2}^2 &= \norm{\mcl T_{D}\bmat{C_v x\\\mcl D\hat{\mbf x}}}_{L_2}^2\le \norm{\mcl T_{D}}^2_{\mcl L(\R L_2)}\left((b-a)^2\norm{C_v x}^2_{\R}+\norm{\mcl D\hat{\mbf x}}_{L_2}^2\right)=K_0\norm{\mcl D\hat{\mbf x}}_{L_2}^2+K_1\norm{x}_{\R}^2&
\end{flalign*}
where $K_0 = \norm{\mcl{T}_{D}}^2 _{\mcl L(\R L_2)}$ and $K_1 = K_0(b-a)^2\bar{\sigma}(C_v)^2$. Recall from Theorem \ref{thm:unitary_T}, for any $\hat{\mbf{x}}\in L_2$ and $v\in \R$, $\norm{\hat{\mcl{T}}\hat{\mbf{x}}+\mcl T_v v}_{X^n}^2 = \norm{\hat{\mbf{x}}}_{L_2}^2$.
Then
\begin{flalign*}
\norm{\bmat{x\\\hat{\mbf x}}}_{\R^{n_x}\times W^n}^2 &= \norm{x}_{\R}^2+K_0\norm{\mcl D\hat{\mbf x}}_{L_2}^2+K_1\norm{x}_{\R}^2 = (1+K_1)\norm{x}_{\R}^2++K_0\norm{\hat{\mcl{T}}\mcl D\hat{\mbf x}+\mcl T_v C_v x}_{X^n}^2 \\
&\le (1+K_1)\norm{x}_{\R}^2+K_0\norm{\hat{\mbf x}}_{X^n}^2\le (1+K_0+K_1)\norm{\bmat{x\\\hat{\mbf x}}}_{\R^{n_x} \times X^n}^2.
\end{flalign*}
\end{proof}
\subsection{Proof of \Cref{thm:stability_equivalence}}\label{app:stability_equivalence}
Now that we have established equivalence of the $X^n$ and $W^n$ norms, we may prove that a GPDE model is internally (exponential, Lyapunov, or asymptotically) stable if and only if the associated PIE is internally stable.



\thmStability*
\begin{proof}
Suppose GPDE defined by $\{n, \mbf G_{\mathrm{o}}, \mbf G_{\mathrm b}, \mbf G_{\mathrm{p}}\}$ is exponentially stable. Then, there exist constants $M$, $\alpha>0$ such that for any $\{x^0,\hat{\mbf{x}}^0\}\in \mcl X_{0,0}$, if $\{x$, $\hat{\mbf{x}},z,y,v,r\}$ satisfies the GPDE defined $\{n, \mbf G_{\mathrm o}, \mbf G_{\mathrm b}, \mbf G_{\mathrm p}\}$ with initial condition $\{x^0,\hat{\mbf x}^0\}$ and input $\{0,0\}$, we have
\begin{align*}
\norm{\bmat{x(t)\\\hat{\mbf x}(t)}}_{\R^{n_x}\times W^n}\le M\norm{\bmat{x^0\\\hat{\mbf{x}}^0}}_{\R^{n_x}\times W^n}e^{-\alpha t} \qquad \text{for all} ~t\ge 0.
\end{align*}

For any $\ubar{\mbf x}^0\in \R L_2^{n_x,n_{\hat{\mbf x}}}$, let $\{\ubar{\mbf x},z,y\}$ satisfy the PIE defined by $\mbf G_{\mathrm{PIE}}$ with initial condition $\ubar{\mbf x}^0 \in \R L_2^{n_x, n_{\hat{\mbf x}}}$ and input $\{0,0\}$. Then, from \Cref{thm:equivalence_2_reverse}, $\{x,\hat{\mbf x},z,y,v,r\}$ satisfies the GPDE defined by $\{n, \mbf G_{\mathrm{o}}, \mbf G_{\mathrm b}, \mbf G_{\mathrm{p}}\}$ with initial condition $\bmat{x^0\\\hat{\mbf x}^0} := \mcl T\ubar{\mbf x}^0 \in \mcl X_{0,0}$ and input $\{0,0\}$ for some $v$ and $r$ where $\bmat{x(t)\\\hat{\mbf x}(t)} := \mcl T \ubar{\mbf x}(t)$.
Then, by the exponential stability of the GPDE, we have
\[
\norm{\bmat{x(t)\\\hat{\mbf x}(t)}}_{\R^{n_x}\times W^n}\le M\norm{\bmat{x^0\\\hat{\mbf x}^0}}_{\R^{n_x}\times W^n}e^{-\alpha t} \qquad \text{for all} ~t\ge 0.
\]

Since $\bmat{x(t)\\\hat{\mbf x}(t)} \in \mcl X_{0,0}$ and $\bmat{x^0\\\hat{\mbf x}^0}\in \mcl X_{0,0}$, from \cref{lem:norm_equivalence}, we have
\[
\norm{\bmat{x(t)\\\hat{\mbf x}(t)}}_{\R^{n_x}\times X^n}\le \norm{\bmat{x(t)\\\hat{\mbf x}(t)}}_{\R^{n_x}\times W^n} \quad \text{and}~\norm{\bmat{x^0\\\hat{\mbf x}^0}}_{\R^{n_x}\times W^n}\le c_0\norm{\bmat{x^0\\\hat{\mbf x}^0}}_{\R^{n_x}\times X^n}.
\]
By \cref{thm:unitary_T}, for any $\mbf x\in \R L_2$ we have $\norm{\mbf x}_{\R L_2} = \norm{\mcl T\mbf x}_{\R^{n_x}\times X^n}$. Thus, we have the following:
\begin{align*}
\norm{\ubar{\mbf x}(t)}_{\R L_2} &= \norm{\mcl T\ubar{\mbf x}(t)}_{\R^{n_x}\times X^n} = \norm{\bmat{x(t)\\\hat{\mbf x}(t)}}_{\R^{n_x}\times X^n}\le \norm{\bmat{x(t)\\\hat{\mbf x}(t)}}_{\R^{n_x}\times W^n} \le M\norm{\bmat{x^0\\\hat{\mbf x}^0}}_{\R^{n_x}\times W^n} e^{-\alpha t} \\
&\le c_0M \norm{\bmat{x^0\\\hat{\mbf x}^0}}_{\R^{n_x}\times X^n} e^{-\alpha t} = c_0M\norm{\mcl T \ubar{\mbf x}^0}_{\R^{n_x}\times X^n}e^{-\alpha t} = c_0 M \norm{\ubar{\mbf x}^0}_{\R L_2}e^{-\alpha t}.
\end{align*}
Therefore, the PIE defined by $\mbf G_{\mathrm{PIE}}$ is exponentially stable.

Suppose the PIE defined by $\mbf G_{\mathrm{PIE}}$ is exponentially stable. Then, there exist constants $M$, $\alpha>0$ such that for any $\ubar{\mbf x}^0\in \R L_2^{m,n}$, if $\ubar{\mbf x}$ satisfies the PIE defined by $\{\mbf G_{\mathrm{PIE}}\}$ with initial condition $\ubar{\mbf x}^0$ and input $\{0,0\}$, we have
\begin{align*}
\norm{\ubar{\mbf x}(t)}_{\R L_2}\le M\norm{\ubar{\mbf x}^0}_{\R L_2}e^{-\alpha t} \qquad \text{for all} ~t\ge 0.
\end{align*}

For any $\{x^0,\hat{\mbf x}^0\} \in \mcl X_{0,0}$, let $\{x,\hat{\mbf x},z,y,v,r\}$ satisfy the GPDE defined by $n$ and $\{\mbf G_{\mathrm{o}}, \mbf G_{\mathrm b}, \mbf G_{\mathrm{p}}\}$ with initial condition $\{x^0,\hat{\mbf x}^0\}$ and input $\{0,0\}$. Then, from \Cref{thm:equivalence_2_reverse}, $\{\ubar{\mbf x},z,y\}$ satisfies the PIE defined by $\mbf G_{\mathrm{PIE}}$ with initial condition $\ubar{\mbf x}^0 \in \R L_2^{n_x, n_{\hat{\mbf x}}}$ and input $\{0,0\}$ where
\[
\ubar{\mbf x}(t) = \bmat{x(t)\\\mcl D \hat{\mbf x}(t)}, \qquad \ubar{\mbf x}^0 = \bmat{x^0\\\mcl D\hat{\mbf x}^0}.
\]
Since $\hat{\mbf x}(t)\in X_{C_v x(t)}$, from \cref{thm:T_map}, we have $\hat{\mbf x}(t) = \hat{\mcl T}\mcl D\hat{\mbf x}(t)+\mcl T_v C_v x(t)$. Therefore,
\[
\bmat{x(t)\\\hat{\mbf x}(t)} = \bmat{x(t)\\\hat{\mcl T}\mcl D\hat{\mbf x}(t)+\mcl T_v C_v x(t)} = \bmat{I&0\\\mcl T_v C_v &\hat{\mcl T}}\bmat{x(t)\\\mcl D \hat{\mbf x}(t)} = \mcl T \ubar{\mbf x}(t).
\]
Similarly, we have $\bmat{x^0\\\hat{\mbf x}^0}=\mcl T \ubar{\mbf x}^0$.

By the exponential stability of the PIE, we have
\[
\norm{\ubar{\mbf x}(t)}_{\R L_2}\le M\norm{\ubar{\mbf x}^0}_{\R L_2}e^{-\alpha t} \qquad \text{for all} ~t\ge 0.
\]
Again, from \cref{lem:norm_equivalence}, we have
\[
\norm{\bmat{x(t)\\\hat{\mbf x}(t)}}_{\R^{n_x}\times W^n}\le c_0\norm{\bmat{x(t)\\\hat{\mbf x}(t)}}_{\R^{n_x}\times X^n}, \quad \norm{\bmat{x^0\\\hat{\mbf x}^0}}_{\R^{n_x}\times X^n}\le\norm{\bmat{x^0\\\hat{\mbf x}^0}}_{\R^{n_x}\times W^n}
\]
and, from \cref{thm:unitary_T}, $\norm{\mcl T\mbf x}_{\R^{n_x}\times X^n}=\norm{\mbf x}_{\R L_2}$ for any $\mbf x\in \R L_2$ which implies
\begin{align*}
\norm{\bmat{x(t)\\\hat{\mbf x}(t)}}_{\R^{n_x}\times W^n} &\le c_0 \norm{\bmat{x(t)\\\hat{\mbf x}(t)}}_{\R^{n_x}\times X^n} = c_0 \norm{\mcl T \ubar{\mbf x}(t)}_{\R^{n_x}\times X^n} = c_0\norm{\ubar{\mbf x}(t)}_{\R L_2}\\
& \le c_0M \norm{\ubar{\mbf x}^0}_{\R L_2}e^{-\alpha t}= c_0M \norm{\mcl T \ubar{\mbf x}^0}_{\R^{n_x}\times X^n}e^{-\alpha t}  \\
&= c_0M \norm{\bmat{x^0\\\hat{\mbf x}^0}}_{\R^{n_x}\times X^n}e^{-\alpha t}\le c_0M \norm{\bmat{x^0\\\hat{\mbf x}^0}}_{\R^{n_x}\times W^n}e^{-\alpha t}.
\end{align*}
Therefore, the GPDE defined by $n$ and $\{\mbf G_{\mathrm{o}}, \mbf G_{\mathrm b}, \mbf G_{\mathrm{p}}\}$ is exponentially stable.%
\end{proof}

We can prove equivalence of stability for the two representations by using other notions of stability as well. For e.g., consider the Lyapunov and asymptotic stability of GPDEs and PIEs are defined as follows.\vspace{-3mm}
\begin{definition}[Lyapunov Stability]\;
\begin{enumerate}
	\item We say a GPDE model defined by $\{n, \mbf G_{\mathrm o}, \mbf G_{\mathrm b}, \mbf G_{\mathrm p}\}$ is Lyapunov stable, if for every $\epsilon>0$ there exists a $\delta>0$ such that for any $\{x^0,\hat{\mbf{x}}^0\}\in \mcl X_{0,0}$ with $\norm{\bmat{x^0\\\hat{\mbf x}^0}}_{\R^{n_x}\times W^n}< \delta$, if $\{x$, $\hat{\mbf{x}}, z, y, v, r\}$ satisfies the GPDE defined by $\{n, \mbf G_{\mathrm o}, \mbf G_{\mathrm b}, \mbf G_{\mathrm p}\}$ with initial condition $\{x^0,\hat{\mbf x}^0\}$ and input $\{0,0\}$, then
	\begin{align*}
	\norm{\bmat{x(t)\\\hat{\mbf{x}}(t)}}_{\R^{n_x}\times W^n}< \epsilon  \qquad \text{for all} ~t\ge 0.
	\end{align*}
	\item We say a PIE model defined by $\mbf G_{\mathrm{PIE}}$ is Lyapunov stable if for every $\epsilon>0$ there exists a constant $\delta>0$ such that for any $\ubar{\mbf{x}}^0\in \R L_2^{m,n}$ with $\norm{\ubar{\mbf x}^0}_{\R L_2^{m,n}}<\delta$, if $\{\ubar{\mbf{x}}, z, y\}$ satisfies the PIE defined by $\mbf G_{\mathrm{PIE}}$ with initial condition $\ubar{\mbf x}^0$ and input $\{0,0\}$, then	$\norm{\ubar{\mbf{x}}(t)}_{\R L_2^{m,n}}< \epsilon$ for all $t\ge 0$.
\end{enumerate}
\end{definition}

\begin{definition}[Asymptotic Stability]\;
\begin{enumerate}
	\item We say a GPDE defined by $\{n, \mbf G_{\mathrm o}, \mbf G_{\mathrm b}, \mbf G_{\mathrm p}\}$ is asymptotically stable, if for every $\{x^0,\hat{\mbf{x}}^0\}\in \mcl X_{0,0}$ and $\epsilon>0$, there exists a $T_{\epsilon}>0$ such that if $\{x$, $\hat{\mbf{x}}, z, y, v, r\}$ satisfies the GPDE defined by $\{n, \mbf G_{\mathrm o}, \mbf G_{\mathrm b}, \mbf G_{\mathrm p}\}$ with initial condition $\{x^0,\hat{\mbf x}^0\}$ and input $\{0,0\}$, then $\norm{\bmat{x(t)\\\hat{\mbf{x}}(t)}}_{\R^{n_x}\times W^n} < \epsilon$ for all $t>T_{\epsilon}$.
	\item We say a PIE model defined by $\mbf G_{\mathrm{PIE}}$ is asymptotically stable, if for every $\ubar{\mbf{x}}^0\in \R L_2^{m,n}$ and $\epsilon>0$, there exists a $T_{\epsilon}>0$ such that if $\{\ubar{\mbf{x}}, z, y\}$ satisfies the PIE defined by $\mbf G_{\mathrm{PIE}}$ with initial condition $\ubar{\mbf x}^0$ and input $\{0,0\}$, then there exists $T_{\epsilon}>0$ such that $\norm{\ubar{\mbf{x}}(t)}_{\R L_2^{m,n}}<\epsilon$ for all $t>T_{\epsilon}$.
\end{enumerate}
\end{definition}

For the above notions of stability, we have the following results.
\begin{restatable}{corollary}{corStability}\label{cor:asymptotic_stability}
Given $\{n,\mbf G_{\mathrm{o}}, \mbf G_{\mathrm b}, \mbf G_{\mathrm{p}}\}$ PIE-compatible, let $\mbf G_{\mathrm{PIE}}:=\mbf M(\{n,\mbf G_{\mathrm b},\mbf G_{\mathrm o},\mbf G_{\mathrm p}\})$. Then
\begin{enumerate}
	\item The GPDE model defined by $\{n,\mbf G_{\mathrm{o}}, \mbf G_{\mathrm b}, \mbf G_{\mathrm{p}}\}$ is Lyapunov stable if and only if the PIE system defined by $\mbf G_{\mathrm{PIE}}$ is Lyapunov stable.
	\item The GPDE model defined by $\{n,\mbf G_{\mathrm{o}}, \mbf G_{\mathrm b}, \mbf G_{\mathrm{p}}\}$ is asymptotically stable if and only if the PIE system defined by $\mbf G_{\mathrm{PIE}}$ is asymptotically stable.
\end{enumerate}
\end{restatable}
\begin{proof}
\textbf{Proof of part 1.} Suppose GPDE defined by $\{n,\mbf G_{\mathrm{o}}, \mbf G_{\mathrm b}, \mbf G_{\mathrm{p}}\}$ is Lyapunov stable. For any $\ubar{\mbf x}^0\in \R L_2^{n_x,n_{\hat{\mbf x}}}$, let $\{\ubar{\mbf x},z,y\}$ satisfy the PIE defined by $\mbf G_{\mathrm{PIE}}$ with initial condition $\ubar{\mbf x}^0 \in \R L_2^{n_x, n_{\hat{\mbf x}}}$ and input $\{0,0\}$. Then, from \Cref{thm:equivalence_2_reverse}, $\{x,\hat{\mbf x},z,y,v,r\}$ satisfies the GPDE defined by $n$ and $\{\mbf G_{\mathrm{o}}, \mbf G_{\mathrm b}, \mbf G_{\mathrm{p}}\}$ with initial condition $\bmat{x^0\\\hat{\mbf x}^0} := \mcl T\ubar{\mbf x}^0 \in \mcl X_{0,0}$ and input $\{0,0\}$ for some $v$ and $r$ where $\bmat{x(t)\\\hat{\mbf x}(t)} := \mcl T \ubar{\mbf x}(t)$.
Suppose $\epsilon>0$, then by the Lyapunov stability of the GPDE, there exists $\delta$ such that
\[
\norm{\bmat{x^0\\\hat{\mbf x}^0}}_{\R^{n_x}\times W^n}< \delta\implies \norm{\bmat{x(t)\\\hat{\mbf x}(t)}}_{\R^{n_x}\times W^n}< \epsilon \qquad \text{for all} ~t\ge 0.
\]
Since $\bmat{x(t)\\\hat{\mbf x}(t)}\in \mcl X_{0,0}$ and $\bmat{x^0\\\hat{\mbf x}^0}\in \mcl X_{0,0}$, from \cref{lem:norm_equivalence}, we have
\[
\norm{\bmat{x(t)\\\hat{\mbf x}(t)}}_{\R^{n_x}\times X^n}\le \norm{\bmat{x(t)\\\hat{\mbf x}(t)}}_{\R^{n_x}\times W^n} \quad \text{and}~\norm{\bmat{x^0\\\hat{\mbf x}^0}}_{\R^{n_x}\times W^n}\le c_0\norm{\bmat{x^0\\\hat{\mbf x}^0}}_{\R^{n_x}\times X^n}.
\]
Let $\norm{\ubar{\mbf x}^0}_{\R L_2}< \frac{\delta}{c_0}$. By \cref{thm:unitary_T}, for any $\ubar{\mbf x}\in \R L_2^{n_x,n_{\hat{\mbf x}}}$ we have $\norm{\ubar{\mbf x}}_{\R L_2} = \norm{\mcl T\ubar{\mbf x}}_{\R^{n_x}\times X^n}$. Thus, we have the following:
\begin{align*}
\norm{\bmat{x^0\\\hat{\mbf x}^0}}_{\R^{n_x}\times W^n}&\le c_0\norm{\bmat{x^0\\\hat{\mbf x}^0}}_{\R^{n_x}\times X^n} = c_0 \norm{\mcl T\ubar{\mbf x}^0}_{\R^{n_x}\times X^n} = c_0 \norm{\ubar{\mbf x}^0}_{\R L_2}< \delta, \quad \text{and}\\
\norm{\ubar{\mbf x}(t)}_{\R L_2} &= \norm{\mcl T\ubar{\mbf x}(t)}_{\R^{n_x}\times X^n} = \norm{\bmat{x(t)\\\hat{\mbf x}(t)}}_{\R^{n_x}\times X^n}\le \norm{\bmat{x(t)\\\hat{\mbf x}(t)}}_{\R^{n_x}\times W^n} < \epsilon.
\end{align*}
Therefore, the PIE defined by $\mbf G_{\mathrm{PIE}}$ is Lyapunov stable.

Suppose the PIE defined by $\mbf G_{\mathrm{PIE}}$ is Lyapunov stable. For any $\{x^0,\hat{\mbf x}^0\} \in \mcl X_{0,0}$, let $\{x,\hat{\mbf x},z,y,v,r\}$ satisfy the GPDE defined by $n$ and $\{\mbf G_{\mathrm{o}}, \mbf G_{\mathrm b}, \mbf G_{\mathrm{p}}\}$ with initial condition $\{x^0,\hat{\mbf x}^0\}$ and input $\{0,0\}$. Then, from \Cref{thm:equivalence_2_reverse}, $\{\ubar{\mbf x},z,y\}$ satisfies the PIE defined by $\mbf G_{\mathrm{PIE}}$ with initial condition $\ubar{\mbf x}^0 \in \R L_2^{n_x, n_{\hat{\mbf x}}}$ and input $\{0,0\}$ where
\[
\ubar{\mbf x}(t) = \bmat{x(t)\\\mcl D \hat{\mbf x}(t)}, \qquad \ubar{\mbf x}^0 = \bmat{x^0\\\mcl D\hat{\mbf x}^0}.
\]
Since $\hat{\mbf x}(t)\in X_{C_v x(t)}$, from \cref{thm:T_map}, we have $\hat{\mbf x}(t) = \hat{\mcl T}\mcl D\hat{\mbf x}(t)+\mcl T_v C_v x(t)$. Therefore,
\[
\bmat{x(t)\\\hat{\mbf x}(t)} = \bmat{x(t)\\\hat{\mcl T}\mcl D\hat{\mbf x}(t)+\mcl T_v C_v x(t)} = \bmat{I&0\\\mcl T_v C_v &\hat{\mcl T}}\bmat{x(t)\\\mcl D \hat{\mbf x}(t)} = \mcl T \ubar{\mbf x}(t).
\]
Similarly, we have $\bmat{x^0\\\hat{\mbf x}^0}=\mcl T \ubar{\mbf x}^0$.
Again, from \cref{lem:norm_equivalence}, we have
\[
\norm{\bmat{x(t)\\\hat{\mbf x}(t)}}_{\R^{n_x}\times W^n}\le c_0\norm{\bmat{x(t)\\\hat{\mbf x}(t)}}_{\R^{n_x}\times X^n}, \quad \norm{\bmat{x^0\\\hat{\mbf x}^0}}_{\R^{n_x}\times X^n}\le\norm{\bmat{x^0\\\hat{\mbf x}^0}}_{\R^{n_x}\times W^n}
\]
and, from \cref{thm:unitary_T}, $\norm{\mcl T\mbf x}_{\R^{n_x}\times X^n}=\norm{\mbf x}_{\R L_2}$ for any $\mbf x\in \R L_2^{n_x,n_{\hat{\mbf x}}}$.
Let $\epsilon>0$. Then, by the Lyapunov stability of the PIE, there exists $\delta$ such that
\[
\norm{\ubar{\mbf x}^0}_{\R L_2}< \delta \implies \norm{\ubar{\mbf x}(t)}_{\R L_2}< \frac{\epsilon}{c_0} \qquad \text{for all} ~t\ge 0.
\]
For any initial condition for the GPDE such that $\norm{\bmat{x^0\\\hat{\mbf x}^0}}_{\R^{n_x}\times W^n}< \delta$, we have
\begin{align*}
\norm{\ubar{\mbf x}^0}_{\R L_2} &= \norm{\mcl T \ubar{\mbf x}^0}_{\R^{n_x}\times X^n} = \norm{\bmat{x^0\\\hat{\mbf x}^0}}_{\R^{n_x}\times X^n}\le \norm{\bmat{x^0\\\hat{\mbf x}^0}}_{\R^{n_x}\times W^n}< \delta, \quad \text{and}\\
\norm{\bmat{x(t)\\\hat{\mbf x}(t)}}_{\R^{n_x}\times W^n} &\le c_0 \norm{\bmat{x(t)\\\hat{\mbf x}(t)}}_{\R^{n_x}\times X^n} = c_0 \norm{\mcl T \ubar{\mbf x}(t)}_{\R^{n_x}\times X^n} = c_0\norm{\ubar{\mbf x}(t)}_{\R L_2} < \epsilon.
\end{align*}
Therefore, the GPDE defined by $\{n, \mbf G_{\mathrm{o}}, \mbf G_{\mathrm b}, \mbf G_{\mathrm{p}}\}$ is Lyapunov stable.%

\textbf{Proof of part 2.} Suppose GPDE defined by $\{n,\mbf G_{\mathrm{o}}, \mbf G_{\mathrm b}, \mbf G_{\mathrm{p}}\}$ is asymptotically stable. For any $\ubar{\mbf x}^0\in \R L_2^{n_x,n_{\hat{\mbf x}}}$, let $\{\ubar{\mbf x},z,y\}$ satisfy the PIE defined by $\mbf G_{\mathrm{PIE}}$ with initial condition $\ubar{\mbf x}^0 \in \R L_2^{n_x, n_{\hat{\mbf x}}}$ and input $\{0,0\}$. Then, from \Cref{thm:equivalence_2_reverse}, $\{x,\hat{\mbf x},z,y,v,r\}$ satisfies the GPDE defined by $n$ and $\{\mbf G_{\mathrm{o}}, \mbf G_{\mathrm b}, \mbf G_{\mathrm{p}}\}$ with initial condition $\bmat{x^0\\\hat{\mbf x}^0} := \mcl T\ubar{\mbf x}^0 \in \mcl X_{0,0}$ and input $\{0,0\}$ for some $v$ and $r$ where $\bmat{x(t)\\\hat{\mbf x}(t)} := \mcl T \ubar{\mbf x}(t)$.
Suppose $\epsilon>0$, then by the asymptotic stability of the GPDE, there exists $T_0$ such that
\begin{align*}
\norm{\bmat{x(t)\\\hat{\mbf x}(t)}}_{\R^{n_x}\times W^n}< \epsilon \qquad \text{for all} ~t\ge T_0.
\end{align*}
Since $\bmat{x(t)\\\hat{\mbf x}(t)}\in \mcl X_{0,0}$, from \cref{lem:norm_equivalence}, we have
\[
\norm{\bmat{x(t)\\\hat{\mbf x}(t)}}_{\R^{n_x}\times X^n}\le \norm{\bmat{x(t)\\\hat{\mbf x}(t)}}_{\R^{n_x}\times W^n} .
\]
By \cref{thm:unitary_T}, for any $\ubar{\mbf x}\in \R L_2^{n_x,n_{\hat{\mbf x}}}$ we have $\norm{\ubar{\mbf x}}_{\R L_2} = \norm{\mcl T\ubar{\mbf x}}_{\R^{n_x}\times X^n}$. Thus, for any $t> T_0$, we have,
\begin{align*}
\norm{\ubar{\mbf x}(t)}_{\R L_2} &= \norm{\mcl T\ubar{\mbf x}(t)}_{\R^{n_x}\times X^n} = \norm{\bmat{x(t)\\\hat{\mbf x}(t)}}_{\R^{n_x}\times X^n}\le \norm{\bmat{x(t)\\\hat{\mbf x}(t)}}_{\R^{n_x}\times W^n} < \epsilon.
\end{align*}
Therefore, the PIE defined by $\mbf G_{\mathrm{PIE}}$ is asymptotically stable.

Suppose the PIE defined by $\mbf G_{\mathrm{PIE}}$ is asymptotically stable. For any $\{x^0,\hat{\mbf x}^0\} \in \mcl X_{0,0}$, let $\{x,\hat{\mbf x},z,y,v,r\}$ satisfy the GPDE defined by $\{n,$ $\mbf G_{\mathrm{o}},$ $\mbf G_{\mathrm b},$ $\mbf G_{\mathrm{p}}\}$ with initial condition $\{x^0,\hat{\mbf x}^0\}$ and input $\{0,0\}$. Then, from \Cref{thm:equivalence_2_reverse}, $\{\ubar{\mbf x},z,y\}$ satisfies the PIE defined by $\mbf G_{\mathrm{PIE}}$ with initial condition $\ubar{\mbf x}^0 \in \R L_2^{n_x, n_{\hat{\mbf x}}}$ and input $\{0,0\}$ where
\[
\ubar{\mbf x}(t) = \bmat{x(t)\\\mcl D \hat{\mbf x}(t)}, \qquad \ubar{\mbf x}^0 = \bmat{x^0\\\mcl D\hat{\mbf x}^0}.
\]
Again, we know $\hat{\mbf x}(t)\in X_{C_v x(t)}$, and hence from \cref{thm:T_map}, we have $\hat{\mbf x}(t) = \hat{\mcl T}\mcl D\hat{\mbf x}(t)+\mcl T_v C_v x(t)$. Therefore,
\[
\bmat{x(t)\\\hat{\mbf x}(t)} = \bmat{x(t)\\\hat{\mcl T}\mcl D\hat{\mbf x}(t)+\mcl T_v C_v x(t)} = \bmat{I&0\\\mcl T_v C_v &\hat{\mcl T}}\bmat{x(t)\\\mcl D \hat{\mbf x}(t)} = \mcl T \ubar{\mbf x}(t).
\]
Again, from \cref{lem:norm_equivalence}, we have
\[
\norm{\bmat{x(t)\\\hat{\mbf x}(t)}}_{\R^{n_x}\times W^n}\le c_0\norm{\bmat{x(t)\\\hat{\mbf x}(t)}}_{\R^{n_x}\times X^n},
\]
and, from \cref{thm:unitary_T}, $\norm{\mcl T\mbf x}_{\R^{n_x}\times X^n}=\norm{\mbf x}_{\R L_2}$ for any $\mbf x\in \R L_2^{n_x,n_{\hat{\mbf x}}}$.
Let $\epsilon>0$. Then, by the asymptotic stability of the PIE, there exists $T_0$ such that
\[
\norm{\ubar{\mbf x}(t)}_{\R L_2}< \frac{\epsilon}{c_0} \qquad \text{for all} ~t\ge T_0.
\]
Then, for any $t\ge T_0$, we have
\begin{align*}
\norm{\bmat{x(t)\\\hat{\mbf x}(t)}}_{\R^{n_x}\times W^n} &\le c_0 \norm{\bmat{x(t)\\\hat{\mbf x}(t)}}_{\R^{n_x}\times X^n} = c_0 \norm{\mcl T \ubar{\mbf x}(t)}_{\R^{n_x}\times X^n} = c_0\norm{\ubar{\mbf x}(t)}_{\R L_2} < \epsilon.
\end{align*}
Therefore, the GPDE defined by $\{n, \mbf G_{\mathrm{o}}, \mbf G_{\mathrm b}, \mbf G_{\mathrm{p}}\}$ is asymptotically stable.%

\end{proof}

\subsection{Set Of PI Operators Forms A $^*$-Algebra}\label{app:algebra_pi}
In this section, we prove that set of PI operators when parameterized by $L_{\infty}$-bounded functions forms a $^*$-algebra, i.e., closed algebraically. Furthermore, the formulae provided here will act as a guideline to perform the binary operations (addition, composition, and concatenation) of PI operators since various formulae in the paper were presented using such binary operation notation. First, we provide a formal definition of the list of properties a set must satisfy to be a $^*$-algebra. A $^*$-algebra must be an associative algebra with an involution operation. Since definition of $^*$-algebra depends on definitions of an algebra that is associative, we introduce those definitions first. Since we only use $4$-PI operators in this Subsection, we will drop the subscript and use $\Pi$ instead of $\Pi_4$ --- the results naturally extend to operators in $\Pi_3$.

\begin{definition}[Algebra]
A vector space, $A$, equipped with a multiplication operation is said to be an algebra if for every $X,Y\in A$ we have $XY\in A$.
\end{definition}
\begin{definition}[Associative Algebra]
An algebra, $A$, is said to be associative if for every $X,Y,Z\in A$
\begin{align*}
X(YZ) = (XY)Z
\end{align*}
where $XY$ denotes a multiplication operation between $X$ and $Y$.
\end{definition}

\begin{definition}[$^*$-algebra]
An algebra, $A$, over the $\R$ with an involution operation $^*$ is called a $^*$-algebra if
\begin{enumerate}
	\item $(X^*)^* = X, \quad \forall X\in A$
	\item $(X+Y)^* = X^*+Y^*, \quad \forall X,Y\in A$
	\item $(XY)^* = Y^*X^*, \quad \forall X,Y\in A$
	\item $(\lambda X)^* = \lambda X^*, \quad \forall \lambda\in\R, X\in A$
\end{enumerate}
\end{definition}

To prove that the set of PI operators $\Pi_{q,q}^{p,p}$ satisfy all the above properties, we prove that $\Pi_{q,q}^{p,p}$ satisfies the requirements of each of the above definitions where
\[
\Pi_{q,q}^{p,p} : = \left\lbrace\mat{ \fourpi{P}{Q_1}{Q_2}{R_0,R_1,R_2}\;\mid \; P\in \R^{p\times p},~Q_1(s),Q_2(s)^T\in \R^{p\times q},&\\&\hspace{-8.5cm} R_0(s),R_1(s,\theta),R_2(s,\theta)\in\R^{q\times q},~\text{and}~Q_1,~Q_2,~R_0,~R_1,~R_2\in L_{\infty},\\
&\hspace{-8.5cm} R_1,~R_2~\text{are separable}}\right\rbrace.
\]
Also, recall that any $\Pi_4$ operator has an associated set of matrix and polynomial parameters which lie in the space
\begin{align*}
[\Gamma]^{m,p}_{n,q}&:=\left\lbrace \mat{{\small\left[\begin{array}{c|c}P & Q_1 \\\hline  Q_2 & \{R_0,R_{1a}R_{1b},R_{2a}R_{2b}\}\end{array}\right]} \;: \\\hspace{5mm} P\in \R^{m\times n},\,Q_1 \in L_{\infty}^{m\times q},\,Q_2\in L_{\infty}^{p\times n},\\\hspace{5mm}R_0\in L_{\infty}^{q\times n},\,R_{ia}\in L_{\infty}^{q\times n_b},R_{ib}\in L_{\infty}^{n_b\times n}}\right\rbrace.
\end{align*}
To prove that the set is an algebra, we need to define two binary operations addition and multiplication, which in case of $\Pi_{q,q}^{p,p}$ will be given be addition of PI operators (as defined in \Cref{lem:pi_add}) and composition of PI operators (as defined in \Cref{lem:pi_comp}). For the set to be a $^*$-algebra we also need an involution operation which is given by the adjoint with respect to $\R\times L_2$ inner-product (as defined in \Cref{lem:adj}).

\begin{lemma}[Addition]\label{lem:pi_add}
For any $A,L\in\R^{m\times p}$ and $B_1,M_1: [a, b]\to\R^{m\times q}$, $B_2, M_2: [a, b]\to\R^{n\times p}$, $C_0,N_0:[a, b]\to\R^{n \times q}$, $C_i,N_i:[a, b]\times [a, b] \to \R^{n\times q}$, for $i \in \{1,2\}$, $L_{\infty}$ bounded, define a linear map $\mbf P_+^4:\Gamma_{n,q}^{m,p}\times \Gamma_{n,q}^{m,p}\to \Gamma_{n,q}^{m,p}$ such that
\begin{align*}
\left[\footnotesize\begin{array}{c|c}P&Q_1\\\hline Q_2&\{R_i\}\end{array}\right]=\mbf P_+^4\left(\left[\footnotesize\begin{array}{c|c}A&B_1\\\hline B_2&\{C_i\}\end{array}\right],\left[\footnotesize\begin{array}{c|c}L&M_1\\\hline M_2&\{N_i\}\end{array}\right]\right)
\end{align*}
where
\begin{align*}
P&=A+L, ~Q_i=B_i+M_i, ~R_i=C_i+N_i.
\end{align*}
If $P, Q_i, R_i$ are as defined above, then, for any $x \in \mathbb{R}^{p}$ and $z \in L_2^q([a, b])$
\begin{align*}
&\pie\left[\mbf P_+^4\left(\footnotesize\left[\begin{array}{c|c}A&B_1\\\hline B_2&\{C_i\}\end{array}\right],\left[\begin{array}{c|c}L&M_1\\\hline M_2&\{N_i\}\end{array}\right]\right)\right]\bmat{x\\z}=\left(\fourpi{A}{B_1}{B_2}{C_i}+\fourpi{L}{M_1}{M_2}{N_i}\right)\bmat{x\\z}.
\end{align*}
\end{lemma}
\begin{proof}
Let $x\in\R^p$ and $\mbf y\in L_2^q[a,b]$ be arbitrary. Then
\begin{align*}
&\fourpi{P}{Q_1}{Q_2}{R_i}\bmat{x\\\mbf{y}}(s)= \bmat{Px + \int_{a}^{b}Q_1(s)\mbf{y}(s)ds\\Q_2(s)x+\threepi{R_i}\mbf{y} (s)}\\
&= \bmat{(A+L)x + \int_{a}^{b}(B_1+M_1)(s)\mbf{y}(s)ds\\(B_2+M_2)(s)x+(\threepi{C_i+N_i})\mbf{y}(s)}= \bmat{Ax + \int_{a}^{b}B_1(s)\mbf{y}(s)ds\\B_2(s)x+\threepi{C_i}\mbf{y} (s)}+\bmat{Lx + \int_{a}^{b}M_1(s)\mbf{y}(s)ds\\M_2(s)x+\threepi{N_i}\mbf{y} (s)}\\
&= \fourpi{A}{B_1}{B_2}{C_i}\bmat{x\\\mbf{y}}(s)+\fourpi{L}{M_1}{M_2}{N_i}\bmat{x\\\mbf{y}}(s) = \left(\fourpi{A}{B_1}{B_2}{C_i}+\fourpi{L}{M_1}{M_2}{N_i}\right)\bmat{x\\\mbf{y}}(s).
\end{align*}
\end{proof}

\begin{lemma}[Composition] \label{lem:pi_comp}
For any matrices $A\in\R^{m\times k}, P\in\R^{k\times p}$ and $L_{\infty}$ bounded functions $B_1: [a, b]\to\R^{m\times l}, Q_1:[a,b]\to\R^{k\times q}$, $B_2: [a, b]\to\R^{n\times k}, Q_2:[a,b]\to\R^{l \times p}$, $C_0:[a, b]\to\R^{n \times l}, R_0:[a,b]\to\R^{l \times q}$, $C_i:[a, b]^2\to \R^{n\times l}, R_i:[a,b]^2\to\R^{l\times q}$, for $i \in \{1,2\}$, define a linear map $\mbf P_\times:\Gamma_{n,l}^{m,k}\times \Gamma_{l,q}^{k,p}\to \Gamma_{n,q}^{m,p}$ such that
\begin{align*}
\left[\footnotesize\begin{array}{c|c}\hat{P}&\hat{Q}_1\\\hline \hat{Q}_2&\{\hat{R}_i\}\end{array}\right]=\picomp{A}{B_1}{B_2}{C_i}{P}{Q_1}{Q_2}{R_i}
\end{align*}
where
\begin{align*}
&\hat{P} = AP + \int_a^b B_1(s)Q_2(s)ds,\quad \hat{R}_0(s) = C_0(s)R_0(s),\\
&\hat{Q}_1(s) = AQ_1(s) + B_1(s)R_0(s)+\int_{s}^b B_1(\eta)R_1(\eta,s)d\eta+\myinta{s}B_1(\eta)R_2(\eta,s)d\eta,\\
&\hat{Q}_2(s) = B_2(s)P + C_0(s)Q_2(s) + \int_a^s C_1(s,\eta)Q_2(\eta)d\eta+\int_s^b C_2(s,\eta)Q_2(\eta)d\eta,\\
&\hat{R}_1(s,\eta) =B_2(s)Q_1(\eta)+C_0(s)R_1(s,\eta)+C_1(s,\eta)R_0(\eta)\\
&+\int_a^{\eta} C_1(s,\theta)R_2(\theta,\eta)d\theta+\int_{\eta}^{s}C_1(s,\theta)R_1(\theta,\eta)d\theta+\int_{s}^bC_2(s,\theta)R_1(\theta,\eta)d\theta,\\
&\hat{R}_2(s,\eta) =B_2(s)Q_1(\eta)+C_0(s)R_2(s,\eta)+C_2(s,\eta)R_0(\eta)\\
&+\int_a^{s} C_1(s,\theta)R_2(\theta,\eta)d\theta+\int_{s}^{\eta}C_2(s,\theta)R_2(\theta,\eta)d\theta+\int_{\eta}^bC_2(s,\theta)R_1(\theta,\eta)d\theta.
\end{align*}
If $\hat{P}, \hat Q_i, \hat R_i$ are as defined above, then, for any $x \in \mathbb{R}^{m}$ and $z \in L_2^n([a, b])$,
\begin{align*}
&\pie\left[\picomp{A}{B_1}{B_2}{C_i}{P}{Q_1}{Q_2}{R_i}\right]\bmat{x\\ z}=\fourpi{A}{B_1}{B_2}{C_i}\left(\fourpi{P}{Q_1}{Q_2}{R_i}\bmat{x\\z}\right).
\end{align*}
\end{lemma}
\begin{proof}
Let $\{A,~B_i,~C_i\}$, $\{P,~Q_i,~R_i\}$ and $\{\hat P,~\hat Q_i,~\hat R_i\}$ be such that
\begin{align*}
&\fourpi{A}{B_1}{B_2}{C_i}\left(\fourpi{P}{Q_1}{Q_2}{R_i}\bmat{x_1\\x_2}\right)(s) = \left(\fourpi{\hat{P}}{\hat{Q}_1}{\hat{Q}_2}{\hat{R}}\bmat{x_1\\x_2}\right)(s),
\end{align*}
for any $x_1\in \R^{p}$ and $x_2 \in L_2^{q}[a,b]$. Since PI operators are bounded operators on $\R\times L_2$, we define
\begin{align*}
\bmat{y_1\\y_2(s)}:=\left(\fourpi{P}{Q_1}{Q_2}{R}\bmat{x_1\\x_2}\right)(s).
\end{align*}
Then, by definition of a PI operator,
\begin{align*}
y_1 &= Px_1 + \myint Q_1(s) x_2(s)ds\\ y_2(s)&= Q_2(s)x_1 + R_0(s)x_2(s)+\myinta{s}R_1(s,\eta)x_2(\eta)d\eta+\myintb{s}R_2(s,\eta)x_2(\eta)d\eta.
\end{align*}
Likewise, let us also define
\begin{align*}
\bmat{z_1\\z_2(s)} = \left(\fourpi{A}{B_1}{B_2}{C_i}\bmat{y_1\\y_2}\right)(s)= \fourpi{A}{B_1}{B_2}{C_i}\fourpi{P}{Q_1}{Q_2}{R_i}\bmat{x_1\\x_2}(s),
\end{align*}
which gives us the equations
\begin{align*}
z_1 &= Ay_1 + \myint B_1(s) y_2(s)ds\\
z_2(s)&= B_2(s)y_1 + C_0(s)y_2(s)+\myinta{s}C_1(s,\eta)y_2(\eta)d\eta+\myintb{s}C_2(s,\eta)y_2(\eta)d\eta.
\end{align*} We will try to find a direct map between $x_i$ and $z_i$ by substituting $y_i$ in the above equation, however, we will perform the substitution by taking one term at a time. First,
\begin{align*}
\myint B_1(s) y_2(s) ds &= \myint B_1(s) \Big(Q_2(s)x_1 + R_0(s)x_2(s)+\myinta{s}R_1(s,\eta)x_2(\eta)d\eta +\myintb{s}R_2(s,\eta)x_2(\eta)d\eta\Big)ds.
\end{align*}
Then
\begin{align*}
z_1 &= Ay_1 + \myint B_1(s) y_2(s)ds \\
&=
APx_1+\int_a^b AQ_1(s)x_2(s)ds+ \myint B_1(s) \Big(Q_2(s)x_1 + R_0(s)x_2(s)\\
&\qquad+\myinta{s}R_1(s,\eta)x_2(\eta)d\eta +\myintb{s}R_2(s,\eta)x_2(\eta)d\eta\Big)ds\\ &=\hat{P}x_1+\myint \hat{Q}_1(s) x_2(s) ds.
\end{align*}
Next, we substitute $y_i$ in the map from $y_i$ to $z_2(s)$ to get
\begin{align*}
z_2(s) &= B_2(s)Px_1 + \myint B_2(s)Q_1(\eta) x_2(s)d\eta + C_0(s)Q_2(s)x_1 + C_0(s)R_0(s)x_2(s)\\
&+\myinta{s}C_0(s)R_1(s,\eta)x_2(\eta)d\eta+\myintb{s}C_0(s)R_2(s,\eta)x_2(\eta)d\eta+\myinta{s}C_1(s,\eta)Q_2(\eta)x_1d\eta \\
&+ \myinta{s}C_1(s,\eta)R_0(s)x_2(s)d\eta+\myinta{s}\myinta{\eta}C_1(s,\eta)R_1(\eta,\beta)x_2(\beta)d\beta d\eta\\
&+\myinta{s}\myintb{\eta}C_1(s,\eta)R_2(\eta,\beta)x_2(\beta)d\beta d\eta+\myintb{s}C_2(s,\eta)Q_2(\eta)x_1d\eta \\
&+ \myintb{s}C_2(s,\eta)R_0(\eta)x_2(\eta)d\eta+\myintb{s}\myinta{\eta}C_2(s,\eta)R_1(\eta,\beta)x_2(\beta)d\beta d\eta\\&+\myintb{s}\myintb{\eta}C_2(s,\eta)R_2(\eta,\beta)x_2(\beta)d\beta d\eta.
\end{align*}
Next, we separate the terms by factoring $x_1$. Then, we change the order of integration in the double integrals (and swap the variable $\beta\leftrightarrow\eta$) to get
{\small
\begin{align*}
z_2(s)&= \Big(B_2(s)P + C_0(s)Q_2(s) + \myinta{s}C_1(s,\eta)Q_2(\eta)\text{d}\eta\myintb{s}+C_2(s,\eta)Q_2(\eta)\text{d}\eta\Big)x_1 \\
&+  C_0(s)R_0(s)x_2(s)+\myint B_2(\eta)Q_1(s) x_2(\eta)\text{d}\eta +\myinta{s}C_0(s)R_1(s,\eta)x_2(\eta)\text{d}\eta\\
&+\myintb{s}C_0(s)R_2(s,\eta)x_2(\eta)\text{d}\eta+ \myintb{s}C_2(s,\eta)R_0(\eta)x_2(\eta)\text{d}\eta+ \myinta{s}C_1(s,\eta)R_0(s)x_2(s)\text{d}\eta\\&+\myinta{s}\left(\int_a^{\eta} C_1(s,\theta)R_2(\theta,\eta)\text{d}\theta+\int_{\eta}^{s}C_1(s,\theta)R_1(\theta,\eta)\text{d}\theta+\int_{s}^bC_2(s,\theta)R_1(\theta,\eta)\text{d}\theta\right)x_2(\eta)\text{d}\eta \\
&+\myintb{s}\left(\int_a^{s} C_1(s,\theta)R_2(\theta,\eta)\text{d}\theta+\int_{s}^{\eta}C_2(s,\theta)R_2(\theta,\eta)d\theta+\int_{\eta}^bC_2(s,\theta)R_1(\theta,\eta)\text{d}\theta\right)x_2(\eta)\text{d}\eta\\
&= \hat{Q}_2(s) x_1+\hat{S}(s) x_2(s)+\myinta{s} \hat{R_1}(s,\eta)x_2(\eta)\text{d}\eta+ \myintb{s} \hat{R_2}(s,\eta)x_2(\eta) \text{d}\eta.
\end{align*}
}
This completes the proof.
\end{proof}

\begin{lemma} \label{lem:adj}
For any matrices $P\in\R^{m\times p}$ and $L_{\infty}$-bounded functions $Q_1: [a, b]\to\R^{m\times q}$, $Q_2: [a, b]\to\R^{n\times p}$, $R_0:[a, b]\to\R^{n \times q}$, $R_1, R_2:[a, b]\times [a, b] \to \R^{n\times n}$, define a linear map $\mbf P_*^4:\Gamma_{n,q}^{m,p}\to \Gamma_{q,n}^{p,m}$ such that
\begin{align*}
\left[\footnotesize\begin{array}{c|c}\hat P&\hat Q_1\\\hline \hat Q_2&\{\hat R_i\}\end{array}\right]=\mbf P_*^4\left(\left[\footnotesize\begin{array}{c|c}P&Q_1\\\hline Q_2&\{R_i\}\end{array}\right]\right)
\end{align*}
where
\begin{align}
\label{adjoint_matrix_pf}
&\hat{P} = P^{\top}, &&\hat{Q}_1(s) = Q_2^{\top}(s),&&\hat{Q}_2(s) = Q_1^{\top}(s), \nonumber\\
&\hat{R}_0(s) = R_0^{\top}(s),&&\hat{R}_1(s,\eta) = R_2^{\top}(\eta,s), &&\hat{R}_2(s,\eta) = R_1^{\top}(\eta,s).
\end{align}
Then, for any $\mbf x \in \R L_2^{m,n}, \mbf y \in \R L_2^{p,q}$, then we have
\begin{align}
&\ip{\mbf x}{\fourpi{P}{Q_1}{Q_2}{R_i}\mbf y}_{\R L_2^{m,n}}=\ip{\pie\left[\mbf P^4_*\left(\left[\footnotesize\begin{array}{c|c}P&Q_1\\\hline Q_2&\{R_i\}\end{array}\right]\right)\right]\mbf x}{\mbf y}_{\R L_2^{p,q}},
\end{align}
\end{lemma}
\begin{proof}
To prove this, we use the fact that for any scalar $a$ we have $a=a^{\top}$. Let $\mbf x(s)=\bmat{x_1\\\mbf x_2(s)}$ and $\mbf y=\bmat{y_1\\\mbf y_2(s)}$.
Then
\begin{align*}
\ip{\mbf x}{\fourpi{P}{Q_1}{Q_2}{R_i}\mbf y}_{\R L_2^{m,n}}&= x_1^{\top}Py_1+\myint x_1^{\top}Q_1(s)\mbf y_2(s)\text{d}s+\myint \mbf x_2^{\top}(s) Q_2(s) y_1 \text{d}s \\
&~~~+ \myint \mbf x_2(s)^{\top} R_0(s) \mbf y_2(s) \text{d}s + \myint\myinta{s}\mbf x_2^{\top}(s) R_1(s,\eta)\mbf y_2(\eta) \text{d}\eta \text{d}s\\
&~~~+\myint\myintb{s}\mbf x_2^{\top}(s) R_2(s,\eta)\mbf y_2(\eta) \text{d}\eta \text{d}s\\
&~~= y_1^{\top}P^{\top}x_1+\myint y_1^{\top}Q_2^{\top}(s)\mbf x_2(s)\text{d}s+\myint \mbf y_2(s) Q_1^{\top}(s) x_1 \text{d}s \\
&~~~+ \myint \mbf y_2^{\top}(s) R_0^{\top}(s) \mbf x_2(s) \text{d}s + \myint\myinta{s}\mbf y_2^{\top}(s) R_2^{\top}(\eta,s)\mbf x_2(\eta) \text{d}\eta \text{d}s\\
&~~~+\myint\myintb{s}\mbf y_2^{\top}(s) R_1^{\top}(\eta,s)\mbf x_2(\eta) \text{d}\eta \text{d}s\\
&~~= y_1^{\top}\hat{P}x_1+\myint y_1^{\top}\hat{Q}_1(s)\mbf x_2(s)\text{d}s+\myint \mbf y_2^{\top}(s) \hat{Q}_2(s) x_1 \text{d}s \\
&~~~+ \myint \mbf y_2^{\top}(s) \hat{R}_0(s)\mbf  x_2(s) \text{d}s + \myint\myinta{s}\mbf y_2^{\top}(s) \hat{R}_1(s,\eta)\mbf x_2(\eta) \text{d}\eta \text{d}s\\
&~~~+\myint\myintb{s}\mbf y_2^{\top}(s) \hat{R}_2(s,\eta)\mbf x_2(\eta) \text{d}\eta \text{d}s\\
&= \ip{\mbf y}{\fourpi{\hat{P}}{\hat{Q}_1}{\hat Q_2}{\hat{R}_i}\mbf x}_{\R L_2^{p,q}}= \ip{\fourpi{\hat{P}}{\hat{Q}_1}{\hat Q_2}{\hat{R}_i}\mbf x}{\mbf y}_{\R L_2^{p,q}}
\end{align*}
where,
\begin{align*}
&\hat{P} = P^{\top},  &&\hat{Q}_1(s) = Q_2^{\top}(s), &&\hat{Q}_2(s) = Q_1^{\top}(s),\\
&\hat{R}_0(s) = R_0^{\top}(s), &&\hat{R}_1(s,\eta) = R_2^{\top}(\eta,s),  &&\hat{R}_2(s,\eta) = R_1^{\top}(\eta,s).
\end{align*}
This completes the proof.
\end{proof}

Now that we have formally defined the binary and involution operations on the set of PI operators, we show that $\Pi_{q,q}^{p,p}$ when equipped with these operations forms a $^*$-algebra.

\begin{lemma}\label{lem:pi_algebra}
The set $\Pi_{q,q}^{p,p}$ equipped with composition operation forms an associative algebra.
\end{lemma}
\begin{proof}
Suppose $\fourpi{P}{Q_1}{Q_2}{R_i}, \fourpi{A}{B_1}{B_2}{C_i} \in \Pi_{q,q}^{p,p}$. From \Cref{lem:pi_comp}, we have that $\fourpi{\hat P}{\hat Q_1}{\hat Q_2}{\hat R_i} = \fourpi{A}{B_1}{B_2}{C_i}\fourpi{P}{Q_1}{Q_2}{R_i}$ with
\begin{align*}
&\hat{P} = AP + \int_a^b B_1(s)Q_2(s)\text{d}s,~\hat{R}_0(s) = C_0(s)R_0(s),\\
&\hat{Q}_1(s) = AQ_1(s) + B_1(s)R_0(s)+\int_{s}^b B_1(\eta)R_1(\eta,s)\text{d}\eta+\myinta{s}B_1(\eta)R_2(\eta,s)\text{d}\eta,\\
&\hat{Q}_2(s) = B_2(s)P + C_0(s)Q_2(s) + \int_a^s C_1(s,\eta)Q_2(\eta)\text{d}\eta+\int_s^b C_2(s,\eta)Q_2(\eta)\text{d}\eta,\\
&\hat{R}_1(s,\eta) =B_2(s)Q_1(\eta)+C_0(s)R_1(s,\eta)+C_1(s,\eta)R_0(\eta)\\
&+\int_a^{\eta} C_1(s,\theta)R_2(\theta,\eta)\text{d}\theta+\int_{\eta}^{s}C_1(s,\theta)R_1(\theta,\eta)\text{d}\theta+\int_{s}^bC_2(s,\theta)R_1(\theta,\eta)\text{d}\theta,\\
&\hat{R}_2(s,\eta) =B_2(s)Q_1(\eta)+C_0(s)R_2(s,\eta)+C_2(s,\eta)R_0(\eta)\\
&+\int_a^{s} C_1(s,\theta)R_2(\theta,\eta)\text{d}\theta+\int_{s}^{\eta}C_2(s,\theta)R_2(\theta,\eta)d\theta+\int_{\eta}^bC_2(s,\theta)R_1(\theta,\eta)\text{d}\theta.
\end{align*}
Since $B_i, C_i, Q_i, R_i$ are all $L_{\infty}$ we have $\hat Q_i, \hat R_i\in L_{\infty}$. Thus, composition of any two PI operators in $\Pi_{q,q}^{p,p}$ is a PI operator in the same set.

Similarly, by using composition formulae from \Cref{lem:pi_comp}, we can show that for any 3 PI operators $\mcl P, \mcl Q, \mcl R \in \Pi_{q,q}^{p,p}$ we have $(\mcl P\mcl Q)\mcl R = \mcl P(\mcl Q \mcl R)$. The steps are omitted here since the proof is a straightforward arithmetic exercise. Thus $\Pi_{q,q}^{p,p}$ is an associative algebra.
\end{proof}
So far we have shown that the set $\Pi_{q,q}^{p,p}$ is closed algebraically, i.e., the binary and involution operations on PI operators also result in PI operators. In the following Lemma, we conclude that $\Pi_{q,q}^{p,p}$ is a $*$-algebra.


\begin{lemma}
The set $\Pi_{q,q}^{p,p}$ equipped with the binary operations of addition and composition along with the involution operation given by the adjoint with respect to $\R\times L_2$ inner product is a $^*$-algebra.
\end{lemma}
\begin{proof}
To prove this, we first show that $\Pi_{q,q}^{p,p}$ when equipped with the adjoint operator satisfies the requirements of a $^*$-algebra.
Since PI operators are operators on a Hilbert space $\R\times L_2$, from Propositions 2.6 and 2.7 in \cite[p~.32]{conway}, we know that for any two such operators $\mcl P$ and $\mcl Q$
\begin{itemize}
	\item $(\mcl P^*)^* = \mcl P$
	\item $(\lambda \mcl P)^* = \lambda \mcl P^*$
	\item $(\mcl P+\mcl Q)^* = \mcl P^*+\mcl Q^*$
	\item $(\mcl P\mcl Q)^* = \mcl Q^*\mcl P^*$.
\end{itemize}
Therefore, since $\Pi_{q,q}^{p,p}$ is a Banach algebra with an involution $^*$ that satisfies all the properties in the definition of a $^*$-algebra, $\Pi_{q,q}^{p,p}$ is a $^*$-algebra.
\end{proof}

\subsection{Concatenation properties used in the paper}\label{app:pi_concat}
The results presented in this subsection are specific to the notational convenience granted by concatenation of PI operators. Note that in this subsection, we assume that two vectors $\mbf x, \mbf y \in \R L_2$ are identical if there exists a permutation matrix $P$ such that $\mbf x = P\mbf y$. This assumption is made to accommodate for the notational convenience that concatenation of PI operators provide because any $\Pi^{m,p}_{n,q}$ PI operator requires inputs to be completely segregated with finite-dimensional part of the vector to be on the top while infinite-dimensional part at the bottom. However, since concatenation of such vectors is likely to lose such a segregation, we think of the vector $\mbf x ,\mbf y\in \R L_2$ as ordered pairs $(x,\mbf x_1), (y,\mbf y_1)$ with $x, y\in \R$ and $ \mbf x_1, \mbf y_1 \in L_2$ with concatenation of two such vectors being performed individually on each element of the ordered pair. This allows us to retain the convenient segregation of finite and infinite dimensional parts of the vector and use concatenation notation of PI operators.

\begin{lemma}[Horizontal concatenation]\label{lem:pi_con}
Suppose $A_j\in\R^{m\times p_j}$ and $B_{1,j}: [a, b]\to\R^{m\times q_j}$, $B_{2,j}: [a, b]\to\R^{n\times p_j}$, $C_{0,j}:[a, b]\to\R^{n \times q_j}$, $C_{i,j}:[a, b]\times [a, b] \to \R^{n\times q_j}$, for $i \in \{0,1,2\},~ j\in\{1,2\}$, are bounded functions. If we define $P$, $Q_1$, $Q_2$ and $R_k$, for $k\in\{0,1,2\}$ as
\begin{align*}
P&=\bmat{A_1&A_2}, ~Q_i=\bmat{B_{i,1}&B_{i,2}}, ~R_i=\bmat{C_{i,1}&C_{i,2}},
\end{align*} then
\[
\fourpi{P}{Q_1}{Q_2}{R_i}=\bmat{\fourpi{A_1}{B_{1,1}}{B_{2,1}}{C_{i,1}}&\fourpi{A_2}{B_{1,2}}{B_{2,2}}{C_{i,2}}}.
\]
\end{lemma}
\begin{proof}
We will prove this identity by a series of equalities. Let $x_1\in\R^{p_1}$, $y_1 \in \R^{p_2}$, $\mbf x_2\in L_2^{q_1}[a,b]$, and $\mbf y_2 \in L_2^{q_2}[a,b]$ be arbitrary. Next, we define $z_1 = \bmat{x_1\\x_2}\in\R^{p_1+p_2}$ and $\mbf z_2 \in L_2^{q_1+q_2}$. Then the following series of equalities hold. We can substitute $\{z_1,~\mbf z_2\}$ in terms of $\{x_1,~y_1,~\mbf x_2,~\mbf y_2\}$ and perform matrix multiplication to get
\begin{align*}
\fourpi{P}{Q_1}{Q_2}{R_i}\bmat{z_1\\\mbf z_2}(s) 
&= \bmat{Pz_1 +\int_a^b Q_1(s) \mbf z_2(s) ds \\Q_2(s)z_1 + R_0(s)\mbf z_2(s)+\int_a^s R_1(s,\theta) \mbf z_2(\theta)d\theta +\int_s^b R_2(s,\theta)\mbf z_2(\theta)d\theta}\\
&= \bmat{P\bmat{x_1\\y_1} +\int_a^b Q_1(s)\bmat{\mbf x_2(s)\\\mbf y_2(s)} ds \\ Q_2(s)\bmat{x_1\\y_1} + R_0(s)\bmat{\mbf x_2(s)\\\mbf y_2(s)}+\int_a^s R_1(s,\theta) \bmat{\mbf x_2(\theta)\\\mbf y_2(\theta)}d\theta +\int_s^b R_2(s,\theta)\bmat{\mbf x_2(\theta)\\\mbf y_2(\theta)}d\theta}\\
&= \bmat{A_1x_1 +\int_a^b B_{0,1}(s)\mbf x_2(s) ds \\ B_{2,1}(s)x_1 + C_{0,1}(s)\mbf x_2(s)+\int_a^s C_{1,1}(s,\theta) \mbf x_2(\theta)d\theta +\int_s^b C_{2,1}(s,\theta)\mbf x_2(\theta)d\theta}\\
&\qquad+ \bmat{A_2y_1 +\int_a^b B_{0,2}(s)\mbf y_2(s) ds \\ B_{2,2}(s)y_1 + C_{0,2}(s)\mbf y_2(s)+\int_a^s C_{1,2}(s,\theta) \mbf y_2(\theta)d\theta +\int_s^b C_{2,2}(s,\theta)\mbf y_2(\theta)d\theta}\\
&= \fourpi{A_1}{B_{1,1}}{B_{2,1}}{C_{i,1}}\bmat{x_1\\\mbf x_2}(s)+\fourpi{A_2}{B_{1,2}}{B_{2,2}}{C_{i,2}}\bmat{y_1\\\mbf y_2}(s)\\
&= \bmat{\fourpi{A_1}{B_{1,1}}{B_{2,1}}{C_{i,1}}&\fourpi{A_2}{B_{1,2}}{B_{2,2}}{C_{i,2}}}\bmat{x_1\\\mbf x_2\\y_1\\\mbf y_2}(s).
\end{align*}
By rearranging the vector $\text{col}(x_1,\mbf x_2, y_1, \mbf y_2)$ we can obtain $\{z_1, \mbf z_2(s)\}$. Thus the horizontal concatenation of two PI maps gives rise to another uniquely defined PI map.
\end{proof}
Note that in the last equality permutation of rows of the vector is needed to obtain $\{z_1,~\mbf z_2\}$ back. However, that does not affect the conversion formulae for PIE since states can be arranged in any order based on convenience.

\begin{lemma}[Vertical concatenation]\label{lem:pi_con_v}
Suppose $A_j\in\R^{m_j\times p}$ and $B_{1,j}: [a, b]\to\R^{m_j\times q}$, $B_{2,j}: [a, b]\to\R^{n_j\times p}$, $C_{0,j}:[a, b]\to\R^{n_j \times q}$, $C_{i,j}:[a, b]\times [a, b] \to \R^{n_j\times q}$, for $i \in \{0,1,2\},~ j\in\{1,2\}$, are bounded functions. If we define $P$, $Q_1$, $Q_2$ and $R_k$, for $k\in\{0,1,2\}$ as
\begin{align*}
P&=\bmat{A_1\\A_2}, ~Q_i=\bmat{B_{i,1}\\B_{i,2}}, ~R_i=\bmat{C_{i,1}\\C_{i,2}},
\end{align*} then
\[
\fourpi{P}{Q_1}{Q_2}{R_i}=\bmat{\fourpi{A_1}{B_{1,1}}{B_{2,1}}{C_{i,1}}\\\fourpi{A_2}{B_{1,2}}{B_{2,2}}{C_{i,2}}}.
\]
\end{lemma}

\begin{proof}
Similar to horizontal concatenation, we will prove this identity by a series of equalities. Let $x_1\in\R^{p}$ and $\mbf x_2\in L_2^{q}[a,b]$ be arbitrary. Then the following series of equalities hold. We can substitute $\{P,~Q_i,~R_i\}$ in terms of $\{A_j,~B_{i,j},~C_{i,j}\}$ and perform matrix multiplication to get
{\small
\begin{align*}
\fourpi{P}{Q_1}{Q_2}{R_i}\bmat{x_1\\\mbf x_2}(s) &= \bmat{Px_1 +\int_a^b Q_1(s) \mbf x_2(s) ds \\Q_2(s)z_1 + R_0(s)\mbf x_2(s)+\int_a^s R_1(s,\theta) \mbf x_2(\theta)d\theta +\int_s^b R_2(s,\theta)\mbf x_2(\theta)d\theta}\\
&= \bmat{\bmat{A_1\\A_2}x_1 +\int_a^b \bmat{B_{1,1}(s)\\B_{1,2}(s)}\mbf x_2(s) ds \\ \bmat{B_{2,1}(s)\\B_{2,2}(s)}x_1 + \bmat{C_{0,1}(s)\\C_{0,2}(s)}\mbf x_2(s)+\int_a^s \bmat{C_{1,1}(s,\theta)\\C_{1,2}(s,\theta)}\mbf x_2(\theta)d\theta +\int_s^b \bmat{C_{2,1}(s,\theta)\\C_{2,2}(s,\theta)}\mbf x_2(\theta)d\theta}\\
&= \bmat{\bmat{A_1\\A_2}x_1 +\int_a^b \bmat{B_{1,1}(s)\\B_{1,2}(s)}\mbf x_2(s) ds \\ \bmat{B_{2,1}(s)\\B_{2,2}(s)}x_1 + \bmat{C_{0,1}(s)\\C_{0,2}(s)}\mbf x_2(s)+\int_a^s \bmat{C_{1,1}(s,\theta)\\C_{1,2}(s,\theta)}\mbf x_2(\theta)d\theta +\int_s^b \bmat{C_{2,1}(s,\theta)\\C_{2,2}(s,\theta)}\mbf x_2(\theta)d\theta}.
\end{align*}}
By rearranging the rows of the above vector, we get
\begin{align*}
&\bmat{\bmat{A_1x_1 +\int_a^b B_{0,1}(s)\mbf x_2(s) ds \\ B_{2,1}(s)x_1 + C_{0,1}(s)\mbf x_2(s)+\int_a^s C_{1,1}(s,\theta) \mbf x_2(\theta)d\theta +\int_s^b C_{2,1}(s,\theta)\mbf x_2(\theta)d\theta}\\\bmat{A_2x_1 +\int_a^b B_{0,2}(s)\mbf x_2(s) ds \\ B_{2,2}(s)x_1 + C_{0,2}(s)\mbf x_2(s)+\int_a^s C_{1,2}(s,\theta) \mbf x_2(\theta)d\theta +\int_s^b C_{2,2}(s,\theta)\mbf x_2(\theta)d\theta}}\\
&= \bmat{\fourpi{A_1}{B_{1,1}}{B_{2,1}}{C_{i,1}}\bmat{x_1\\\mbf x_2}(s)\\\fourpi{A_2}{B_{1,2}}{B_{2,2}}{C_{i,2}}\bmat{x_1\\\mbf x_2}(s)}\\
&= \bmat{\fourpi{A_1}{B_{1,1}}{B_{2,1}}{C_{i,1}}\\\fourpi{A_2}{B_{1,2}}{B_{2,2}}{C_{i,2}}}\bmat{x_1\\\mbf x_2}(s).
\end{align*}
Thus the vertical concatenation of two PI maps gives rise to another uniquely defined PI map.
\end{proof}

\subsection{Additional Examples}\label{app:more_examples}
In this section, we present additional examples explaining the process of identification of GPDE parameters and finding PIE representation to illustrate the PIE representation for a wide variety PDE systems.

\begin{illus} [ODE coupled with PDE at the Boundary]
In this example, we consider a thermal reactor, $T_r(t)$, which is modeled as an ODE and which is coupled to a cooling jacket, $T_c(t,s)$ which is modeled as a PDE. The dynamics of the reactor and jacket are given by
\begin{align}\label{eqn:example3_PDE_2}
\dot{T}_r(t) &= \lambda T_r(t)-C(T_r(t) - T_c(t,0)),\notag\\
\dot{T}_c(t,s) &= k\partial_{s}^2 T_c(t,s), \qquad s\in (0,1),\notag\\
T_c(t,0) &= T_r(t), \qquad\qquad \partial_{s}T_c(t,1) = 0
\end{align}
where $\lambda$ is the reaction coefficient of the reactor, $C$ is the specific heat of the reactor, and $k$ is a diffusivity parameter for the coolant. In this case, we first model the ODE, where the influence of the PDE on the ODE is isolated in the signal $r(t)=T_c(t,0)$ and the influence of the PDE on the ODE is isolated in the signal $v(t)=T_r(t)$. The state of the ODE subsystem is $x(t)=T_r(t)$ with the following dynamics.
\begin{align}
\dot{x}(t) &= (\lambda-C) x(t)+ r(t),\qquad v(t) = x(t),\qquad
\dot{T}_c(t,s) &= k\partial_{s}^2 T_c(t,s), \qquad s\in (0,1),\notag\\
T_c(t,0) &= v(t), \qquad\qquad \partial_{s}T_c(t,1) = 0
\end{align}

Examining the ODE dynamics
\begin{flalign*}
\bmat{\dot{x}(t)\\\hline z(t)\\y(t)\\v(t)} &= \bmat{A&\vlines&B_{xw}&B_{xu} & B_{xr}\\\hline C_z&\vlines&D_{zw}&D_{zu}&D_{zr}\\C_y&\vlines&D_{yw}&D_{yu}&D_{yr}\\C_v&\vlines&D_{vw}&D_{vu}&0}\bmat{x(t)\\\hline w(t)\\u(t)\\r(t) }.\notag\\
\end{flalign*}
we may parameterize the ODE subsystem, $\mbf G_{\mathrm o}$ as
\[
\mbf G_{\mathrm o}\;:\;\qquad A = \lambda -C, \qquad B_{xr} = C, \qquad C_v = 1.
\]
Now, examining the PDE subsystem, we have a system similar to Illustration in Section~\ref{ill:entropy-PIE} so that the continuity parameter is
\[
n\;:\; \qquad n=\{0,0,1\}\qquad N=2
\]
with $\hat{\mbf x}_2(t,s)=T_c(t,s)$. Again, the BCs appear in the form
\begin{align*}		
&0=\int_{a}^{b} B_{I}(s)\mcl F\hat{\mbf x}(t,s)ds +\bmat{B_v & -B}\bmat{v(t)\\ (\mcl B \hat{\mbf x})(t)}\\
&=\int_{0}^{1} B_{I}(s)\bmat{T_c(t,s)\\T_{c,s}(t,s)\\T_{c,ss}(t,s)}ds -B\bmat{T_c(t,0)\\ T_{c,s}(t,0)\\T_c(t,1)\\ T_{c,s}(t,1)}+B_v v(t).
\end{align*}

\begin{align}
v(t) = x(t)
\dot{T}_c(t,s) &= k\partial_{s}^2 T_c(t,s), \qquad s\in (0,1),\notag\\
T_c(t,0) &= v(t), \qquad\qquad \partial_{s}T_c(t,1) = 0
\end{align}

By inspection of the BCs, we may now define the parameters for $\mbf G_{\mathrm b}$  as
\[
\mbf G_{\mathrm b}\;:\;\qquad B = \bmat{1&0&0&0\\0&0&0&1}, \qquad B_v = \bmat{1\\0}.
\]
To define the parameters of the PDE dynamics, we again ignore integral terms, yielding
\begin{flalign*}
\bmat{\dot{\hat{\mbf{x}}}(t,s)\\ r(t)} =	\bmat{T_c(t,s)\\ r(t)} &=  \bmat{A_{0}(s)\bmat{T_c(t,s)\\T_{c,s}(t,s)\\T_{c,ss}(t,s)}\\
0	} +\bmat{B_{xv}(s) & B_{xb}(s)\\ 0& D_{rb} }\bmat{v(t)\\ (\mcl B \hat{\mbf x})(t)}.
\end{flalign*}
By inspection of Eq.~\eqref{eqn:example3_PDE}, the only non-zero parameter in this expression is
\[	
\mbf G_{\mathrm p}\;:\;\qquad A_{0} = \bmat{0&0&k}
\]
which becomes the entire parameter set for $\mbf G_{\mathrm p}$.

\end{illus}

\begin{illus} [Second Order Time Derivatives]\label{ill:waveequation}
For our next illustration, we consider wave motion
\begin{align}
&\ddot{\eta}(t,s) = c^2\partial_{s}^2\eta(t,s), \qquad\qquad s\in (0,1),\notag\\
&z(t,s) = \int_0^1 \eta(t,s) ds,\notag
\end{align}
where $z$ is a regulated output (the average displacement of the string)	with a general form of BCs (Sturm-Liouville type BCs) given by
\[
\eta(t,0) -k \partial_{s}\eta(t,0) =0, \qquad \eta(t,1)+ l \partial_{s}\eta(t,1) =w(t),
\]
where $\eta$ stands for lateral displacement, $c$ is the speed of propagation of a wave in the string, and $w$ is external disturbance acting on the boundary. The constants $k$ and $l$ represent reflection and mirroring of the wave at the boundary.

To rewrite this PDE model as a state-space GPDE model, we must first eliminate the second order time-derivative. As is common in state-space representation of ODEs, we eliminate the 2nd order time-derivative by creating a new state $\zeta_2 = \dot{\eta}$ with $\zeta_1 = \eta$. This change of variable leads to a coupled PDE of the form
\begin{align}
\dot{\zeta}_1(t,s) &= \zeta_2,\qquad\qquad s\in(0,1),\notag\\
\dot{\zeta}_2(t,s) &= c^2\partial_{s}^2\zeta_1(t,s),\notag\\
z(t,s) &= \int_0^1 \zeta_1(t,s) ds,\label{eqn:example2_PDE}
\end{align}
with BCs
\begin{align}\label{eqn:example2_BC}
\zeta_1(t,0)-k\partial_{s}\zeta_1(t,0)=0,\qquad \zeta_1(t,1)+l\partial_{s}\zeta_1(t,1)=w(t).
\end{align}

Here we note that, the ODE subsystem has the parameters related to outputs $z$ and inputs $w$, however, there is no ODE state. Thus, we only have parameters related to $z$ and $w$. First, we include the influence of PDE on the ODE into the interconnection signal as $r(t) = \int_0^1\zeta_1(t,s)ds$ whereas the influence of the ODE on the PDE is routed through $v$ where $v(t) = w(t)$. Then, by inspection, the output $z$ can be written as $z(t) = r(t)$. Consequently, we find that $D_{zr}=1$, while the remaining parameters related to $z$ are zero. Likewise, we note that $D_{vw}=1$ and leave the remaining parameters of $v$ as empty. This completes the definition of the ODE subsystem.
\[
\mbf G_{\mathrm o}\;:\;\qquad D_{zr} =1 \qquad D_{vw} = 1.
\]
Examining the partial derivatives and boundary values used in Eqs.~\eqref{eqn:example2_PDE} and~\eqref{eqn:example2_BC}, we first define the continuity equation using $n_0=1$ so that $\hat{\mbf x}_0=\zeta_2$ and $n_2=1$ so that $\hat{\mbf x}_2=\zeta_1$.
\[
n\;:\; \qquad n=\{1,0,1\}\qquad N=2.
\]
For this definition of $n$, we have $n_{\hat{\mbf x}}=n_{S_0}=2$ and $n_{S_1}=n_{S_2}=1$ -- there are two $0^{th}$ order and one $1^{st}$ and $2^{nd}$ order partial derivatives. In addition, $n_{S}=2$, indicating there are 2 possible partial derivatives. Thus
\[
S^0\hat{\mbf x}=\bmat{\hat{\mbf x}_1\\ \hat{\mbf x}_2}=\bmat{\zeta_1\\ \zeta_2} \qquad S\hat{\mbf x}=S^2\hat{\mbf x}=\hat{\mbf x}_2=\zeta_1.
\]
Next, we construct $(\mcl B \hat{\mbf x})$ -- the vector of all possible boundary values of $\hat{\mbf x}$ allowable for the given $n$.
\[		
(\mcl B \hat{\mbf x})=\bmat{(\mcl C\hat{\mbf x})(0) \\ (\mcl C\hat{\mbf x})(1)}=\bmat{\hat{\mbf x}_2(0)\\ \hat{\mbf x}_{2,s}(0)\\\hat{\mbf x}_2(1)\\ \hat{\mbf x}_{2,s}(1)}=\bmat{\zeta(0)\\ \zeta_{1,s}(0)\\\zeta(1)\\ \zeta_{1,s}(1)}
\]
We may now define the BCs. There is no ODE state, however, there is a disturbance $w$ that influences the PDE via the signal $v$ which can be chosen as $v(t)=w(t)$. Then, the BCs appear in the form
{\small
\begin{align*}		
&\bmat{0\\v(t)}= \int_{a}^{b} B_{I}(s)\mcl F\hat{\mbf x}(t,s)ds -B(\mcl B \hat{\mbf x})(t)\\
&=\int_{0}^{1} B_{I}(s)\bmat{\hat{\mbf x}_1(t,s)\\\hat{\mbf x}_2(t,s)\\\hat{\mbf x}_{2,s}(t,s)\\\hat{\mbf x}_{2,ss}(t,s)}ds -B(\mcl B \hat{\mbf x})(t)=\int_{0}^{1} B_{I}(s)\bmat{\zeta_1(t,s)\\\zeta_2(t,s)\\\zeta_{1,s}(t,s)\\\zeta_{1,ss}(t,s)}ds -B\bmat{\zeta_1(t,0)\\ \zeta_{1,s}(t,0)\\\zeta_1(t,1)\\ \zeta_{1,s}(t,1)}
\end{align*}}
By inspection of Eq.~\eqref{eqn:example2_BC}, we may now define the parameters for $\mbf G_{\mathrm b}$ and hence $X_{v(t)}$ as
\[
\mbf G_{\mathrm b}\;:\;\qquad B =-\bmat{1&-k&0&0\\0&0&1&l} \qquad B_v = -\bmat{0\\1}.
\]

The final step is to define the parameters of the PDE dynamics. Ignoring the integral terms for simplicity, and noting that $\hat{\mbf x}=\bmat{\hat{\mbf x}_0\\\hat{\mbf x}_2}=\bmat{\zeta_2\\\zeta_1}$, $S\hat{\mbf x}=S^2\hat{\mbf x}=\hat{\mbf x}_2=\zeta_1$ and $r=v=\emptyset$, we have
\begin{align*}
\bmat{\dot{\zeta}_2(t,s)\\\dot{\zeta}_1(t,s)} &=  A_{0}(s)\bmat{\zeta_2(t,s)\\\zeta_1(t,s)\\\zeta_{1,s}(t,s)\\\zeta_{1,ss}(t,s)}+B_{xb}(s) \bmat{\zeta_1(t,0)\\ \zeta_{1,s}(t,0)\\\zeta_1(t,1)\\ \zeta_{1,s}(t,1)}.
\end{align*}
By inspection of Eq.~\eqref{eqn:example2_PDE}, we have two non-zero parameters in $\mbf G_{\mathrm p}$. However, the interconnection signal $r$ has an integral term, which can be written as
\begin{align*}
r(t) &= \int_{a}^{b} C_{r}(s)\mcl F\hat{\mbf x}(t,s)ds +D_{rb}(\mcl B \hat{\mbf x})(t)\\
&=\int_{0}^{1} C_{r}(s)\bmat{\hat{\mbf x}_1(t,s)\\\hat{\mbf x}_2(t,s)\\\hat{\mbf x}_{2,s}(t,s)\\\hat{\mbf x}_{2,ss}(t,s)}ds +D_{rb} (\mcl B \hat{\mbf x})(t)=\int_{0}^{1} C_{r}(s)\bmat{\zeta_1(t,s)\\\zeta_2(t,s)\\\zeta_{1,s}(t,s)\\\zeta_{1,ss}(t,s)}ds +D_{rb}\bmat{\zeta_1(t,0)\\ \zeta_{1,s}(t,0)\\\zeta_1(t,1)\\ \zeta_{1,s}(t,1)}.
\end{align*}
Clearly, only $C_{r,0} = \bmat{1&0}$ is non-zero, whereas the remaining terms are zero which gives us the final set of parameters for the PDE subsystem as
\[	
\mbf G_{\mathrm p}\;:\;\qquad A_{00} = \bmat{0&0\\1&0}, \qquad A_{02} = \bmat{c^2\\0},\qquad C_{r,0} = \bmat{1&0}.
\]
This completes the definition of the GPDE.
\end{illus}

\begin{illus}{Chemical Reactor with Cooling Jacket}
Consider a chemical reactor with a cooling jacket as described in~\cite{Karafyllis2019}. In this model the reactor temperature is a lumped parameter system while the coolant temperature is a distributed state that varies along the length of the pipe. Assuming conduction inside the cooling jacket to be negligible, we obtain the following coupled ODE-PDE.
\begin{align}
\dot{x}(t) &= kx(t) + \mu\int_{0}^{1}\mbf{x}(t,s)ds\notag\\
\dot{\mbf{x}}(t,s) &= -c\partial_s\mbf{x}(t,s)- \zeta\mbf{x}(t,s)+\zeta x(t)\qquad
\mbf{x}(t,0)=w(t)
\end{align}
where $x$ is the reactor temperature, $\mbf{x}$ is the temperature in the cooling jacket, $w(t)$ is a disturbance, $\mu,c,\zeta$ are positive constants, and $k$ is a negative constant.
In this model, the distributed state $\mbf{x}$ has a single boundary condition and highest spatial derivative of order 1, so $n=\{0,1\}$. In order to retain the parameter dependencies, we use the formulae in~\Cref{fig:Gb_definitions,fig:PIE_subsystem_equation} to obtain the following PIE representation.
\begin{align}
\dot{x}(t) &= kx(t) + \int_{0}^1\mu(1-s)\ubar{\hat{\mbf x}}(t,s)ds,\notag\\
\int_{0}^ss\ubar{\dot{\hat{\mbf x}}}(t,\theta)d\theta &= \zeta x(t)-c\ubar{\hat{\mbf x}}(t,s)-\int_{0}^s\zeta \ubar{\hat{\mbf x}}(t,\theta)d\theta-\dot w(t)-\zeta w(t),
\end{align}
or, alternatively,
\begin{align}
\fourpiFull{1}{0}{0}{0}{s}{0}\bmat{\dot x(t)\\\ubar{\dot{\hat{\mbf x}}}(t)}+\fourpi{0}{0}{1}{0}\dot{w}(t)=\fourpiFull{k}{\mu (1-s)}{\zeta}{-c}{-\zeta}{0}\bmat{x(t)\\\ubar{\hat{\mbf x}}(t)}+ \fourpi{0}{0}{-\zeta}{0}w(t),
\end{align}
where $\ubar{\hat{\mbf x}} = \partial_s\mbf{x}$.

\end{illus}

\end{document}